\documentclass[preprint,12pt]{elsarticle}

\usepackage{amssymb}
 \usepackage{amsthm}

\usepackage{amsmath}

\usepackage[colorinlistoftodos]{todonotes}

\usepackage{mathabx}
\usepackage{graphicx}
\usepackage{subfig}
\usepackage{tabularx}

\usepackage{mathrsfs}     
\usepackage{helvet}         
\usepackage{courier}        
\usepackage{type1cm}      
\usepackage{xcolor}
\usepackage{lipsum} 
\usepackage{url}

\usepackage{geometry,calc,color}
\usepackage{amsmath,amssymb}
\usepackage{enumerate}
\usepackage{arydshln}
\usepackage{siunitx}
\usepackage{booktabs}
\usepackage{multirow}
\usepackage{tikz}

\usepackage{cleveref}
\crefrangelabelformat{equation}{(#3#1#4--#5\crefstripprefix{#1}{#2}#6)}
\crefrangelabelformat{subequation}{(#3#1#4--#5\crefstripprefix{#1}{#2}#6)}

\usepackage{pgfplots}
\pgfplotsset{compat=newest}

\usepackage{textcomp}

\allowdisplaybreaks[4]

\renewcommand{\epsilon}{\varepsilon}
\renewcommand{\phi}{\varphi}
\renewcommand{\theta}{\vartheta}
\newcommand{\p}{k} 

\newcommand{\el}{K}

\newcommand{\Omh}{\Omega_h}

\newcommand{\PinablaK}{\Pi^ {\nabla\!,k}_K}

\newcommand{\K}{K}

\newcommand{\hK}{\h_\K}

\newcommand{\h}{H}
\newcommand{\mh}{h}
\newcommand{\microTess}{\mathcal{T}_\mh}
\newcommand{\microTh}{\mathcal{T}_\mh}

\newcommand{\Poly}[1]{\mathbb{P}_{#1}}
\newcommand{\PolyQ}[1]{\mathbb{Q}_{#1}}

\newcommand{\proj}[1]{\pi_{#1}}

\newcommand{\edge}{e}

\newcommand{\hu}{\widehat u}

\newcommand{\Th}{\mathcal{T}_H}

\newcommand{\Tess}{\Th}

\newcommand{\EK}{\mathcal{E}^\el}

\newcommand{\x}{{\mathbf x}}

\newcommand{\TriK}{\widetilde{\mathcal{T}}_K}
\newcommand{\tVK}{W^K}

\newcommand{\tuI}{\widetilde u_I}

\newcommand{\Pinabla}{\Pi^{\nabla}}

\newcommand{\Eb}{{\mathcal{E}^\partial }}
		
\newcommand{\dn}{\partial_\nu}
\newcommand{\Svem}{S_a^K}

\newcommand{\Vh}{V_\h}
 
\newcommand{\VKp}{\VEMK}

\newcommand{\Pn}{\Pi^\nabla\!}

\newcommand{\xK}{\mathbf{x}_\el}

\newcommand{\PolyDom}{{\Omega}}

\newcommand{\pwPoly}{\Poly{k}(\Tess) }

\newcommand{\Harm}{\mathcal{H}}

\newcommand{\cstar}{c_*}
\newcommand{\Cstar}{C^*}

\newcommand{\Ah}{\mathcal{A}_{H}}

\newcommand{\Pnphi}{\Pi^\nabla\! (\phi)}
\newcommand{\Pnpsi}{\Pi^\nabla\!(\psi)}

\newcommand{\Oh}{\Omega_{\mh}}

\newcommand{\nuhx}{\nu_h}

\newcommand{\dnh}{\partial_{\nu_h}}
\newcommand{\tg}{g^\star}

\newcommand{\sx}{\sigma}

\newcommand{\nK}{\nu_K}

\newcommand{\VEM}{V_\h}
\newcommand{\VEMK}{V^{K,\p}}
\newcommand{\VEMo}{\mathring V_\h}

\newcommand{\ah}{a_\h}

\newcommand{\mE}{E}

\newcommand{\VE}{V_E}

\newcommand{\cv}{\widecheck v}
\newcommand{\cVEM}{\widecheck V_\h }
\newcommand{\hv}{\widehat v}
\newcommand{\hVEM}{\widehat V_\h }

\newcommand{\cu}{\widecheck u_h}

\renewcommand{\hu}{\widehat u_h}
\newcommand{\stab}{S}

\newcommand{\uh}{u_h}
\newcommand{\bK}{{\partial K}}
\newcommand{\bOh}{{\partial\Oh}}

\newcommand{\Cuno}{C_1}

\newcommand{\Cdue}{C_2}
\newcommand{\Ctre}{C_3}

\newcommand{\dsh}{\partial_\sigma}

\newcommand{\stabweight}{\beta}
\newcommand{\uI}{u_I}

\newcommand{\bO}{\partial\O}
\renewcommand{\O}{\Omega}

\newcommand{\interval}[2]{(#1,#2)}

\newcommand{\Ekv}{\EK_\text{ver}}
\newcommand{\Ekh}{\EK_\text{hor}}

\newcommand{\SBM}{SBM}
\newcommand{\BDT}{BPP}
\newcommand{\BHL}{BHL}

\newcommand{\kstar}{k}
\newcommand{\khat}{{\widehat k}}

\newcommand{\Corrstar}[1]{\mathscr C \left[#1\right]}
\newcommand{\Corrhat}[1]{\widehat{\mathscr C} \left[#1\right]}
\newcommand{\dist}{\mathrm{d}}

\newcommand{\Franke}{u_\mathrm{Frnk}}

\pgfplotsset{compat=1.9}

\usepackage{textcomp}
\renewcommand{\epsilon}{\varepsilon}

\newcommand\Tstrut{\rule{0pt}{2.6ex}}         

\newtheorem{theorem}{Theorem}[section]

\newtheorem{lemma}[theorem]{Lemma}

\theoremstyle{definition}

\newtheorem{assumption}[theorem]{Assumption}
\newtheorem{comment}[theorem]{Comment}
\theoremstyle{remark}
\newtheorem{remark}[theorem]{Remark}

\begin{document}
	
	\begin{frontmatter}
%
%
		\title{The virtual element method  \\ on polygonal pixel--based tessellations\tnoteref{proj}}
		\tnotetext[proj]{This paper has been realized in the framework of the ERC Project CHANGE, under the EU’s Horizon 2020  programme (grant agreement No 694515). It was co-funded by the MIUR Progetti di Ricerca di Rilevante Interesse Nazionale (PRIN) Bando 2017 (grant 201744KLJL) 
		 and Bando 2020 (grant 20204LN5N5).  M. Pennacchio has been partially supported by ICSC—Centro Nazionale di Ricerca in High Performance Computing, Big Data, and Quantum Computing funded by European Union—NextGenerationEU. The authors are members of  INDAM-GNCS.}
		\author[imati]{S. Bertoluzza} 
			\ead{silvia.bertoluzza@cnr.it}
			 	 \author[unipv,imati]{M. Montardini}
			\ead{monica.montardini@unipv.it}
		 \author[imati]{M. Pennacchio}
		 	\ead{micol.pennacchio@imati.cnr.it} \author[imati]{D. Prada\corref{cor1}}
			\ead{daniele.prada@imati.cnr.it}
		\affiliation[imati]
		{organization={IMATI ``E. Magenes'', CNR},
			addressline={via Ferrata, 5A},
			city={Pavia},
			postcode={27100},
			state={PV},
			country={Italy}}

		\affiliation[unipv]
		{
			organization={Dipartimento di Matematica, Universita di Pavia},
			addressline={via Ferrata, 5},
			city={Pavia},
			postcode={27100},
			state={PV},
			country={Italy}}
%
%
\begin{abstract}
We analyze and validate the virtual element method combined with a boundary correction similar to the one in \cite{SBMho,VEM_weakly}, to solve problems on two dimensional domains with curved boundaries approximated by polygonal domains. We focus on the case of approximating domains obtained as the union of squared elements out of a uniform structured mesh, such as the one that naturally arises when the domain is issued from an image. We show, both theoretically and numerically, that resorting to  polygonal elements allows the assumptions required for stability to be satisfied for any polynomial order. This  allows us to fully exploit the potential of higher order  methods. Efficiency is ensured by a novel static condensation strategy acting on the edges of the decomposition.
\end{abstract}
%
%
%
\end{frontmatter}

\section{Introduction}
The simplest (and cheapest!) meshes that can be used to approximate a complex domain are the ones whose elements coincide with elements of a sufficiently fine squared/cubic uniform structured grid. This holds particularly true when the domain is retrieved from the results of some imaging procedure, as it often happens, for instance, in medical applications. In such a framework (an approximation of) the physical domain is already given as the union of pixels/voxels that a segmentation procedure tags as belonging to the physical object of interest. These can be seen as elements of a very fine structured quadrangular/hexahedral mesh. The polygonal/polyhedral domain obtained as their union can then be used as an approximation of the physical domain, to be used in the numerical solution of a PDE, modeling some physical behavior.

It is a well known fact that approximating the solution of a problem in a physical domain by simply solving, by a finite element method, a problem in an approximated polygonal domain, with a boundary condition  somehow ``copied'' from the physical boundary data, yields, for methods of higher order, a suboptimal result. For  homogeneous Dirichlet boundary value problems,
the fact itself of approximating the physical domain (with curved boundary)   introduces an error that, even in the best of cases, can be of the order $\delta^{3/2}$, $\delta$ being the distance between physical and approximate boundary \cite{strang1973change}. In and of itself this  leads to suboptimality whenever the order $k$ of the method is greater than or equal to two. In 
 the framework we are considering the situation is worse, and the method turns out to be suboptimal also for $k=1$, with a convergence of order only $h^{1/2}$, $h$ denoting the meshsize, as observed, both theoretically and numerically, in \cite{ramiere2008convergence}. Nevertheless, such an approach is currently used by many practitioners. Resorting to a so called {\em  microFEM} approach (\cite{microFEM1}), they use the mesh whose elements coincide with the voxels in a microCT scan, in order to simulate  some physical phenomenon taking place in an underlying (unknown) domain.
 Of course, given the extreme fineness of the mesh, the results obtained by such an approach turn out to be sufficiently accurate. The  cost  is however much higher, when compared to the cost of the finite element method of the same accuracy on polygonal physical domains. 

Different options exist to counter the sources of error related to the approximation of the domain, and thus obtain a more efficient method, provided, of course, we can rely on information on the actual physical domain.   By what means, and how accurately, such information can be retrieved from available imaging data is a crucial question that is, however, out of the scope of this paper (we refer to \cite{EdgeDetectionSurvey22} and the references therein for an up to date survey on edge detection methods that can be used to this aim). Once this information is available, one option is to work on a possibly coarser mesh  and either use a fictitious domain approach, as in the finite cell method \cite{parvizian2007finite}, or state the problem on the actual domain by ``cutting'' the elements that cross the boundary, as done in the cutFEM method \cite{burman2015cutfem}. Remark that preprocessing the image by changing its resolution is an easy way to obtain a coarser versions of the domain approximation. These approaches can also  be combined with different discretization methodologies such as isogeometric analysis \cite{Schillinger}.

A different approach consists in resorting, while working on the approximate domain, to techniques specifically designed to take into account the fact that its boundary  does not coincide with the actual curve/surface where  boundary conditions should be prescribed.
The first example of this strategy was introduced already in the early '70s in the seminal paper by Bramble, Dupont and Thom\'ee \cite{BDT}. There, for convex domains in 2D, a Taylor extrapolation along the direction normal to the boundary of the approximating polygonal domain was leveraged, within a Nitsche approach, to weakly impose the correct boundary conditions. Introduced in the late 2010s, and already tested in different application fields (see, e.g., \cite{SBMmechanics,SBMfreesurface,SBMhyperbolic,SBMreduced}), the {\em shifted boundary method} (SBM, see \cite{SBM1, SBM2})  overcomes many of the limitations of the original Bramble, Dupont and Thom\'ee approach (BDT), by a careful  choice of the extrapolation direction and a clever design of the stabilization term for the underlying  Nitsche method. Initially limited to order up to two, an high order version of the SBM, allowing, in principle, to attain any arbitrary order $k$, has been recently introduced and analyzed in  \cite{SBMho}. 
In \cite{burman2020dirichlet} a variant of the BDT method is proposed, based on the Lagrange multiplier method for Dirichlet boundary conditions \cite{babuska}. In the same paper, the relation with  Nitsche's method with boundary correction is also discussed,  in the spirit of \cite{stenberg1995some}. Nitsche's method with boundary correction is also studied in the framework of the cutFEM method in \cite{burman2018cut}, where the analysis is carried out for an a priori arbitrary choice of the extrapolation direction. Such direction is, however, in practice, chosen to be, as in SBM, normal to the boundary of the physical domain. The {\em polynomial extension finite element method} (\cite{cheung2019optimally}) avoids the problem of choosing an extrapolation direction by replacing Taylor extrapolation with an  {\em averaged Taylor polynomial extrapolation} to also attain arbitrary order under suitable conditions. Similar ideas can be also leveraged in the discontinuous Galerkin framework (\cite{cockburn2012solving}).

Unfortunately, as the order $k$ of the method gets higher, the approximate domains that we are considering eventually fail,  when discretized by finite elements, to satisfy one of the assumptions required for the stability of such methods. Indeed, the analysis of boundary correction methods requires that the ratio between the distance from the approximate to the true boundary and the diameter of the elements is lower than a constant that decreases to zero as $k$ increases.
A way to overcome this limitation is to replace the fine mesh with a sufficiently coarser polygonal/polyhedral mesh, obtained by agglomeration.  This allows to make the diameter of the elements larger, while keeping the distance of the two boundaries fixed. On the new polygonal/polyhedral mesh, one can then combine a boundary correction strategy with one of the discretization methods capable of handling polytopal meshes. In this direction, in \cite{liu2022weak} the authors combine a polynomial extension method with a {\em weak Galerkin} finite element method. In \cite{VEM_curvo,VEM_weakly} a boundary correction method similar to the one in \cite{burman2018cut} is combined, { both in 2 and 3D}, with  {\em virtual elements} (VEM, \cite{basicVEM,hitchVEM,3DVEM}). 
Thanks to  its flexibility, its robustness (in particular with respect to the shape of the elements), and its potential for high accuracy (the discretization can be designed to be of arbitrarily high order, \cite{3DVEM2}), virtual elements
have, since their introduction in the early 2010's, rapidly gathered the attention of the scientific community. This resulted  in extensions to deal with different type of equations (\cite{VEM_second_order,beirao_parab,navier-stokes2d,stokes3d,Antonietti_VEM_Stokes,beirao_elastic,beirao_linear_elasticity,VEM_3D_elasticity,Antonietti_VEM_Cahn,perugia_Helmholtz,ABPV:minsurf,VEM_Laplace_Beltrami,Moraetal15,wriggers2017efficient,brenner2021ac}), 
different formulations (\cite{VEM_mixed,de2016nonconforming,BMPP_VEM_nonconforming_stab,MASCOTTO_Trefftz,curved_Trefftz}) and numerous applications   (\cite{VEM_discrete_fracture,BEIRAODAVEIGA2015327,CHI2017148,wriggers2017efficient,
	wriggers2016virtual,BEIRAODAVEIGA2021}).

 The aim of this paper is to propose and validate, initially in two dimensions, the use of the combination of virtual elements and boundary correction methods  such as BDT or SBM, as a way to obtain an efficient solver in the context presented previously. To this aim we will need to adapt the theoretical results obtained in \cite{VEM_weakly} to the present framework. Indeed, tessellations whose elements are obtained as the union of squares of a structured uniform mesh,  do not satisfy the standard shape regularity assumptions under which virtual element methods are usually analyzed.  We will also need to show how to efficiently handle elements with a large number of small edges, deriving from the agglomeration process.

The paper is organized as follows. In Section \ref{sec:notation} we introduce a unified notation for a class of boundary correction methods that includes, among others, the BDT method and some of its variants, as well as the high order SBM.  In Section \ref{sec:VEM} we recall definition and properties of the plain virtual element method, and we present the virtual element method for curved domains  \cite{VEM_weakly} obtained by combining the VEM with the previously introduced boundary correction methods. In Section \ref{sec:Lazy} we present a static condensation procedure allowing to keep the  resulting  linear system small, also  in the presence of elements with a large number of small edges. In Section \ref{sec:expes} we present the numerical results attesting to the validity of our proposal and comparing the performance attainable with virtual elements with those obtained by the same boundary correction approach in the finite element framework. 
The proofs of several theoretical bounds and estimates, needed to extend known virtual element results to the present framework, are presented in three appendixes.

\newcommand{\umh}{u_\mh}
\newcommand{\nh}{\nu_\mh}
\newcommand{\map}{B_\mh}

\section{Boundary correction methods in the finite element framework} \label{sec:notation}

For the sake of simplicity, we focus on a simple model problem, namely on the Poisson equation 
\begin{equation}\label{prob_mod}
	-\Delta u = f, \text{ in }\PolyDom, \qquad u = g \text{ on }\Gamma = \partial \PolyDom,
\end{equation}
where $f \in L^2(\Omega)$, $g \in H^{1/2}(\partial \Omega)$, and where $\Omega \subseteq \mathbb{R}^2$ is a  domain with a curved boundary $\Gamma = \partial\Omega$, assumed, for the sake of convenience, to be of class $C^\infty$. We consider here a context where the domain is not directly available but is obtained as the result of an imaging process. We will then only have an approximation of $\Omega$, which we will denote by $\Oh$, naturally decomposed as the union of (tiny) squared elements (the pixels) of size $h \times h$, with edges parallel to the axes. We will not address here the segmentation problem that needs to be solved in order to single the approximated domain $\Oh$ out of the image.  It is natural, in this context, to assume that for all $\x \in \bOh$, the distance $\dist(\x,\bO)$ of $\x$ to $\bO$ verifies $\dist(\x,\bO) \lesssim h$ (this is another way to say that $\Oh$ is an approximation to $\O$). This happens for instance if $\Oh$ is constructed  by retaining all the squares contained in $\O$.

\  

\newcommand{\Qk}{Q$_k$}

The simplest  method to solve \eqref{prob_mod} (and one of the most used by practitioners) is the finite element method on the quadrilateral mesh $\microTess$, whose elements are the pixels, in which $\Oh$ is naturally split. The approximation space  is  chosen as the standard \Qk\ finite element space
\[
V_\mh = \{ v \in C^0(\Oh): \ v|_K \in \PolyQ{k}\ \forall K \in \microTess \},
\]
$\PolyQ{k}$ denoting the space of polynomials $p(x,y)$ of degree at most $k$ both in $x$ and in $y$. Using Nitsche's method (see \cite{Nitsche}) to impose the boundary conditions, the approximate solution is defined as the element $\umh \in V_\mh$ such that for all $v \in V_\mh$ it holds
\begin{multline}\label{NitscheFEM}
\int_{\Oh} \nabla \umh \cdot \nabla v - \int_{\partial \Oh}  \dnh \umh\, v 
- \int_{\partial\Oh}   \umh \left (\dnh v  - \gamma \mh^{-1} v \right)\\
 = \int_{\Oh }f v - \int_{\partial\Oh} \tg ( \dnh v - \gamma 
\mh^{-1} v),
\end{multline}
with $\tg$ suitably defined, where \(\nh\) is the outer unit normal to \(\partial \Oh\), $\gamma > 0$ is a suitable, sufficiently large, scalar constant, and where $\dnh v \in L^2(\bOh)$ denotes the $L^2(\bOh)$ function coinciding with $\nabla u_h \cdot \nu_h$ on each boundary edge of the mesh $\microTess$. 
It is however  known that, in our framework,  such a method is suboptimal already for $k=1$ (see \cite{ramiere2008convergence}). Several strategies have been proposed in order to retrieve the optimal order of approximation for the finite element solution on approximating polygonal domains, relying on the 
 idea of suitably correcting the Nitsche's formulation \eqref{NitscheFEM} so that the boundary condition is somehow imposed on the original boundary, while still maintaining all computations on the approximating domain. 
 
 \
 
\begin{figure}
	\centering
	\includegraphics[height=4cm]{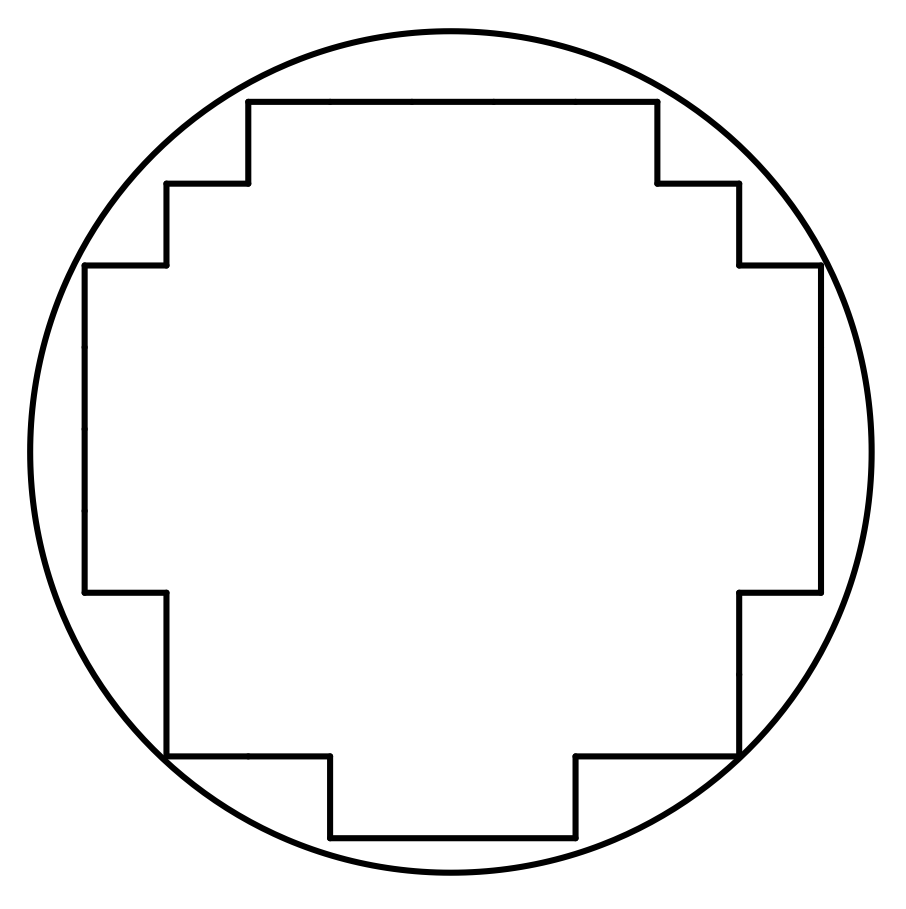} \qquad
	\includegraphics[height = 3.5cm]{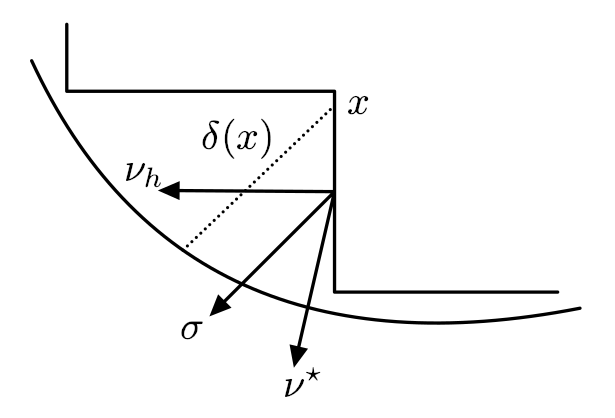}
	\caption{ An approximate domains $\Omega_h$ falling in  our framework. The theoretical framework does not in principle require the extrapolation direction $\sigma$ to coincide with either the normal $\nu_h$ to the approximate boundary or the normal $\nu^\star$ to the physical boundary, though the latter is generally the best choice.
	}\label{fig:sigma}
\end{figure}

 We consider here a general formulation that encompasses a number of such strategies.
%
 %
We let
 $\Eb$ denote the set of edges of $\microTh$ lying on $\partial\Oh$, and, for $\x \in \partial\Oh$,  we let $\sigma(\x)$ denote an outward unit vector. 
  Assuming that $\Oh \subseteq \O$, for $\x \in \partial\Oh$ we let $\delta(\x) > 0$ denote the distance to $\partial \Omega$ along the direction $\sigma$, that is, the smallest non negative scalar such that 
 \[
 \x + \delta(\x) \sigma(\x) \in \partial \Omega,
 \]
see Figure \ref{fig:sigma}. Letting $\dsh^j u = ({\partial}_{\sigma})^j u$ denote the $L^2(\partial\Oh)$ function coinciding, on each edge, with the $j$-th partial derivative of $u$ in the $\sigma$ direction, and letting, for $\x \in \partial\Oh$, $\tg(\x)$ be defined as \[\tg (\x) = g(\x+\delta(\x) \sigma(\x)),\]
we look for $\uh \in V_h$ such that for all $v \in V_h$ it holds that 
\begin{multline}\label{bilinear_Nitsche}
 \int_{\Oh} \nabla \uh\cdot \nabla v - \int_{\partial \Oh} \dnh \uh   v
	 -   \sum_{e\in \Eb }
	\int_e 
	\left(\uh+ \Corrstar {u_h}\right)
	\left( \dnh v  - \gamma h^{-1} (v + \Corrhat v ) \right)  \\= \int_{\Oh} f v -
	\sum_{e \in \Eb} \int_e \tg  	\left( \dnh v  - \gamma h^{-1} (v + \Corrhat v  )\right),
\end{multline}
where the ``correction'' term for the trial and test functions are defined as
\begin{equation}
\label{definecorrections}
\Corrstar w = \sum_{j=1}^{\kstar} \frac{\delta^j}{j!} \dsh^j w,\qquad  
\Corrhat w = \sum_{j=1}^{\khat} \frac{\delta^j}{j!} \dsh^j w. 
\end{equation}
Different choices have been proposed in the literature for the extrapolation direction $\sigma$ and for the parameter $\khat$ involved in the definition of $\Corrhat{\cdot}$ in \eqref{definecorrections}, resulting in different correction methods. The choice $\khat = 0$ (which is to be interpreted as $\Corrhat v = 0$) and $\sigma = \nh$ yields the method originally proposed in the seminal work \cite{BDT} by Bramble, Dupont and Thom\'ee.
Choosing $\khat = 1$ and $\sigma(\x) = \nabla \dist(\x,\bO) = \nu^\star$, (that is  $\sigma(\x) = \nu(\x + \delta(\x)\sigma(\x))$, $\nu$ denoting the outer unit normal to $\partial \Omega$) yields the high order shifted boundary method (\SBM) proposed in \cite{SBMho}.
The choice $\khat = 0$ and $\sigma$ a priori arbitrary is analyzed in \cite{burman2015cutfem} and in \cite{Nitsche,VEM_curvo}, where it is respectively exploited in the context of the cut finite element method, and  in the virtual element framework. In both cases $\sigma$ is in practice also chosen as the gradient of the distance function to the boundary, so that $\delta(\x)$ is as small as possible.


\

\begin{figure}
	\centering
	\includegraphics[width=13cm]{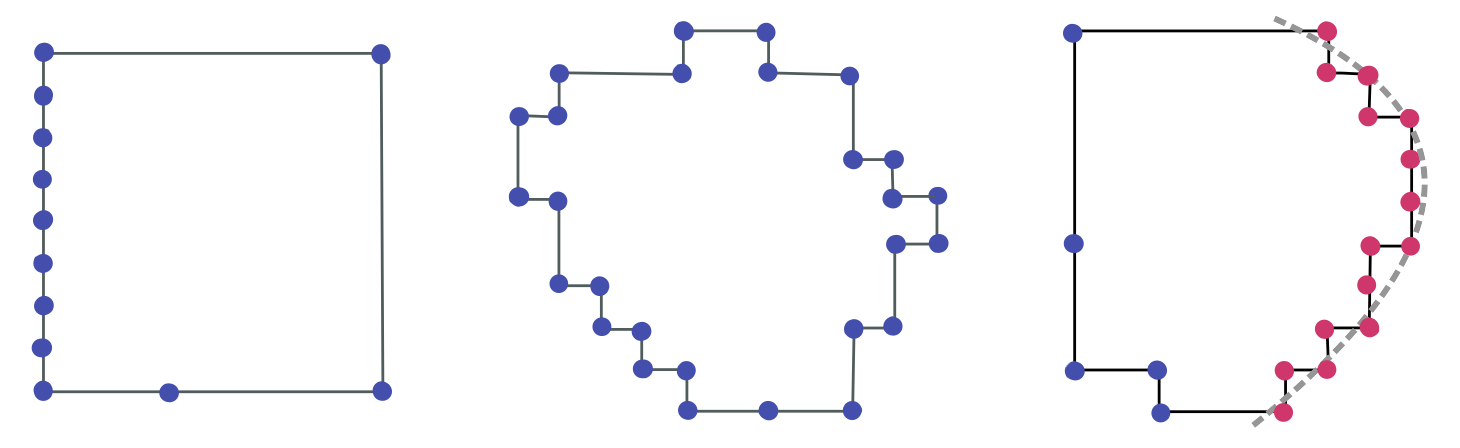}
	\caption{Three possible elements of the tessellation $\Tess$. For the sake of the exposition, boundary edges (with vertices marked in red) are never agglomerated to form larger edges, even when this is possible. Agglomeration of interior edges (vertices in blue) into larger edges fits instead in our exposition.} \label{fig:elements}
\end{figure}

\begin{figure}
	\centering
	\includegraphics[width=13cm]{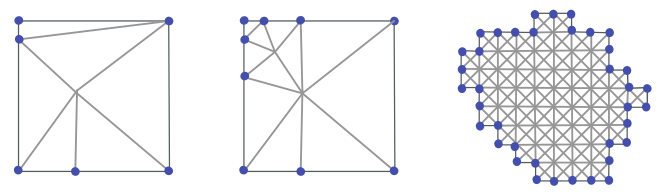}
	\caption{
		Three examples of the auxiliary triangulation $\TriK$. As one can see in the leftmost example, the presence of two adjacent boundary edges with very different length results in a badly shaped triangle. Adding few nodes, as in the central example, may improve the shape regularity of the triangulation. As the  $\TriK$ is allowed to have a number of elements as large as needed, the presence of a large number of very small edges does not, in itself, result in badly shaped triangulation (see the rightmost example).	
	}	\label{fig:additional}
\end{figure}

 For all these choices it is possible to prove that, under the condition that $\delta_h = \max_{x\in \partial\Oh} |\delta(x)| < \tau h$, with $\tau > 0$ a constant depending on the order $k$, and provided $\gamma$ is large enough, the method is stable and converges with optimal order. The (small) constant $\tau$  decreases to $0$ as $k$ increases. 
 { If we construct the domain and choose the  extrapolation direction $\sigma$,  in such a way that $\delta_h = o(h)$,} for all these choices and for all order $k$ there exists a $h_0$ such that, provided $\gamma$ is large enough, for all $h< h_0$ the method is stable and converges with optimal order.

%

{ Unfortunately for the class of polygons that we are considering in our framework it is not possible to choose $\sigma$ in such a way that  condition $\delta_h = o(h)$ is satisfied.} Indeed, in general, for approximating domains issued from imaging, which are the union of equal square elements (pixels), we have that\footnote{\label{note1}Unless the specific value of the constant $C$ is explicitely  needed, throughout the paper we will write $A \lesssim B$ (resp. $A \gtrsim B$) to indicate that the quantity $A$ is bounded from above (resp. from below) by a constant $C$ times the quantity $B$, with $C$ possibly depending on $\Omega$ as well as on the parameters $\alpha_0$ and $N_0$ appearing in the shape regularity assumption \ref{shape_regular}, but otherwise independent of the shape and size of the elements of the tessellations. The notation $A \simeq B$ will stand for $A \lesssim B \lesssim A$.
} $\delta_h \simeq h$, { even for the best choice of $\sigma$}. Consequently, there exist a $\bar k$ such that for $k > \bar k$ our tessellations will not satisfy the condition $\delta_h \leq \tau h$, and, as $k$ increases all the considered methods will eventually lose stability.
\label{page1}


 \begin{remark}
	Other boundary correction strategies can be found in the literature that we could include in the unified formulation \eqref{bilinear_Nitsche}, provided we allow more general forms for the correction terms $\mathscr{C}$ and $\widehat{\mathscr{C}}$, than the ones in \eqref{definecorrections}. We recall the polynomial extension method (\cite{cheung2019optimally}, see also \cite{cockburn2012solving}).
\end{remark}

\section{The Virtual Element Discretization}\label{sec:VEM}

The main idea behind the method we are proposing is to discretize the polygonal approximate domain $\Oh$ with a polygonal tessellation $\Th$, with meshsize $\h$, whose elements are obtained as union of quadrilateral elements of the fine tessellation $\microTess$. On the tessellation $\Th$ we can then use the virtual element Nitsche's method with boundary correction proposed in \cite{VEM_curvo, VEM_weakly}.  By taking particular care in the implementation, this will result in a method with a much more favorable cost/performance ratio than the one obtained by plain finite elements as used in the  microFEM approach, without the need for modifying the bulk bilinear form near the boundary. To this aim, we start by reviewing the definition of the method we will be employing, as proposed in \cite{VEM_weakly}. 
 
\subsection{The tessellation} 
%
We assume that the tessellation $\Tess$ of  $\Omh$ into polygons,  obtained by agglomeration of the square elements of $\microTess$, satisfies the following Assumption.  
%
%
%

\begin{assumption}\label{shape_regular} All elements $K \in \Tess$ are simply connected union of squares of the cartesian mesh $\microTess$ of meshsize $h$.
Moreover, letting for $K \in \Tess$
\[\hK = \max_{(x,y),(x',y') \in K} \max\{| x-x' |,|y - y'|\},\] we have $\hK \simeq H$,{\ with $H=\max_{K\in \Tess} \hK$}, and there exist constants $\alpha_0<1$ and $N_0\geq 2$ such that:
	\begin{enumerate}[(i)]
		\item 
		each element $K \in \Th$ verifies $S(\xK,\alpha_0 \h _K) \subset K \subset S(\xK,\h_K)$, $\xK = (x_K,y_K)$ denoting the geometrical center of $K$ and $S(\xK,d)$ the square of center $\xK$ and side $d$;
	\item for all $x\in (x_K - \hK/2,x_K+\hK/2)$ (resp. $y \in (y_K - \hK/2,y_K+\hK/2)$) there exist at most $N_0$ values $y\in \mathbb{R}$ (resp. $x\in \mathbb{R}$) such that $(x,y) \in \bK_\text{hor}$ (resp. $(x,y)\in \bK_\text{ver}$), $\bK_\text{hor}$ and $\bK_\text{ver}$ respectively denoting the union of horizontal and vertical edges of $K$.
	\end{enumerate}
\end{assumption}

In order to prove the virtual element approximation estimate,  we will also need to make the following additional assumption, where
	the shape regularity constant of a triangulation is intended as the maximum over all triangles $T$ of the ratio $h_T/\rho_T$, $h_T$ and $\rho_T$ respectively denoting the diameter of the circumscribed and inscribed circle. 
	\begin{assumption}\label{additional}
		There exists a constant $\alpha_1$ such that for each element $K \in \Th$ there exists a conforming triangulation $\TriK$ of $K$ with shape regularity constant at most $\alpha_1$, whose set of boundary edges  coincides with the set of edges of $K$.
	\end{assumption}
	In the above assumption we do not require the number of elements of $\TriK$ to be uniformly bounded. We underline that the triangulation $\TriK$ only plays a role in the theoretical analysis of the method, and its actual construction 
is not needed for the design and implementation of the method. 
{  Remark that the virtual element method allows to handle  polygons with aligned edges. Then, whenever a geometrical straight edge of the element $K$ is the union of two or more edges of the fine mesh $\microTess$, the set $\EK$ of computational edges of $K$  can be defined in different ways, depending on which vertices of $\microTess$ lying within the geometrical edge are retained as vertices of $\Tess$. }

	\begin{remark}\label{rem:additional}
		In the present framework, it is always possible to build a triangulation $\TriK$ of $K$, whose set of boundary edges coincides with $\EK$. In general, however, we can only guarantee that for $T$ in $\TriK$, $h_T/\rho_T \lesssim H/h$.
		Remark that the presence of very small edges is not problematic in itself, but, as illustrated in figure \ref{fig:additional}, the bad situations are rather the ones where very small edges are adjacent to large edges.	
	\end{remark}

\

Assumption \ref{shape_regular} is not sufficient to imply the validity of the standard shape regularity assumption on which the analysis of the virtual element method relies. In particular, it does not imply star shapedness of the element with respect to all points in a ball of radius of order $H$. We can however show that it is sufficient to obtain a number of bounds which are usually proven under more restrictive assumptions. In particular, the following bounds, which we prove in \ref{appendix:A}, hold with constants only depending on $\alpha_0$, $\alpha_1$ and $N_0$.

\begin{lemma}\label{lem:trace} Under Assumption \ref{shape_regular}, for all $\phi \in H^1(K)$ it holds that
	\begin{equation}\label{trace}
		\| \phi \|_{0,\bK} \lesssim \h ^{-1/2} \| \phi \|_{0,K} + \h ^{1/2} | \phi |_{1,K}.
	\end{equation}
	Additionally, for all $p \in \Poly{k}$  we have that
	\begin{equation}\label{tracep}
		\| p \|_{0,\partial K} \lesssim \h ^{-1/2} \| p \|_{0,K}.
	\end{equation}
Moreover, provided Assumption \ref{additional} also holds, for all $v \in H^{r}(K)$, $r \geq 1$, we have that
\begin{equation}\label{trace-r}
\sum_{e \in \EK} | v |_{r-1/2,e}^2 \lesssim | v |_{r,K}^2,
\end{equation}
where $\EK$ denotes the set of edges of the polygon $K$. 
\end{lemma}

For the sake of notational simplicity we also make the minor assumption that the set of boundary edges of $\Tess$ coincides with  $\Eb$, which, we recall, is the set of boundary edges of the fine mesh $\microTess$, that is,  all boundary nodes of $\microTess$ are also boundary nodes of $\Tess$ (see Figure \ref{fig:elements}). Remark however that all the result presented here carry over to the case where we allow also the boundary edges of $\Tess$ to be agglomerations of boundary edges of $\microTess$.



\

%

\subsection{The virtual element space} 

We will consider the standard order $\p \geq 2$ enhanced virtual element discretization space \cite{basicVEM}, whose definition we briefly recall.
For each  polygon $K \in \Th$ we let the space $\mathbb{B}_\p(\partial K)$ be defined as
\begin{gather*}
	\mathbb{B}_\p(\partial K) = \{ v \in C^0(\partial K): v|_{e} \in \mathbb{P}_\p\ \forall e \in \EK \},
\end{gather*}
where, we recall, 
$\EK$ denotes 
the set of edges of the polygon $K$.   A local space $\widetilde V^{K,\p}$ is defined as
\[
\widetilde V^{K,\p} =   \{ v \in H^1(K):\ v|_{\partial K} \in \mathbb{B}_\p(\partial K),\ \Delta v \in \mathbb{P}_{\p}\},
\]
and we introduce the operator $\PinablaK: H^1(K) \to \Poly{\p}$ defined as 
\[
\int_K \nabla \PinablaK v \cdot \nabla q = \int_K \nabla v \cdot \nabla q, \quad \forall q \in \Poly{\p}, \qquad \int_{K} \PinablaK v = \int_K v.
\]
The local VE space is then defined (\cite{3DVEM}) as 
\begin{equation}\label{defVEMK}
	\VEMK =   \{ v \in \widetilde V^{K,\p}:\ \int_K v q = \int_K \PinablaK v q, \ \forall q \in \Poly{\p} \cap  \Poly{\p-2}^\perp\},
\end{equation}
where $\Poly{\p} \cap\Poly{\p-2}^\perp$ denotes the $L^2(K)$ orthogonal complement of $\Poly{k-2}$ in $\Poly{\p }$. 
 The global discrete VE space $\VEM$ is finally obtained by glueing the local spaces continuously:
\begin{gather}\label{VEMm}
	\VEM = \{ v \in H^1(\PolyDom): v|_{K} \in \VEMK \ \forall K \in \Th \}.
\end{gather}
A function in $\VEM$ is uniquely determined by the following degrees of freedom 
\begin{itemize}
	\item its  values at the vertices of the tessellation;
	\item for each edge $e$, its values at the $k-1$ internal points of the $k+1$-points Gauss-Lobatto quadrature rule on $e$;
	\item  for each element $K$, its moments in $K$ up to order $k-2$. 
\end{itemize} 

The following  lemma, usually proven under stronger shape regularity assumptions, also holds under Assumption \ref{shape_regular}, as we show in  \ref{appendix:C}.

\begin{lemma}\label{lem:interpolation}
	For any given function $u \in H^2(\Omega)$ we can  define a unique function $u_I \in \Vh$ such that, if $u \in H^s(K)$ with $2 \leq s \leq \p+1$ we have
\begin{equation}\label{VEMapproxI}
 | u - u_I |_{1,K} \lesssim H_K^{s-1} | u |_{s,K}.
\end{equation}
	\end{lemma}

We recall now that a key concept in the  definition of  virtual element methods is the one of {\em computability}. Essentially,  operators or bilinear forms, acting on virtual element functions, are said to be computable if the knowledge of the degrees of freedom of the argument functions is sufficient for the direct evaluation of the operator/bilinear form, without the need of solving the PDEs implicitly involved in the definition of $\VEMK$. The elliptic projector $\PinablaK: \VEMK \to \Poly{\p }$ is computable \cite{basicVEM}, while 
%
  the bilinear form $a: \VEM \times \VEM \to \mathbb{R}$ and its local counterpart $a^K: \VEMK \times \VEMK \to \mathbb{R}$
\[ a(\phi,\psi) =  \int_\PolyDom \nabla \phi \cdot \nabla \psi, \qquad a^K(\phi,\psi) = \int_K\nabla \phi \cdot \nabla \psi,\]  are not. In the definition of the virtual element discretization, the bilinear form $a$ is replaced by a computable  approximate bilinear form $\ah : \VEM \times \VEM  \to \mathbb{R}$
\[
\ah (\phi,\psi) = \sum_K \ah ^K(\phi,\psi),
\]
where the elemental approximate bilinear form $\ah ^ K: \VKp \times\VKp \to \mathbb{R}$ is defined as
\[
\ah (\phi,\psi) = a^K(\PinablaK  (\phi), \PinablaK  (\psi) ) +  \stabweight\Svem(\phi - \PinablaK  (\phi), \psi - \PinablaK  (\psi)),
\]
$\stabweight > 0$ being a constant parameter to be chosen later, and the stabilizing bilinear form  $ \Svem$ being   any computable  symmetric  bilinear form satisfying
%
\begin{equation}
\label{defSK}
\cstar a^K(\phi,\phi) \leq  \Svem(\phi,\phi) \leq \Cstar a^K(\phi,\phi),\quad \forall \phi \in \VKp \ \text{ with } \Pinabla  \phi=0,
\end{equation}
with $\cstar$ and $\Cstar$ two positive constants independent of $K$. 
Different choices for the bilinear form $\Svem$ are available in the literature (see \cite{beirao_stab}), several of which reduce, when expressed in terms of the degrees of freedom, to a suitably scaled euclidean scalar product \cite{3DVEM2}. 
%

Letting $H^ 1 (\Tess)$ and $\Poly{\p}(\Tess)$  respectively denote 
 the spaces of discontinuous piecewise $H^1$ functions and of discontinuous piecewise polynomials of order up to $\p$,  defined on the tessellation $\Tess$:
\begin{gather*}
	H^ 1 (\Tess) = \{ v \in L^2(\Oh): \ v|_K \in H^1(K)\ \text{for all }K \in \Tess
	\},\\
	\Poly{\p}(\Tess) = \{
	v \in L^2(\Oh): \ v|_K \in \Poly{\p}(K)\ \text{for all }K \in \Tess
	\},
\end{gather*}
it will be convenient in the following to introduce the global projector $\Pn : H^1(\Tess) \to \Poly{\p}(\Tess)$ defined by $\Pn(v)|_K = \PinablaK(v|_K)$ for all $K \in \Tess$. Moreover, we let $\Pi^0 : L^2(\Oh) \to \Poly{\p}(\Tess)$ denote the $L^2(\Oh)$ orthogonal projection onto the space of discontinuous piecewise polynomials of degree at most $\p$.

\



\subsection{Nitsche's method with boundary correction for VEM}\label{sec:BDT}
%

In order to discretize Problem \eqref{prob_mod} we substantially follow \cite{VEM_weakly}, with some minor differences that will allow us to make the method more efficient (see Remark \ref{rem:condstat}). 
More precisely, we assume that $\Oh \subseteq \O$ and we choose, on $\bOh$, an outward direction $\sigma$, not necessarily normal to $\bO$ or $\bOh$, which we assume to be constant on each edge $\edge$  (see Figure \ref{fig:sigma}). For $u \in L^2(\Oh)$ with $u \in C^m(\bar  K)$ for all $K \in \Tess$, we let $\partial_\sigma^m u$ denote the $L^2(\bOh)$ function coinciding, on each boundary edge of $\Tess$, with the  $m$-th derivative of $u$ in the direction $\sigma$.
We recall that for 
$x \in \partial \Oh$, $\delta(x) > 0$ denotes the smallest non negative scalar such that 
\[\x + \delta(\x) \sx(\x) \in \partial \Omega.\]
We look for $\uh \in \VEM$ such that for all $v \in \VEM$ it holds that
\begin{multline}\label{problem_BDT_VEM}
\ah (\uh,v)  
- \int_{\bOh} \dn  \Pn(\uh) v  
\\ -  \int_{\bOh}  \left(
\Pn(\uh)+ \Corrstar {\Pn(\uh)}\right) \left(\dnh \Pn(v) - \gamma \h^{-1}(\Pn(v) + \Corrhat{\Pn(v)}) \right)
\\ = \int_{\Oh} f \Pi^0 (v) -
\int_{\bOh} \tg   \left(\dnh \Pn(v) - \gamma \h^{-1}(\Pn(v) + \Corrhat{\Pn(v)}) \right)
\end{multline}
with $\tg(\x) = g(\x + \delta(\x)\sigma)$, and with $\Corrstar{w}$ and $\Corrhat{w}$ defined in \eqref{definecorrections}.

\

%
The analysis of equation \eqref{problem_BDT_VEM} relies on the assumption that the discrete boundary $\bOh$ is sufficiently close to the true boundary $\partial \Omega$. More precisely, in order for \eqref{problem_BDT_VEM} to be well posed and yield an optimal error estimate, we need to assume that for some constant $\tau \in (0,1)$, sufficiently small, we  have that
\begin{equation}\label{deftau}
\max_{K\in \Tess} \max_{\x \in \partial K \cap \partial\Oh} \frac{\delta(\x)}{\h} \leq \tau.
\end{equation}
When $\Omega_h$ is constructed as the union of all elements in the fine mesh $\microTess$ wich are included in $\Omega$, such assumption reduces to 
\[
\widehat \tau :=	\frac h H \leq C \tau,
\]
the constant $C$ being the constant such that $\delta/\hK \leq C h/H$. 
If such an assumption is satisfied, existence and uniqueness of the solution of (\ref{problem_BDT_VEM}) can be proven by an identical argument to the one in \cite{VEM_curvo,VEM_weakly}, which also yields an error estimate (see \cite{VEM_weakly}) in the norm $\vvvert \cdot \vvvert_{\Oh}$, defined as
\[\vvvert \phi \vvvert_{\Oh}^2 = | \Pn \phi |^2_{1,\Tess} + | \phi - \Pn \phi |^2_{1,\Tess} + {H^{-1}}\| \Pn \phi + \Corrhat{\Pn \phi} \|^2_{0,\bOh},\ 
\]
where we set $| \phi |^2_{1,\Tess} = \sum_{K\in \Tess }| \phi |_{1,K}^2$.
 More precisely, we have the following theorem.

\begin{theorem}\label{cor:nitsche:curved}
	There exists $\stabweight_0> 0$ and $\gamma_0 > 0$ such that, provided $\stabweight> \stabweight_0$, for all $\gamma > \gamma_0$, the following holds: there exists a constant $\tau_\gamma$, depending on $\gamma$, such that if \eqref{deftau} holds for $\tau  < \tau_\gamma$, Problem \eqref{problem_BDT_VEM} is well posed, and,
	if  $u\in H^{k+1}(\Omega)\cap W^{m,\infty}(\Omega)$, $m \leq k+1$ we have the following error bound:
	\[
	\vvvert u - u_h \vvvert_{1,\Oh}  \lesssim H^{k} | u |_{k+1,\Omega} + H^{-1/2} h^{m} | u |_{m,\infty,\Omega}.
	\]
\end{theorem}

Observe that  as $\Pn + \widehat{\mathscr C} \circ \Pn$ preserves the constants,
$\vvvert \cdot \vvvert_{1,\Oh}$  is indeed a norm on $H^1(\Oh)$.
Theorem \ref{cor:nitsche:curved} has been proven in \cite{VEM_curvo} for the case $\khat = 0$ with a slightly different formulation of equation \eqref{problem_BDT_VEM}. The proof for the new, more general, formulation is fundamentally the same, but, for the sake of completeness we report it in \ref{appendix:B}.

\begin{remark}
	We point out that, as $\gamma_0$ is independent of the parameter $\cstar$ appearing in the condition \eqref{defSK} on the stabilization bilinear form $\Svem$, the  parameter $\gamma$ in \eqref{problem_BDT_VEM} can be chosen independently of the choice of $\Svem$.
\end{remark}

\begin{remark}
	By assuming $\Omega$ to be of class $C^\infty$, and allowing $\sigma$ to be an arbitrary direction, { we only displace the problems related with the choice of the extrapolation direction (particularly in the vicinity of the singular points of $\O$). Theorem \ref{cor:nitsche:curved} is valid for any choice of $\sigma$, but of course, if $\sigma$ is badly chosen, the conditions on $\delta$ required by such theorem for stability will be more difficult to satisfy.}
	The actual choice of $\sigma$ will be guided both by the result of Theorem \ref{cor:nitsche:curved} (in particular, by the condition $\tau < \tau_\gamma$ that suggests that $\sigma$ be chosen so that $\delta$ is as small as possible), and by practical considerations related to the numerical evaluation of the edge integrals in the bilinear form \eqref{bilinear_Nitsche}. These suggest to choose $\sigma$ so that $\sigma|_e$ is smooth for all edges even in the presence of corners in the physical domain. In this last case, we refer to \cite{atallah2021analysis} for a strategy for the choice of the extrapolation direction.
\end{remark}

\begin{remark}
	We recall that other approaches exist to adapting the VEM to deal with curved boundaries in 2D. Both \cite{VEM_curved_beirao} and \cite{VEM_curved_brezzi} propose versions of the VEM that handle elements with curved edges, which can be fitted to the actual physical boundary. Several extensions and applications are found in \cite{VEMcurvoext1,VEMcurvoapp1,VEMcurvoapp2,VEMcurvoapp3}.  As always, there are advantages and disadvantages in both approaches. In particular,  the use of elements with curved edges also allows to deal with curved interfaces, which we do not presently treat. Conversely, the boundary correction approach has the main advantage that, in a geometry adaptive procedure where the knowledge of the actual physical boundary (implicitely defined by the image) is iteratively improved, modifications of the boundary only imply  the modification of the correction term, with no need of recomputing the bulk integrals in the boundary adjacent elements. 
\end{remark}

\section{Static condensation of the ``lazy'' degrees of freedom}\label{sec:Lazy}
\begin{figure}
	\centering
	\includegraphics[width=5cm]{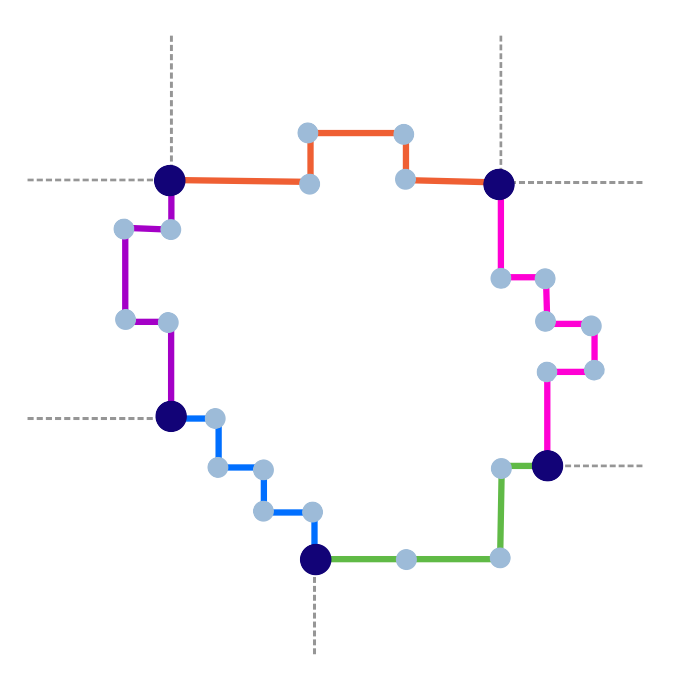}
	\caption{One element with 27 edges and 5 macro edges.}
	\label{fig:macroedges}
\end{figure}

While, depending on the ratio $H/h$, the number of degrees of freedom for the VE space $\VEM$ can be  much lower than the one for the finite element space $V_\mh$, it might still be quite high and, more importantly,  the discretization \eqref{problem_BDT_VEM} leads to a linear system with a matrix that may have some large dense blocks. This happens whenever an element $K$ presents a large number of small edges.  { When $\widehat \tau$ is large,} this always happens in our framework, at least for boundary elements. While in 2D this does not yet pose a major problem, in order for our approach to be viable also in 3D, we need to tackle this issue. { We observe that resorting to a classical static condensation approach \cite{guyan1965reduction}, or to other strategies aimed at reducing the number of interior degrees of freedom (\cite{da2016serendipity,wriggers2020serendipity,beirao2018serendipity,chen2023stabilization}), while feasible, would not really yield a solution to this problem. Indeed,  the number of  interior degrees of freedom only  depends on the polynomial degree but is independent of the refinement level of the underlying fine grid. On the other hand, the number of degrees of freedom lying on the boundary of the elements increases with both the polynomial degree $k$ and the refinement level $H/h$, and soon becomes dominant as $\widehat \tau$ increases. We then need to somehow locally eliminate  those degrees of freedom.}
In order to do so, we start with the observation that the bilinear form on the left hand side of  \eqref{problem_BDT_VEM} can be split into two components: a consistency component that only sees the test and trial functions after the action of the operator $\Pn$, and  a stability component, only needed for ensuring the well posedness of the discrete equation. The consistency component has, in our framework, a large kernel, whose elements are only seen by the stability part of  the bilinear form. In principle, we could then factor out such a kernel and restrict the discrete problem to a smaller space without losing the approximation properties. Unfortunately, the splitting of $\VEM$ into $\ker(\Pn) \oplus \VEMo$ (with $\VEMo \subset \VEM$ being a complement, not necessarily orthogonal, of $\ker(\Pn)$ into $\VEM$) is, a priori, a non local operation that would require a singular value decomposition of a large matrix. Then, a full computation of such a splitting is not a viable option. We can however single out  blocks of degrees of freedom on which such an operation can be performed locally. To do so, we start by introducing the macro edges of the tessellation:  
a macro edge $\mE$ is either a maximal connected component of $\partial K \cap \partial K'$, or a maximal connected component of $\partial K \cap \partial \Oh$,  with $K, K' \in \Tess$ (see Figure \ref{fig:macroedges}). We observe that, in our framework, it generally happens that, at least  macro edges on the boundary, but possibly also interior macro edges, are the union of a significant number of edges of the tessellation. 

\

Let then $E$ be a macro edge, and consider the subset $\VE \subset \VEM$ defined as
\[
\VE = \{ v \in \VEM : v = 0 \ \text{ on } \cup_K \partial K \setminus E, \ \int_ K v p = 0, \ \forall K \in \Tess,\ p \in \Poly{\p-2} \}.
\]
We easily see that, for a non zero $v \in \VE$,  $v|_K \not = 0$ if and only if $E$ is a macro edge of $K$. 
For $v \in \VE$, $E$ macro edge of $K$, it is not difficult to write down a necessary and sufficient condition for $\PinablaK v = 0$. Indeed, for $v \in \VE$ we have
\[
\int_K \nabla v \cdot \nabla p = \int_E v \nabla p \cdot \nK
\] 
whence 
\begin{equation*}
v \in \ker{\PinablaK} \Longleftrightarrow \int_E v \nabla p \cdot \nK = 0\ \text{ for all } p \in \Poly{\p }.
\end{equation*}
As $\nabla \Poly k \subseteq (\Poly{k-1})^2$, we also have that 
\begin{equation}\label{condlazy}
\int_E v \vec{p }\cdot \nK = 0\quad \text{ for all } \vec{p}=(p_1,p_2) \in (\Poly{k-1})^2 \Longrightarrow v \in \ker \PinablaK.
\end{equation}
We can immediately see that the condition \eqref{condlazy} is completely local to the macro edge. Moreover we observe that if $E$ is a macro edge common to $K$ and $K'$, then, as, on $E$,  $\nu_{K'} = - \nK$, we have that
\[
v \in \ker{\Pi^{\nabla,\p}_{K'}} \Longleftarrow \int_E v \vec p \cdot \nu_{K'} =  - \int_E v \vec p \cdot \nK = 0\ \text{ for all } \vec{p} \in (\Poly{\p - 1})^2,
\]
that is the sufficient conditions for having $v|_K \in \ker{\PinablaK}$ and $v|_{K'} \in \ker{ \Pi^{\nabla,\p}_{K'}}$ coincide. We can then split $\VE$ as 
\[
\VE = \widehat V_E \oplus \widecheck V_E\quad \text{ with }\quad \widecheck V_E = \{ v \in \VE : \int_E v  \vec p \cdot \nu_E = 0, \ \forall \vec p \in (\Poly{k-1})^2 \},
\]
where on a macroedge $E = \partial K \cap \partial K'$ we choose a normal direction, setting either $\nu_E = \nK$ or $\nu_E = \nu_{K'}$, and where $\widehat V_E$ satisfies $\dim(\widehat V_E) + \dim(\widecheck V_E)  = \dim(\VE)$. 
 In turn, this results in a splitting of $\VEM$ as 
\[
\VEM = \widehat V_\h \oplus\widecheck V_\h, \qquad \text{ with }  \widecheck V _\h= \oplus_{E} \widecheck V_E, \qquad \widehat V_\h = V^I_\h \oplus V^X_\h \oplus_E \widehat V_E,
\]
where $V^I$ and $V^X$ are, respectively, the space of functions in $\VEM$ identically vanishing on all edges of the tessellation $\Tess$ (that is the subspace of functions with vanishing vertex and edges degrees of freedom), and the space spanned by the basis functions corresponding to macro vertex (that is of vertexes of macro edges) degrees of freedom. We observe that for $\cv \in \cVEM$ not only we have that $\Pn(\cv) = 0$ (this is by construction), but we also have that $\Pi^0(\cv) = 0$.
In a way the functions in $\cVEM$ do not (at least not directly) contribute to the consistency/accuracy of the method, and we then dub them (and the corresponding set of degrees of freedom) as ``lazy''.

 If we  test equation \eqref{problem_BDT_VEM} with $v = \cv \in \cVEM$ we obtain the following identity
\begin{equation}\label{static}
\stab ( \hu - \Pn(\hu), \cv) + \stab(\cu,\cv) = 0,\quad \text{ where } \quad
\stab(\phi,\psi) = \sum_{K} S_a^K(\phi,\psi), 
\end{equation}
from which we immediately obtain
\[
\cu = -\Pi_S(\hu - \Pn(\hu)),
\]
where $\Pi_S$ is the projection onto $\cVEM$, orthogonal with respect to the scalar product $\stab(\cdot,\cdot)$. We can then plug back this expression in equation \eqref{problem_BDT_VEM} and set $v = \hv \in \hVEM$, obtaining the following reduced problem in $\hVEM$
\begin{multline}
	\int_{\Tess} \nabla \Pn(\hu) \cdot\nabla \Pn(\hv) - \int_{
		\bOh} \dn  \Pn(\hu) \hv  
\\	
-  \int_{\bOh} \left(
	\Pn(\hu)  + \Corrstar{\Pn(\hu)}\right) \left( \dnh \Pn(\hv)  - \gamma \h^{-1} (\Pn(\hv) + \Corrhat{\Pn(\hv)}) \right)\\
	+ s((1-\Pi_S)(\hu-\Pn(\hu)),(1-\Pi_S)(\hv-\Pn(\hv)))  \\
= \int_{\Oh} f \Pi^0 (\hv) -
\sum_{e \in \Eb} \int_e \tg   \left( \dnh \Pn(\hv)  - \gamma \h^{-1} (\Pn(\hv) + \Corrhat{\Pn(\hv)}) \right),
\end{multline}
where, to make the term involving the stabilization symmetric,  we exploited the fact that $1-\Pi_S$ is the self adjoint projector with respect to the scalar product $\stab(\cdot,\cdot)$. Remark that, compared to the full formulation in $\VEM$, the only term that is modified is the stabilization term. Moreover, if we choose $\stab$ in such a way that the corresponding matrix is diagonal, the computation of $\Pi_S$ turns out to be particularly cheap.

\begin{remark}
	\label{rem:condstat}
As already mentioned before, the method proposed in  \cite{VEM_weakly} is slightly different. More precisely, the discrete equation considered in such a paper is the following
\begin{multline*}
\ah (\uh,v)  
- \int_{\bOh} \dn  \Pn(\hu) v  
\\ - \int_{\bOh}  \left(
\uh + \Corrstar{\Pn(\uh)}\right) \left(\dnh \Pn(v) - \gamma \h^{-1} v \right)
\\ = \int_{\Oh} f \Pi^0 (v) -
\int_{\bOh} \tg  \left( \dnh \Pn (v) - \gamma H^{-1} v \right).
\end{multline*}  
With respect to the above formulation, besides considering a more general design of the Nitsche stabilization term, the equation \eqref{problem_BDT_VEM} replaces $u_h$ and $v$ on $\bOh$ with $\Pn(u_h)$ and $\Pn(v)$. Without this modification, the static condensation of the ``lazy'' degrees of freedom would require solving an equation of the form
\begin{multline*}
\stab(\cu, \cv)  + \gamma \h^{-1}\int_{\bOh} \cu \cv 
\\ = \int_{\bOh} \dn  \Pn(\hu) \cv  
- \gamma \h^{-1} \int_{\bOh}  \left(
\hu - \Corrstar{ \Pn(\hu)}\right)   \cv 
+ \gamma H^{-1}
\int_{\bOh} \tg   \cv 
\end{multline*}  	
instead of the simpler and much cheaper equation \eqref{static}. 
\end{remark}

{
	\begin{remark}
 We point out that, while we focus here on the Poisson equation, the underlying ideas can be applied to a much larger class of differential operator, including general elliptic operators with non constant coefficients. Moreover, the idea carries over also to the case of different, possibly non diagonal, VEM stabilization strategies, such as the ones proposed in \cite{beirao_stab}. In such case, however, the computation of $\Pi_S$ will require solving a block diagonal linear system and will not be as cheap.
	\end{remark}
}

\renewcommand{\SBM}{(A)}
\renewcommand{\BHL}{(B)}
\renewcommand{\BDT}{(C)}

\section{Numerical results}\label{sec:expes}
We devote this section to test the performance of the virtual element method with boundary correction  for increasing values of $k$, and to compare its performance with the performance of  analogous methods in the finite element framework.  We consider different variants of the proposed method, as illustrated in Table \ref{tab:SBMvariants}. { Remark that choosing $\sigma = \nu_h$, as proposed in the original version of the BDT method, is extremely disavantageous in the present framework. Consequently, we do not consider such an option.
For variants \SBM~ and \BHL~, the extrapolation direction will follow the SBM recipe, and will be chosen as the gradient of the distance to the boundary. For the case \BDT~, we slightly modify the SBM recipe, and let $\sigma$ be defined to be constant on each edge of $\Oh$, equals to  the gradient of the distance to the boundary  evaluated at the midpoint of each edge.
}
 \begin{table}
	\centering
{ 	\begin{tabular}{l|l||c|  l | c }
		 & Extrapolation direction& Stabilization & Ref.\\
			\hline
		\SBM	& SBM ($\sigma(\x) =    \nabla \dist(\x,\partial\O)$ )\Tstrut  &  SBM ($\widehat k = 1$)  & \cite{SBMho}\\
		\BHL & SBM ($\sigma(\x) =    \nabla \dist(\x,\partial\O)$ )\Tstrut &BDT   ($\widehat k = 0$) & \cite{burman2018cut} \\
		\BDT& p.w. constant SBM ($\sigma(\x) =    \nabla \dist(\x_e,\partial\O)$ )\Tstrut  & BDT   ($\widehat k = 0$)  & \cite{VEM_weakly}
		\end{tabular}
\caption{The three  boundary correction strategies considered for the numerical tests. Following the SBM recipe, for all three methods the extrapolation direction  is given by the gradient of the distance to the boundary, which for the \BDT~ case is evaluated at the midpoint $\x_e$ of each edge. }	\label{tab:SBMvariants}}
\end{table}

We tested the proposed method on two different curved domains, with different characteristics, namely a disk and a {  non convex} bean shaped domain with curved boundary, and a reentrant corner with interior angle equals to $3\pi/2$. 
We consider different meshes, obtained by agglomerations of squared elements from uniform structured meshes of meshsize $h$, for different values of the ratio $\widehat \tau = h/H$,  which, unless otherwise stated, is the same for all the elements of a given  mesh. We underline that, as stated in Remark \ref{rem:additional}, the considered tessellations automatically satisfy Assumption \ref{additional} with $\alpha_1 \lesssim \widehat \tau^{-1}$.

\

Letting $u_h$ denote the discrete solution obtained by the order $\p$ VEM method proposed in the previous section, for all the tests we consider the  relative error in the $H^1(\Th)$ seminorm, as well as in the $L^2(\Oh)$ norm. For the VEM case, these are, as usual, approximated as
\begin{align}\label{eS}
	e^u_1 &:= \frac{
		\| \nabla u - \Pi^0_{k-1}(\nabla u_h) \|_{0,\Oh}
	}{| u |_{1,\Oh}}, & e^u_0 &:= \frac{\| u - \Pi^0_{k}  u_h\|_{0,\Oh}}{\| u \|_{0,\Oh}},
\end{align}
where $\Pi^0_\ell : L^2(\Oh) \to \mathbb{P}_\ell(\Th)$ is the $L^2$ orthogonal projection onto the space of discontinuous piecewise polynomials of order up to $\ell$.
For both test cases, we consider values of $k$ between $1$ and $6$.
To facilitate the interpretaton of the results, which, for all tests, we plot in log-log scale, for each set of data we also plot a dotted line which we fit, by linear regression,
to the central part of the data set (we exclude one or two values, that we estimate less relevant, at each end). The slope of such a line is reported in the plots.

\subsection{Test 1 -- Disk domain}\label{sec:test1}
For the first test, the domain is the disk of center $(0.5,0.5)$ and radius $0.5$. We solve Problem \eqref{prob_mod} with  data chosen in such a way that the solution $u$ of the problem is the Franke function \cite{franke1979critical}
\begin{multline}
	\Franke(x,y):=\frac 3 4 e^{-\left(
			(9x-2)^2+(9y-2)^2
			\right)/4} + \frac 3 4 e^{-
			\left(
			(9x+1)^2/49 + (9y+1)/10\right)}\\
	+ \frac 1 2 e^{-
			\left(
			(9x-7)^2 + (9y-3)^2	
			\right)/4
		}
	+\frac 1 5 e^{-
			\left(
			(9x-4)^2 + (9y-7)^2
			\right)	
		} . \label{franke}
\end{multline}

As for the boundary correction strategy, we test both the shifted boundary method (strategy \SBM) and the Bramble, Dupont and Thomée method with piecewise constant closest point extrapolation direction, as proposed in \cite{VEM_weakly} (strategy \BDT).
For the sake of comparison, we also present the results of the shifted boundary method and of the Bramble, Dupont, Thomée method with closest point extrapolation direction (strategy \BHL, proposed in \cite{burman2018cut}), in combination with order $k$ rectangular finite elements on the structured uniform grid corresponding to the $\widehat \tau=1$ case. For the first set of tests we set the stabilization parameter for the Nitsche's method $\gamma = 100$.

\begin{figure}[htpb]
	\centering
	\includegraphics[width=5.5cm,height=4cm]{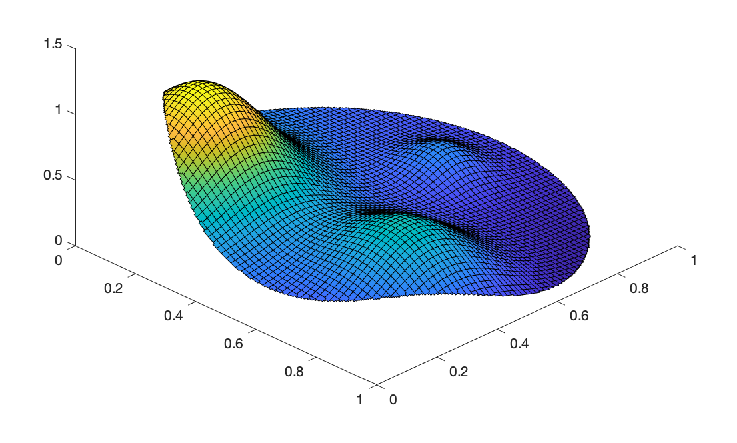}
				\includegraphics[width=4.5cm,height=4cm]{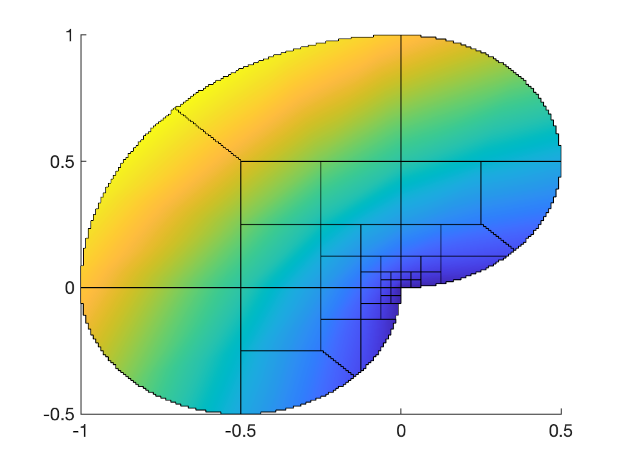}
		\includegraphics[width=4.5cm,height=4cm]{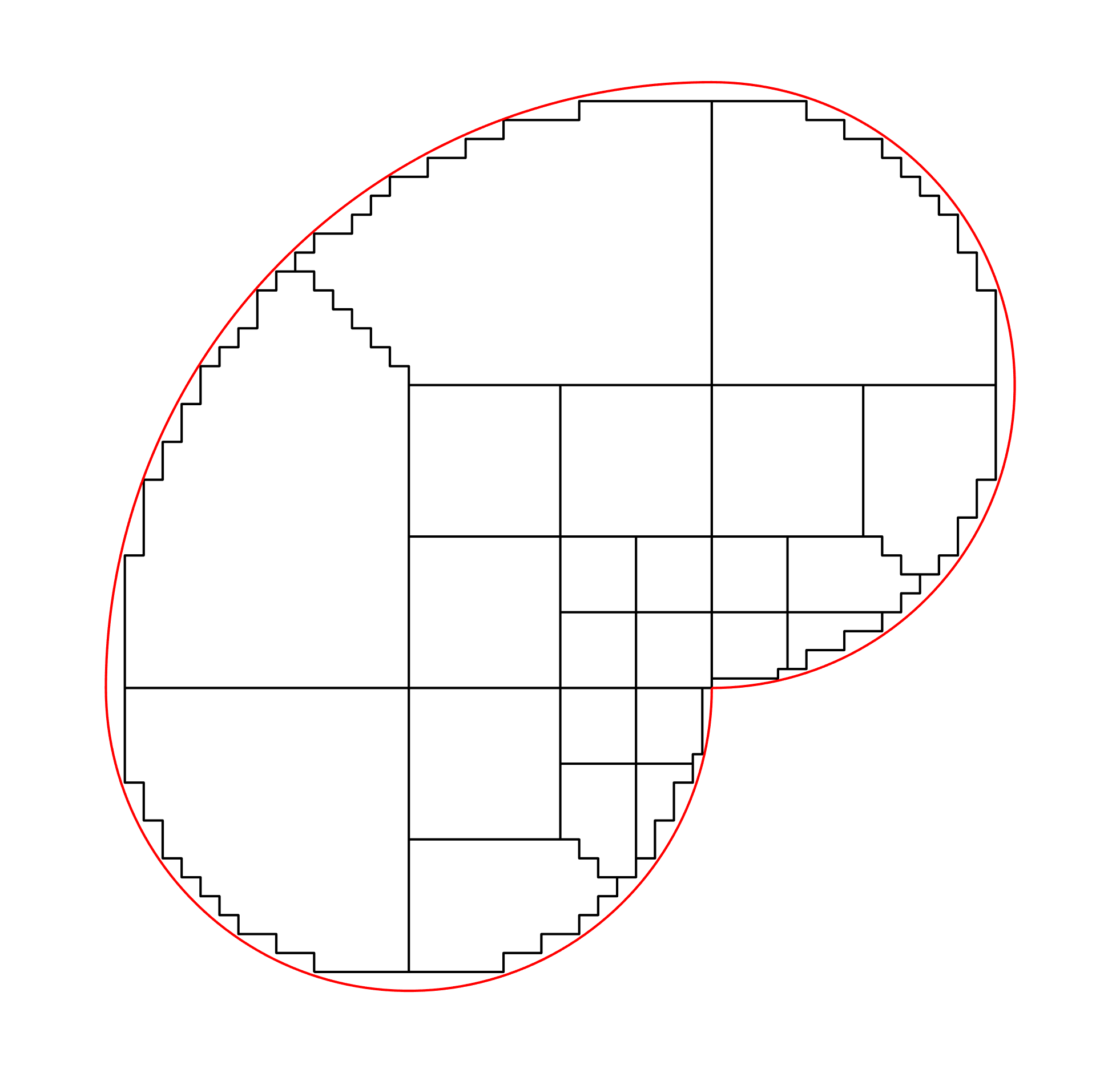}	
	\caption{The solution to test problem 2 (top) and one of the meshes used in the tests. Remark that the approximate domain $\Omega_h$ is included in $\Omega$.}\label{fig:solution}
	\end{figure}

\subsubsection{Effect of the refinement parameter $\widehat \tau$}
\begin{figure}[h]
	\centering
	\includegraphics[width=4.5cm]{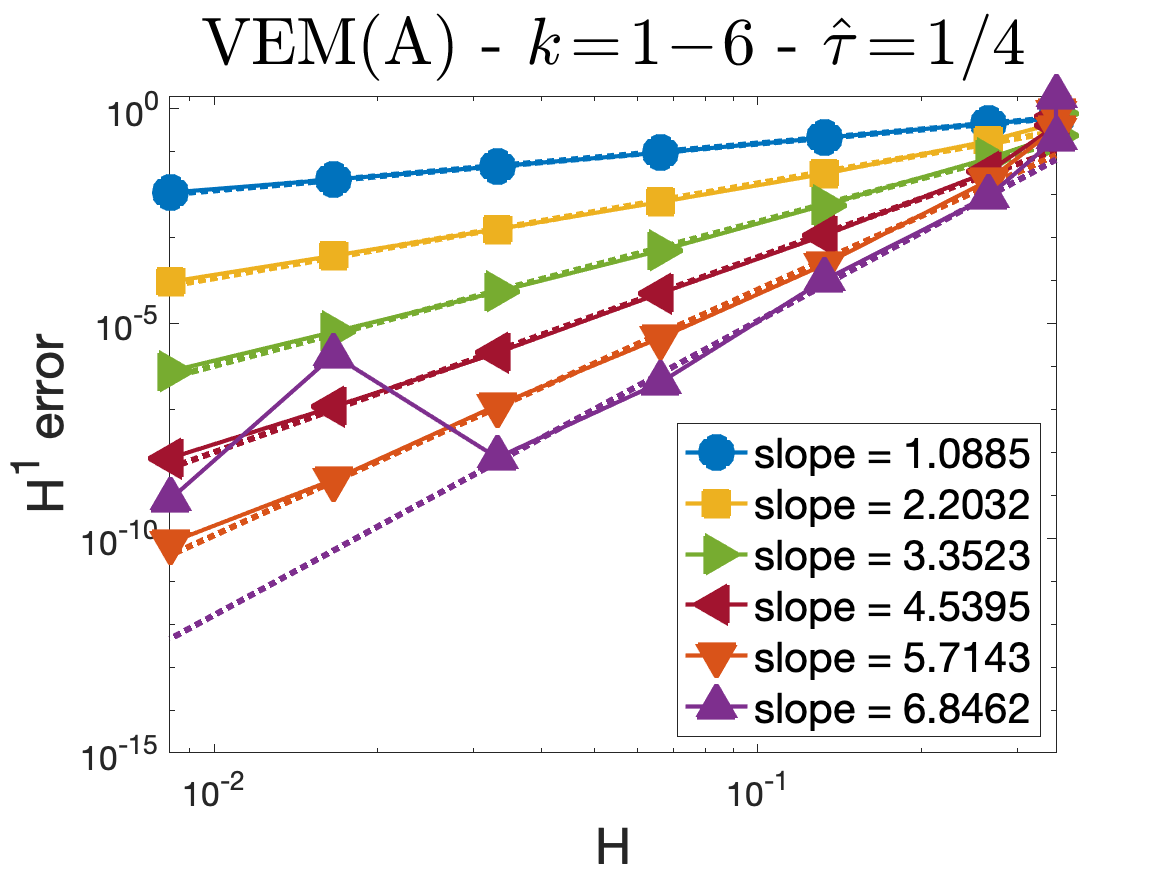}
	\includegraphics[width=4.5cm]{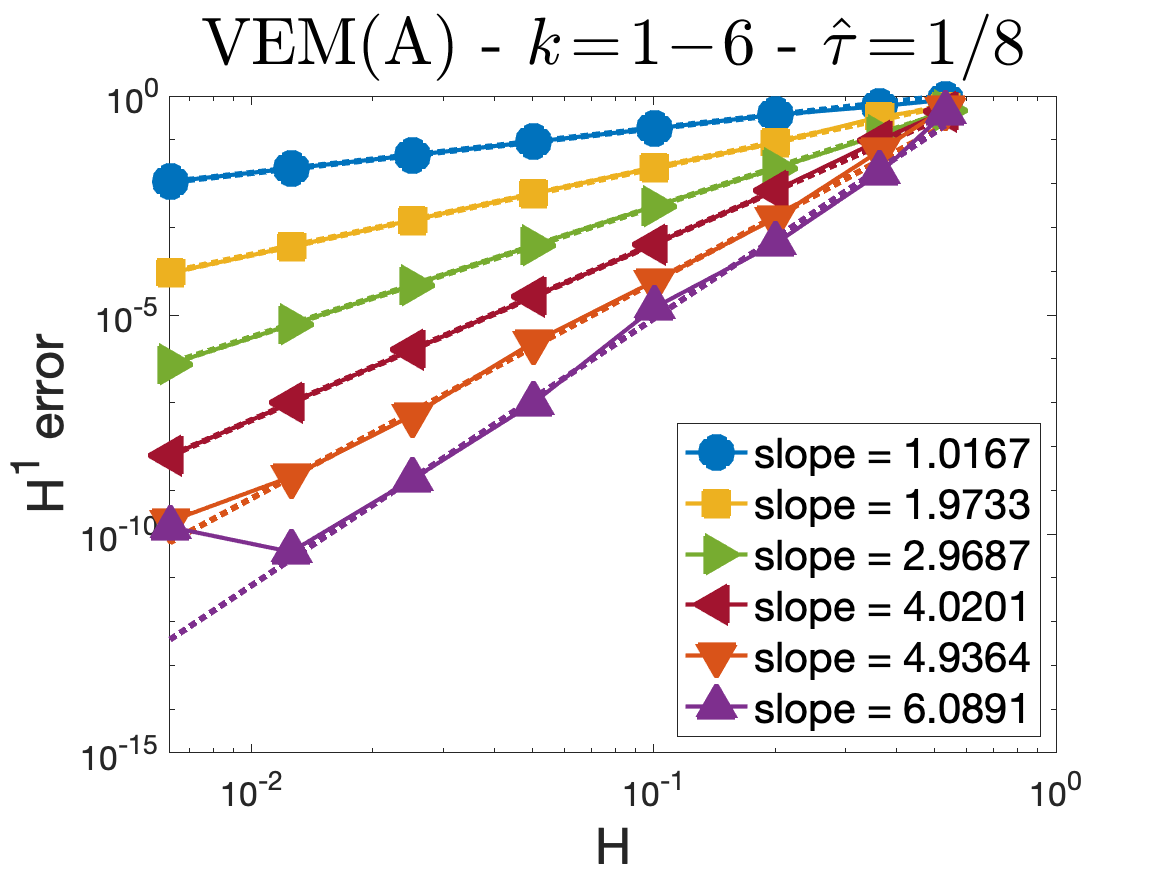}
	\includegraphics[width=4.5cm]{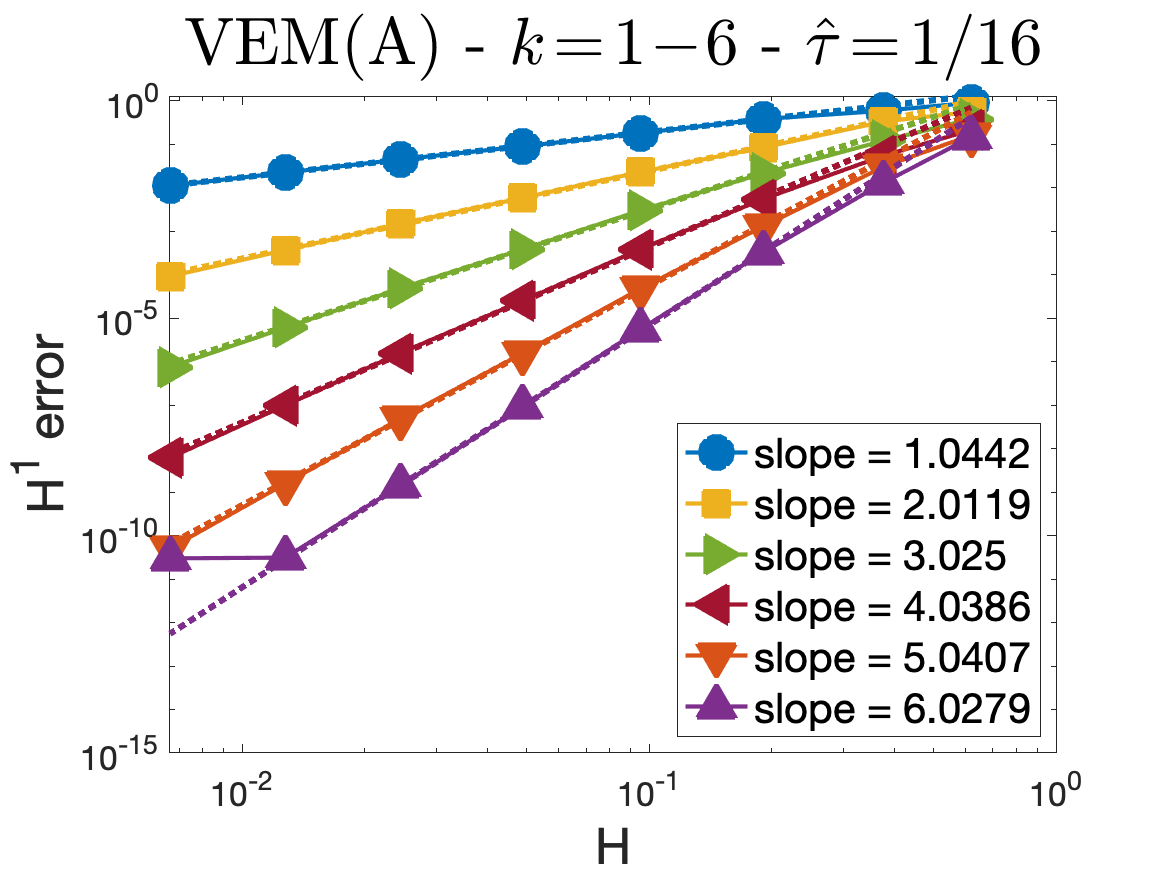}\\
	\includegraphics[width=4.5cm]{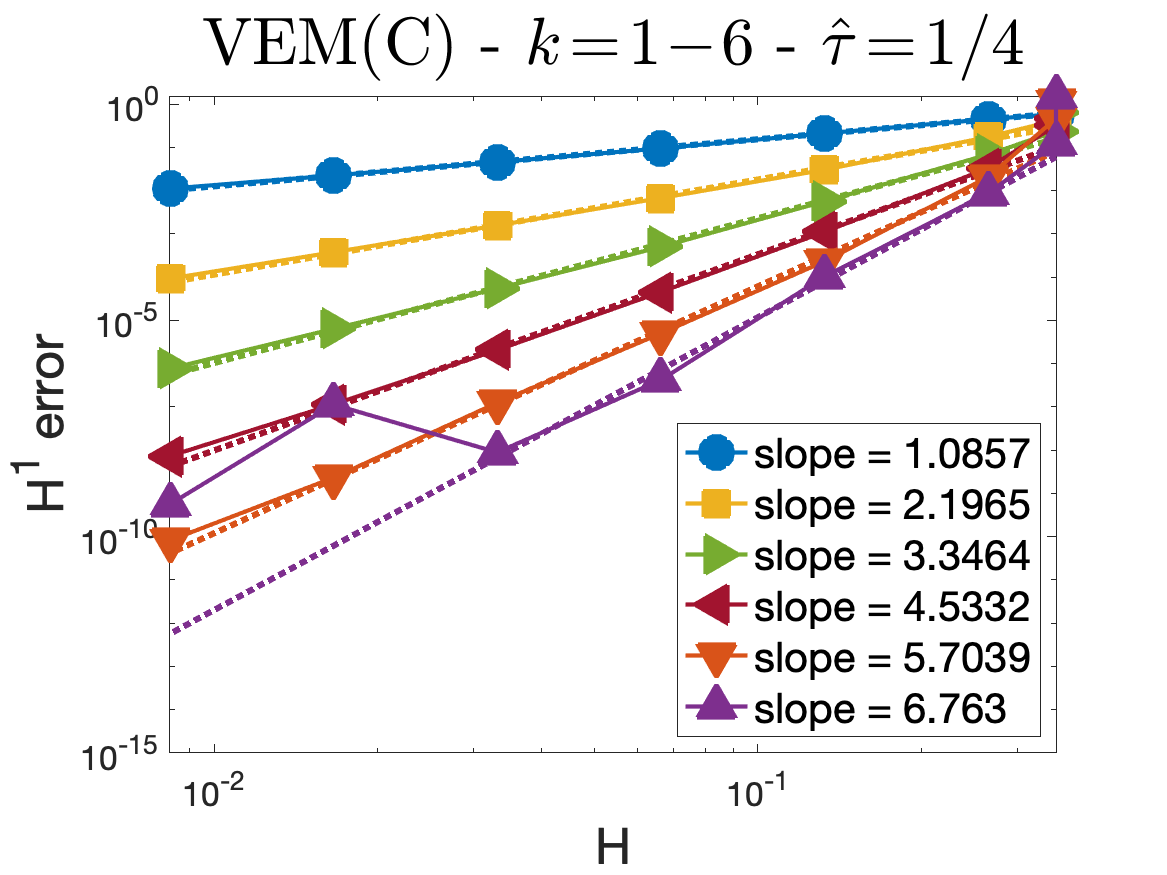}
	\includegraphics[width=4.5cm]{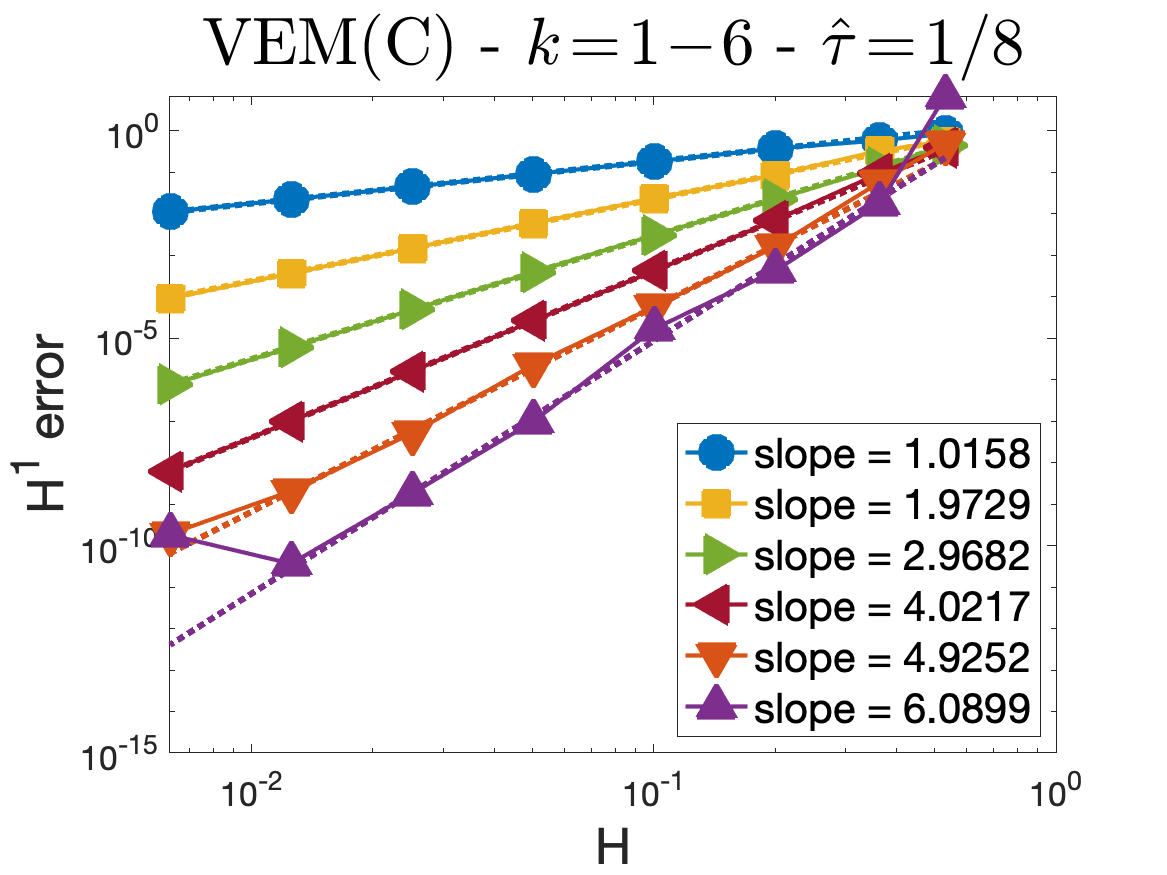}
		\includegraphics[width=4.5cm]{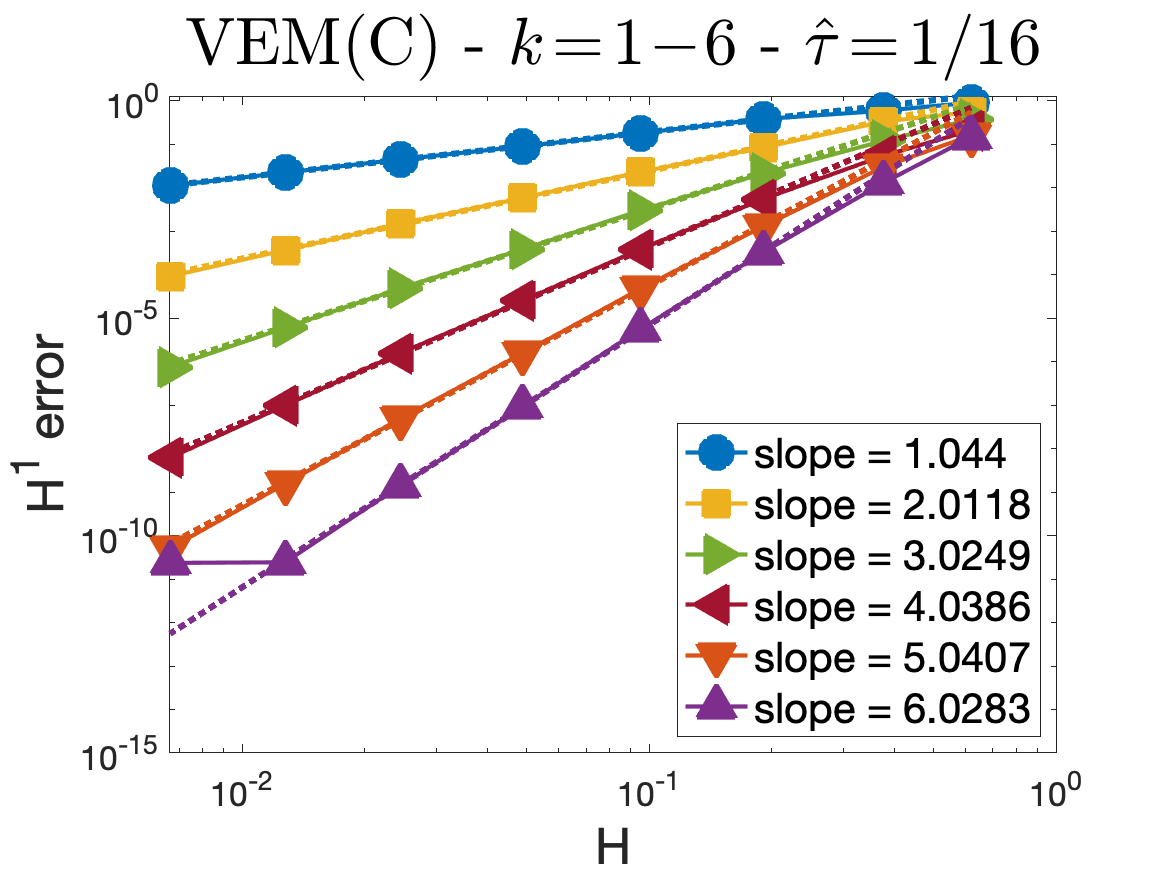}
	\caption{
			Test case 1,  $H^1$ convergence of the VEM method with boundary correction strategies \SBM~and \BDT~  for $k=1,\cdots,6$, with different values of $\widehat \tau= h/H$. 	{For $\widehat \tau = 1$ and $k\geq 5$ the method displays evident instabilities. For $\widehat\tau = .5$ the method is only slightly suboptimal. For $\widehat \tau \leq .25$ we observe an optimal behavior for all $k\leq 5$, and only some very mild oscillations for $k = 6$.} The dotted lines are obtained by linear regression fitting to a subset of the data that excludes the  coarsest as well as the two finest meshes. 
	 }\label{fig:test1-g100-refs-H1}
\end{figure}

\begin{figure}[h]
	\centering
	\includegraphics[width=4.5cm]{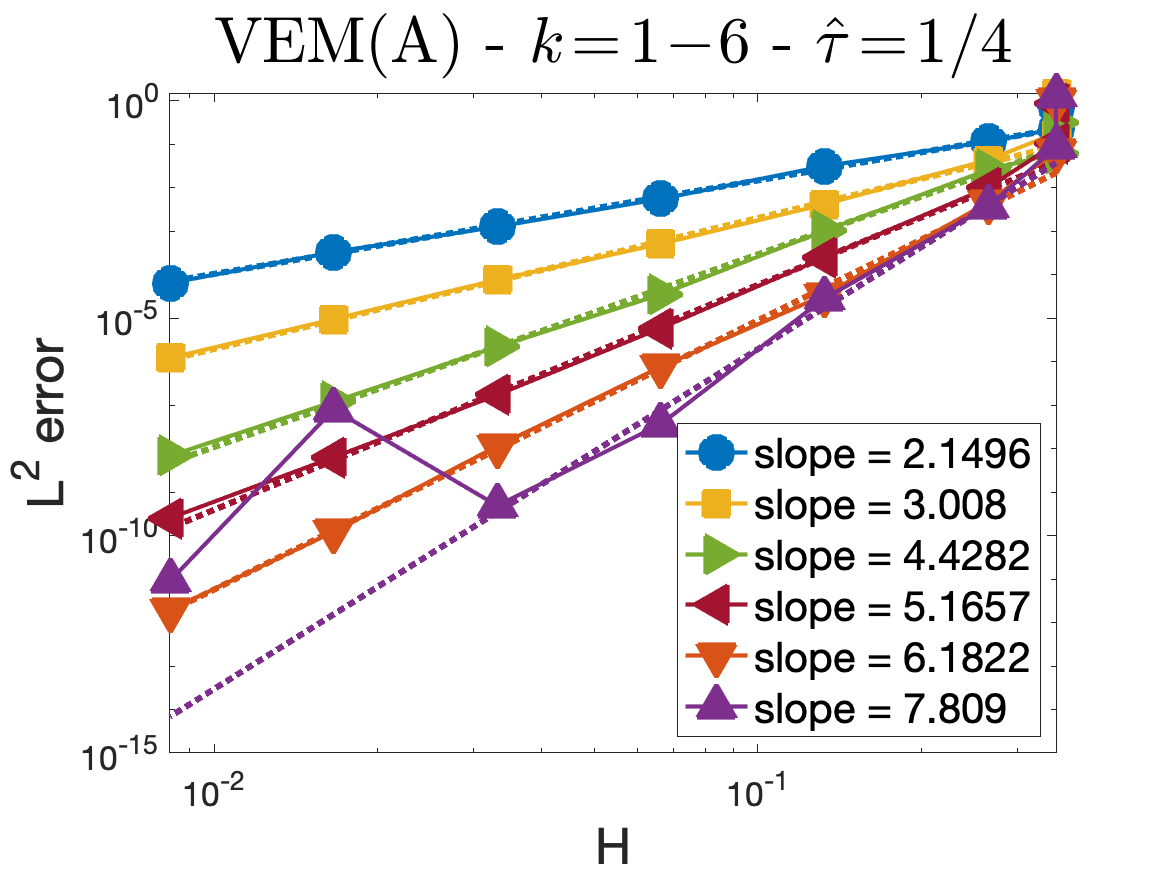}
	\includegraphics[width=4.5cm]{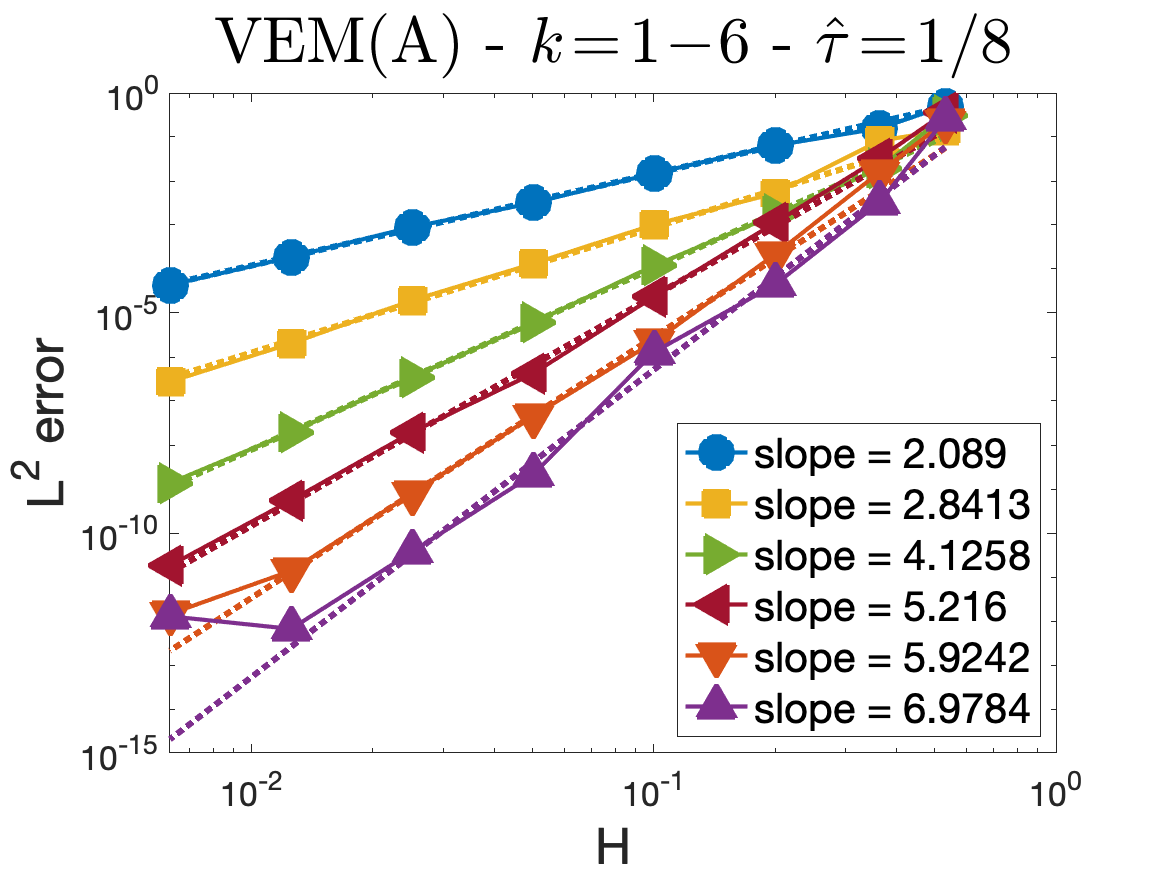}
	\includegraphics[width=4.5cm]{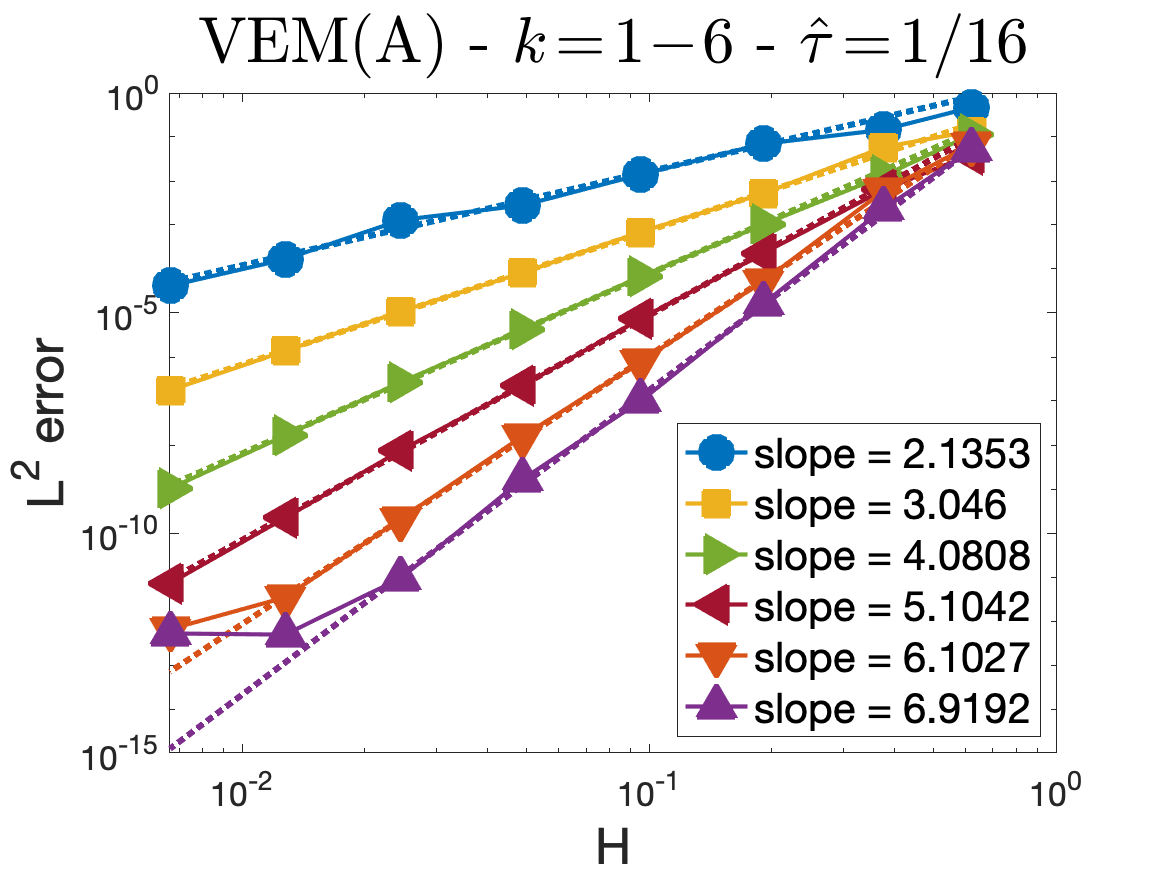}\\
	\includegraphics[width=4.5cm]{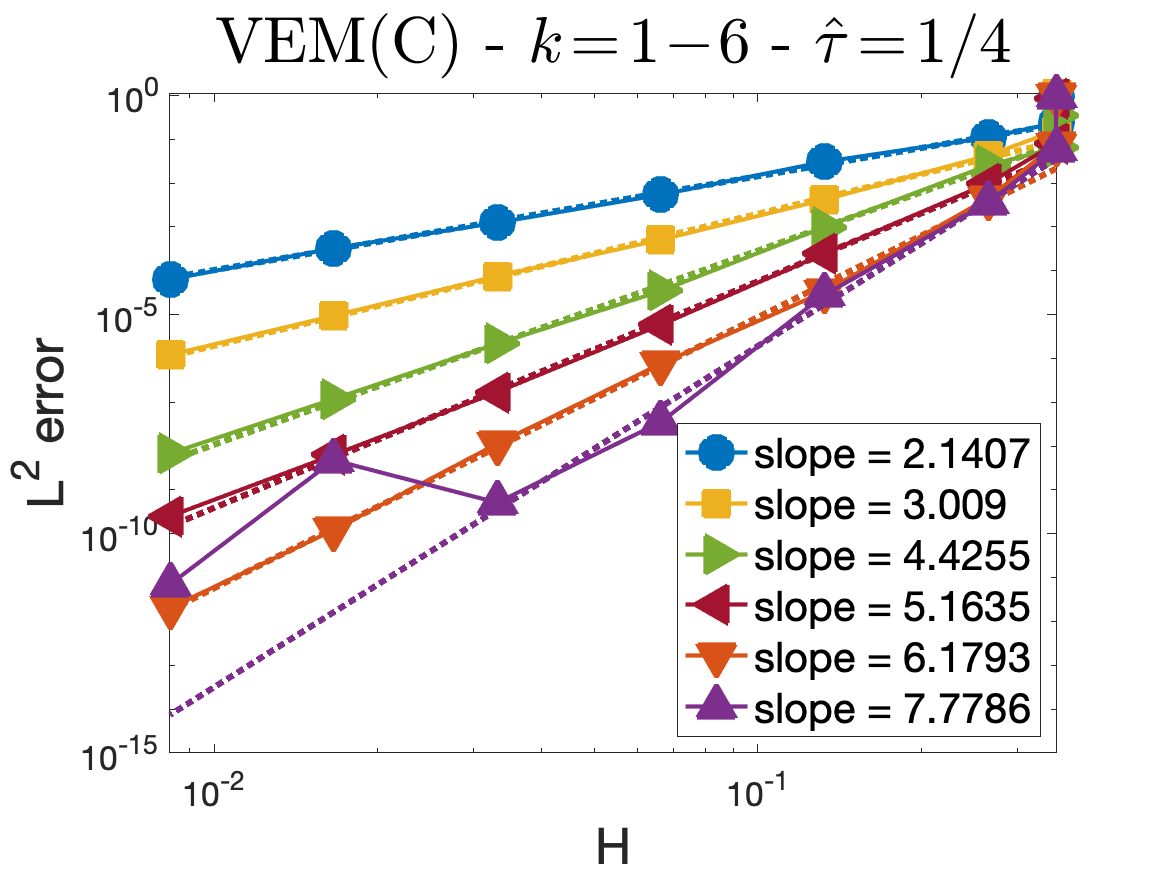}
	\includegraphics[width=4.5cm]{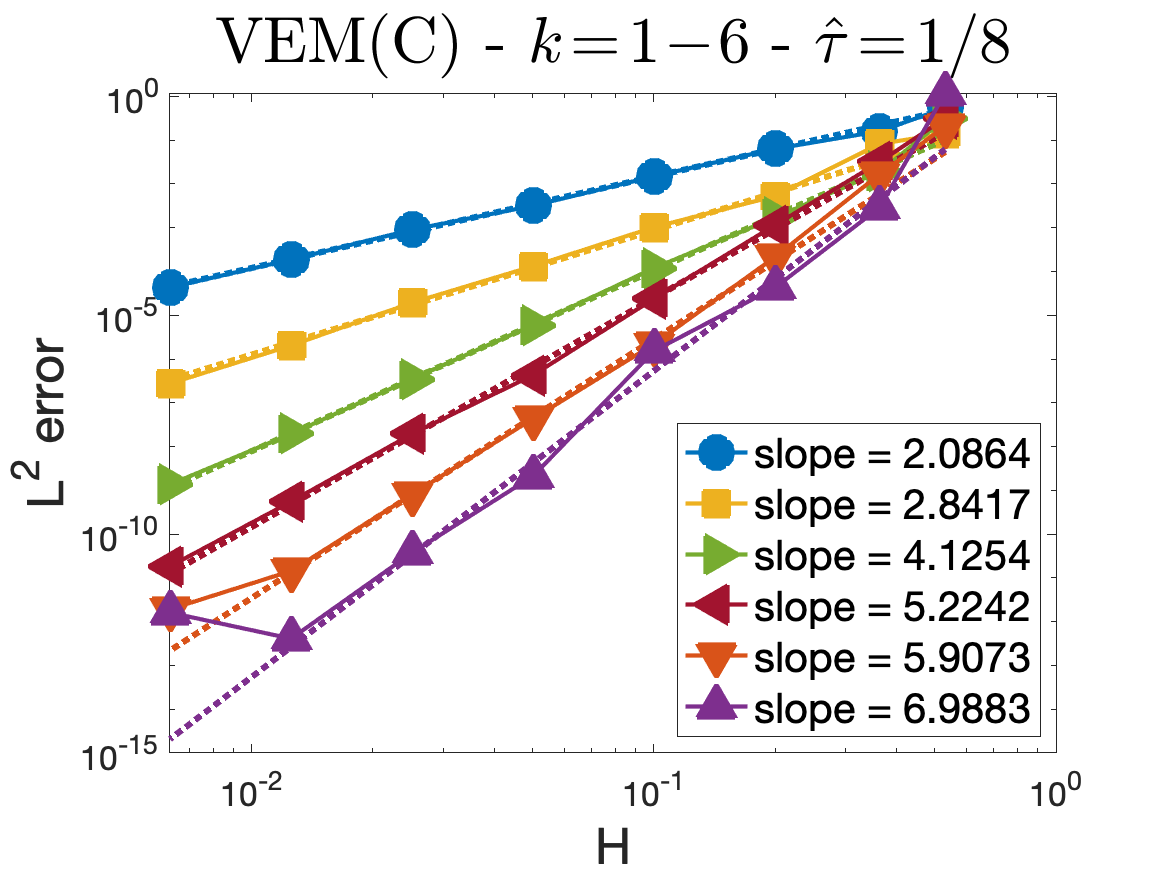}
	\includegraphics[width=4.5cm]{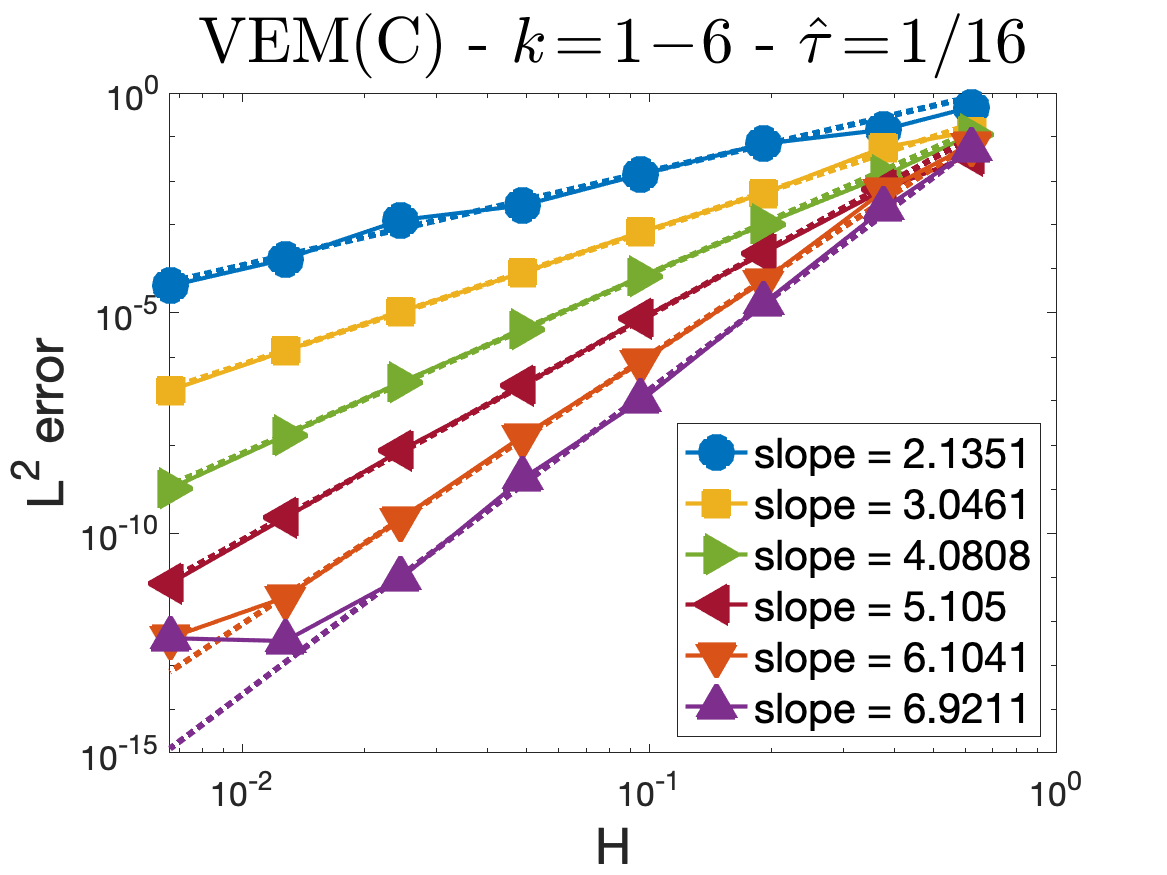}
	\caption{
		Test case 1,  $L^2$ convergence of the VEM method with boundary correction strategies \SBM~and \BDT~  for $k=1,\cdots,6$, with different values of $\widehat \tau= h/H$.  The dotted lines are obtained by linear regression fitting to a subset of the data that excludes the  coarsest as well as the two finest meshes. The behavior is similar to the one of the $H^1$ norm.
	}\label{fig:test1-g100-refs-L2}
\end{figure}

  To demonstrate the effect of the choice of the refinement parameter $\widehat\tau$, we start by  considering the case of no refinement ($\widehat \tau = 1$), and we test order $k$ finite elements and virtual elements combined with different boundary correction versions, { on standard triangular and quadrilateral meshes.  For all the versions of the method we choose the parameter $\gamma$ according to the recipe proposed in \cite[Section 10.1]{SBMho}, where, for quadrangular elements, the inverse inequality constant $C_{inv}$ is computed according  to  \cite{harari1992c}. 
We observe that, for $k \leq 3$, all versions of the boundary correction method yield optimal convergence. On the other end of the polynomial order range, for $k=6$, all versions lose optimality, in a more or less pronounced way,  as $H = h$ becomes small. In general, the \SBM~ version (which for finite elements coincides with the high order SBM proposed in \cite{SBMho}) appears to behave better than the \BHL~and \BDT~versions, due, we believe, to the better property of the SBM stabilization term. Also in general the $Q_k$ version of FEM behaves more poorly than the $P_k$ FEM and the order $k$ VEM. We believe that this is due to the fact that the inverse inequality constant for the space $\PolyQ k$ of polynomials of degree less than or equal to $k$ in each variable is strictly greater then the one for $\Poly{k}$.

Switching to  polygonal elements, and focusing on VEM on quadrangular meshes, the situation  improves as we take smaller values of $\widehat \tau$. In such case, however, due to the difficulty of evaluating the inverse inequality constant, we were unable to use the recipe of  \cite[Section 10.1]{SBMho} to evaluate the best value for $\gamma$. In the following tests we took $\gamma =100$.}  In Figures \ref{fig:test1-g100-refs-H1} and \ref{fig:test1-g100-refs-L2} we present the results of our test on polygonal meshes obtained by agglomeration of uniform square meshes, for $\widehat \tau = 0.5$, $0.25$ and $0.125$. As $\widehat \tau$ decreases, the instabilities  progressively disappear and, for $\widehat \tau = .125$, the behavior of the $H^1$ error for both strategies \SBM~ and \BDT~ appears optimal for all $k$s in our range (we believe the increase of the error for $k=6$ on the finest mesh to be the result of round-off errors). 	We can however still see some slight oscillations in the $L^2$ norm for strategy \BDT. We also tested lower values of $\widehat \tau$, namely $0.0625$ and $0.03125$, with similar results, that we do not report here.

\begin{figure}[h]
	\centering
	\includegraphics[width=3.65cm]{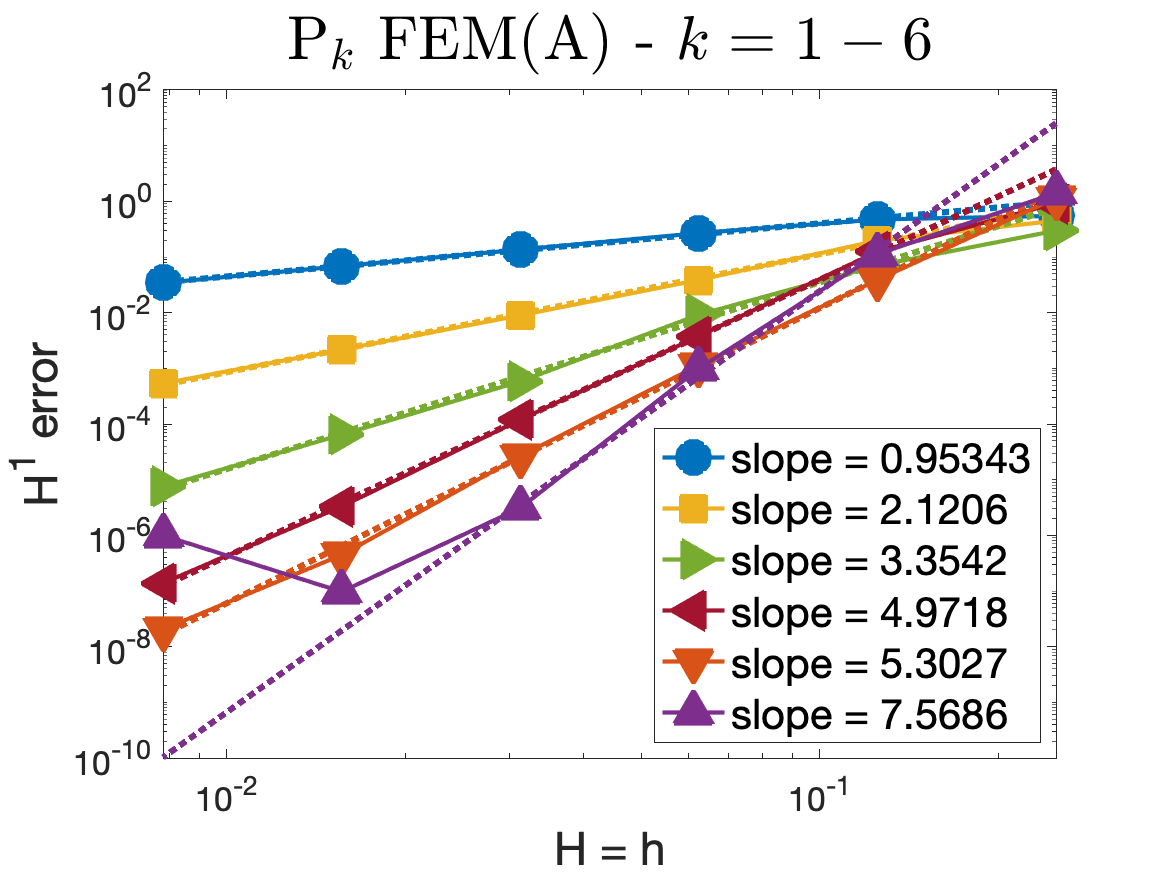}
	\includegraphics[width=3.65cm]{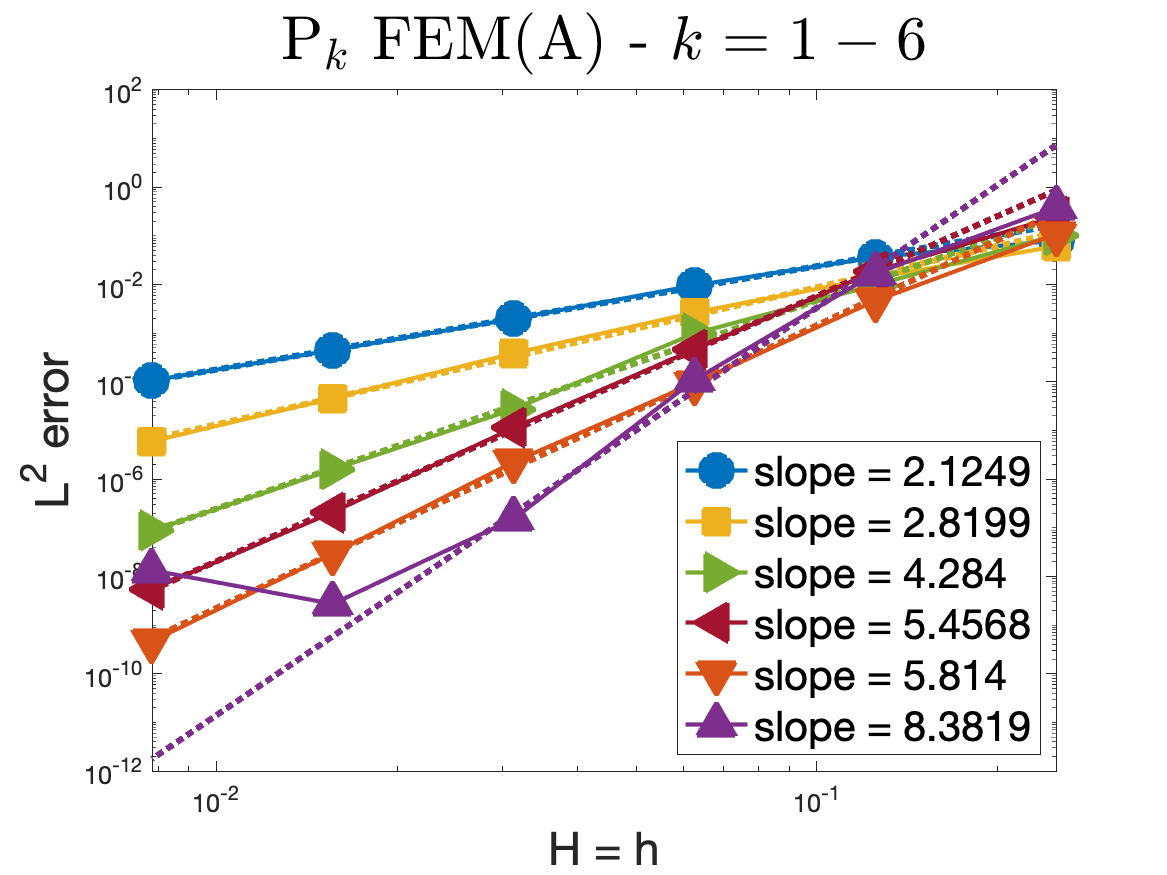}
	\includegraphics[width=3.65cm]{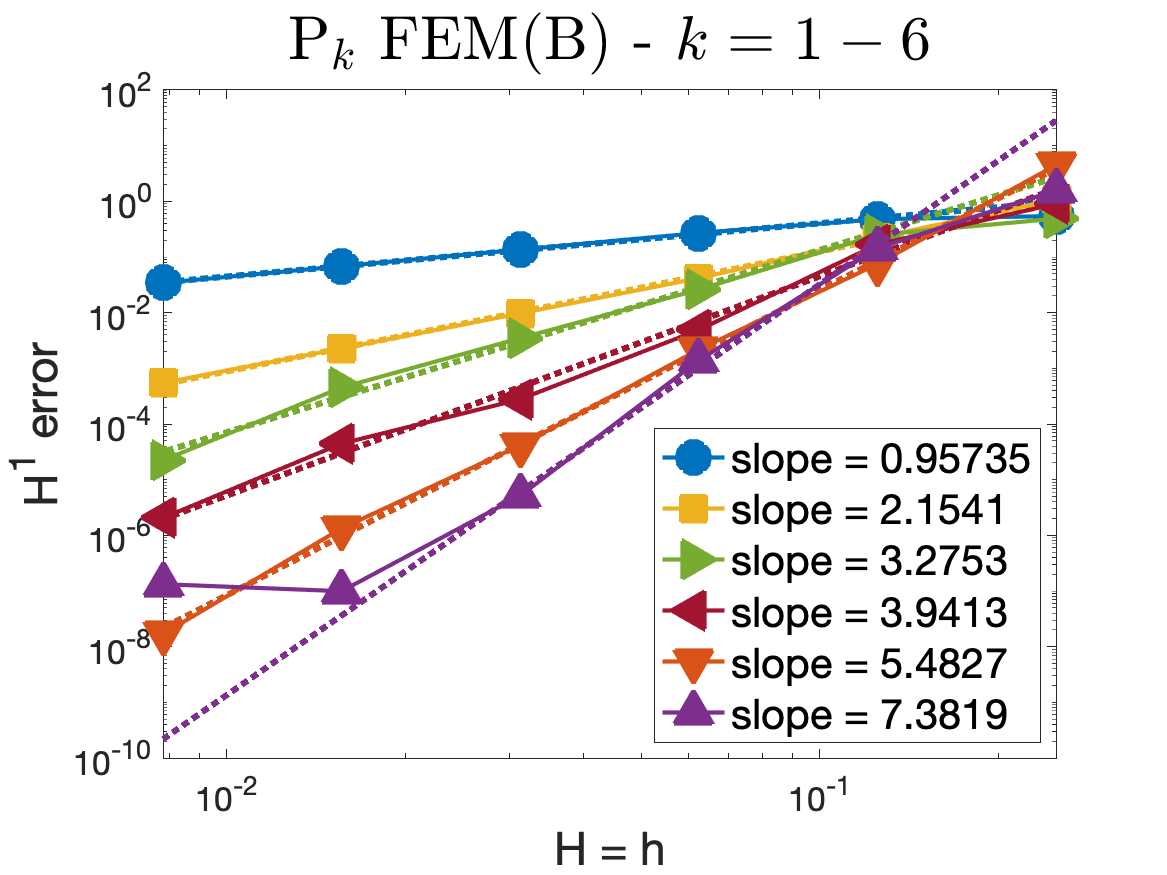}
	\includegraphics[width=3.65cm]{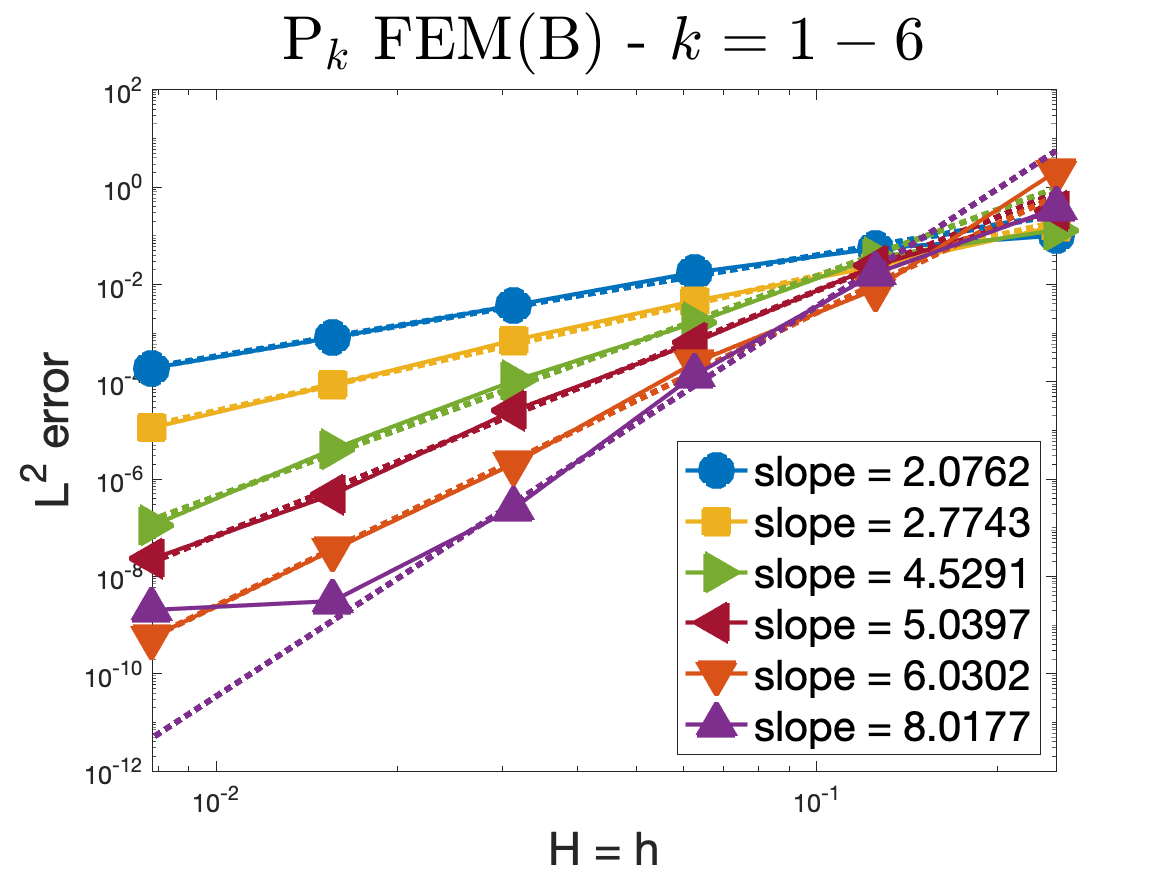}\\
	\includegraphics[width=3.65cm]{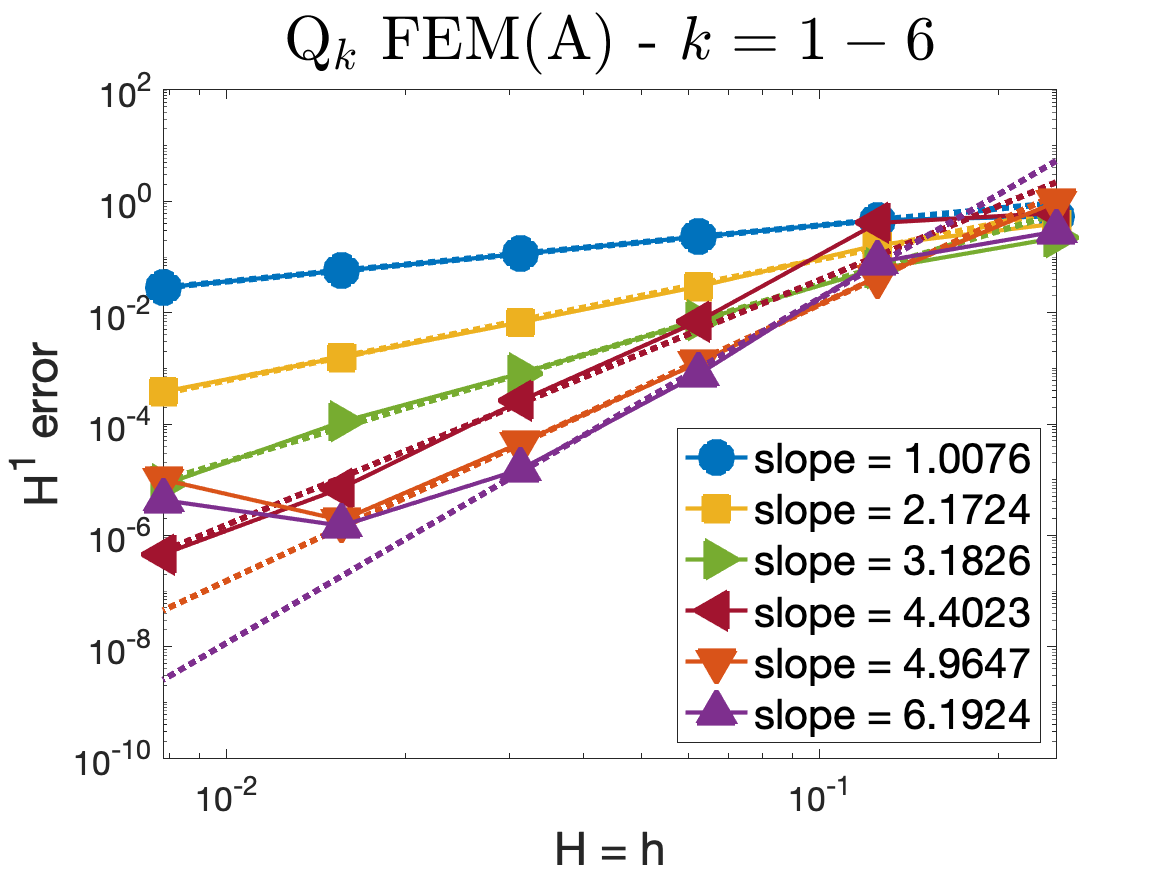}
	\includegraphics[width=3.65cm]{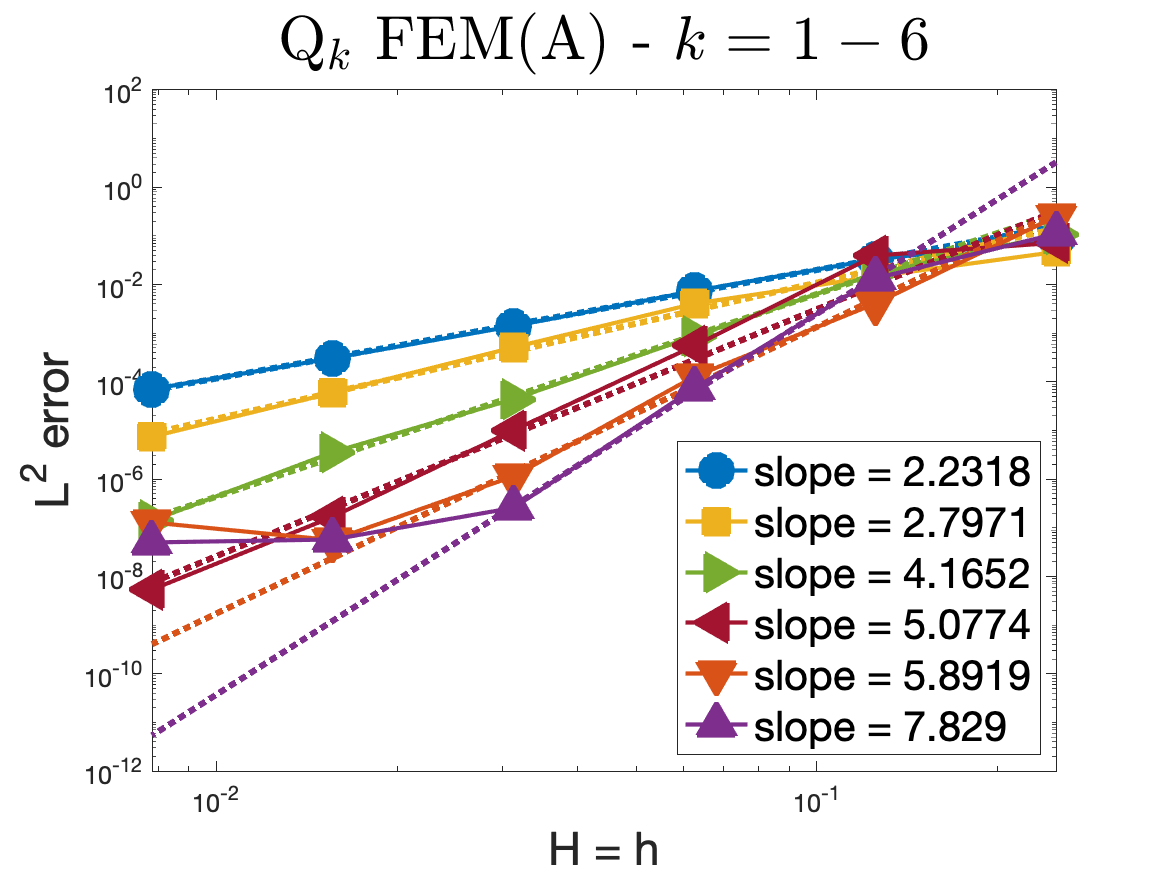}
	\includegraphics[width=3.65cm]{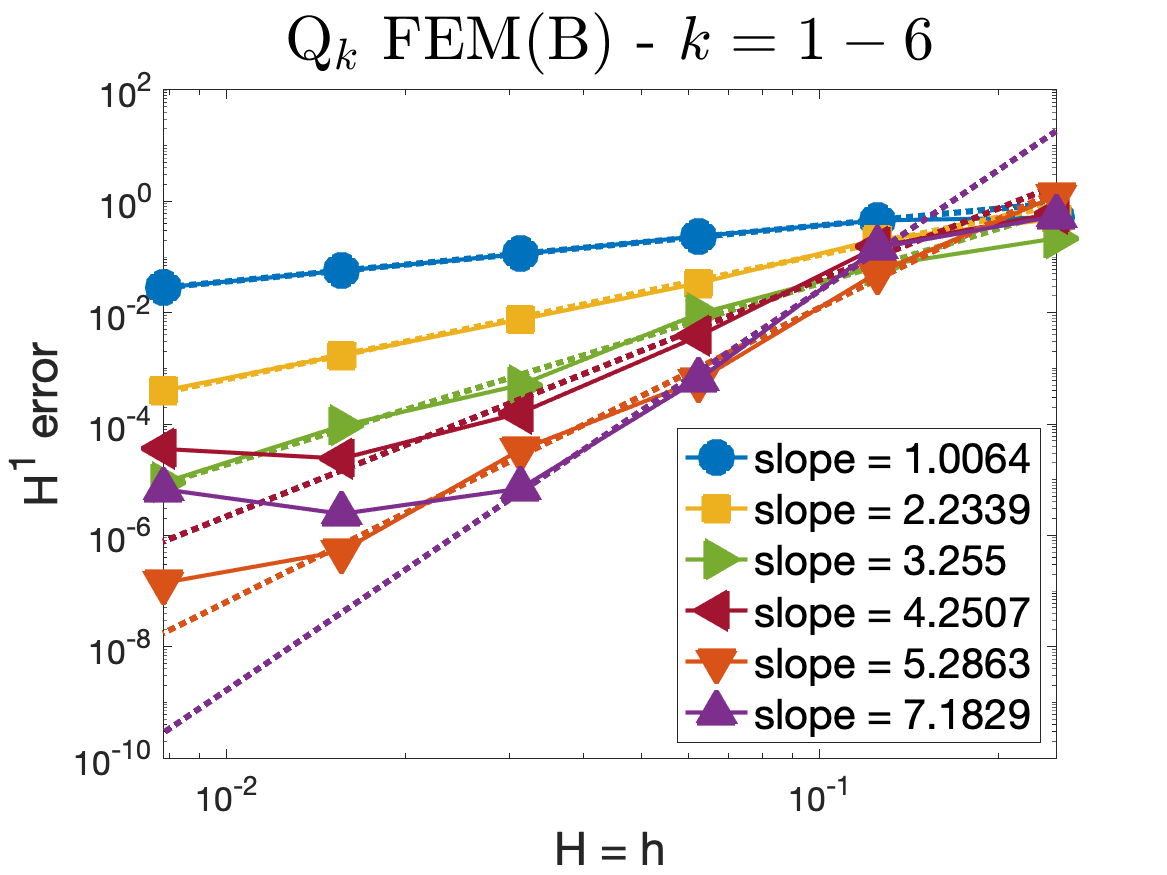}
	\includegraphics[width=3.65cm]{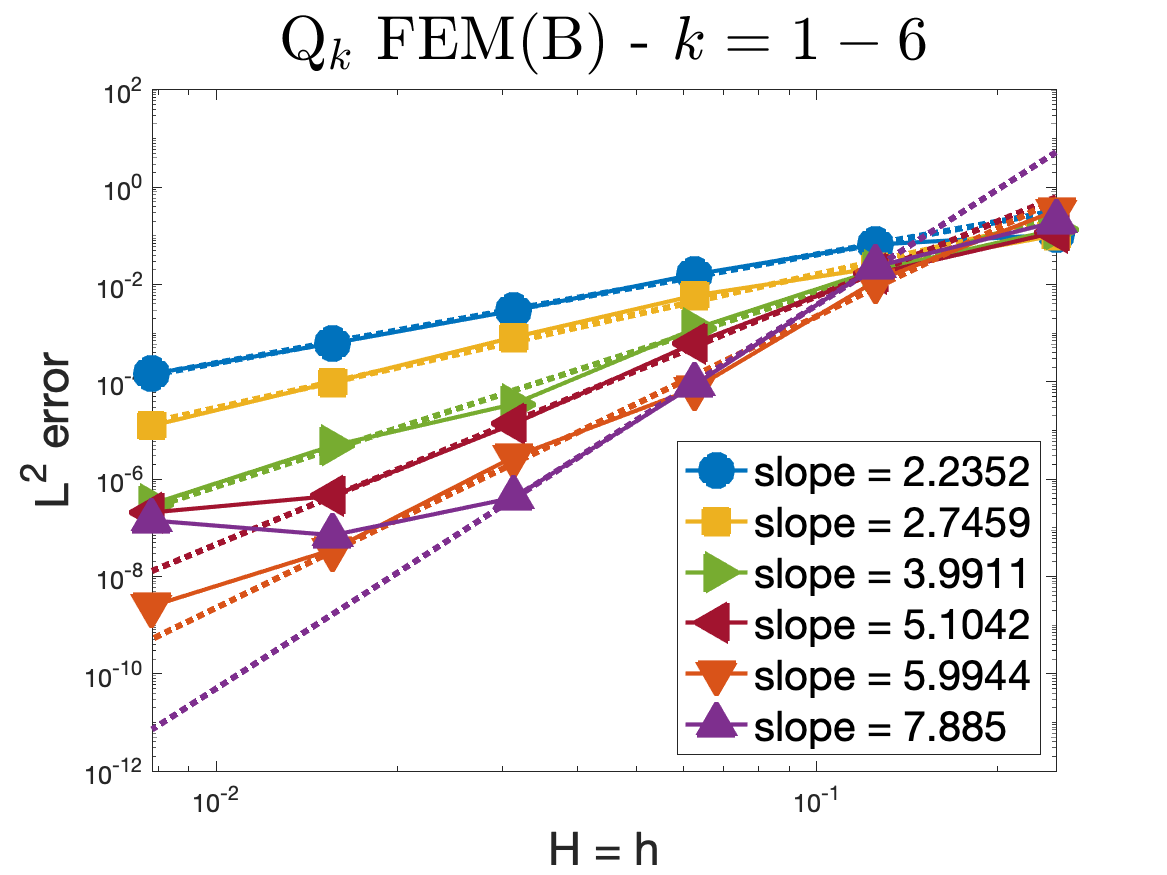}\\
	\includegraphics[width=3.65cm]{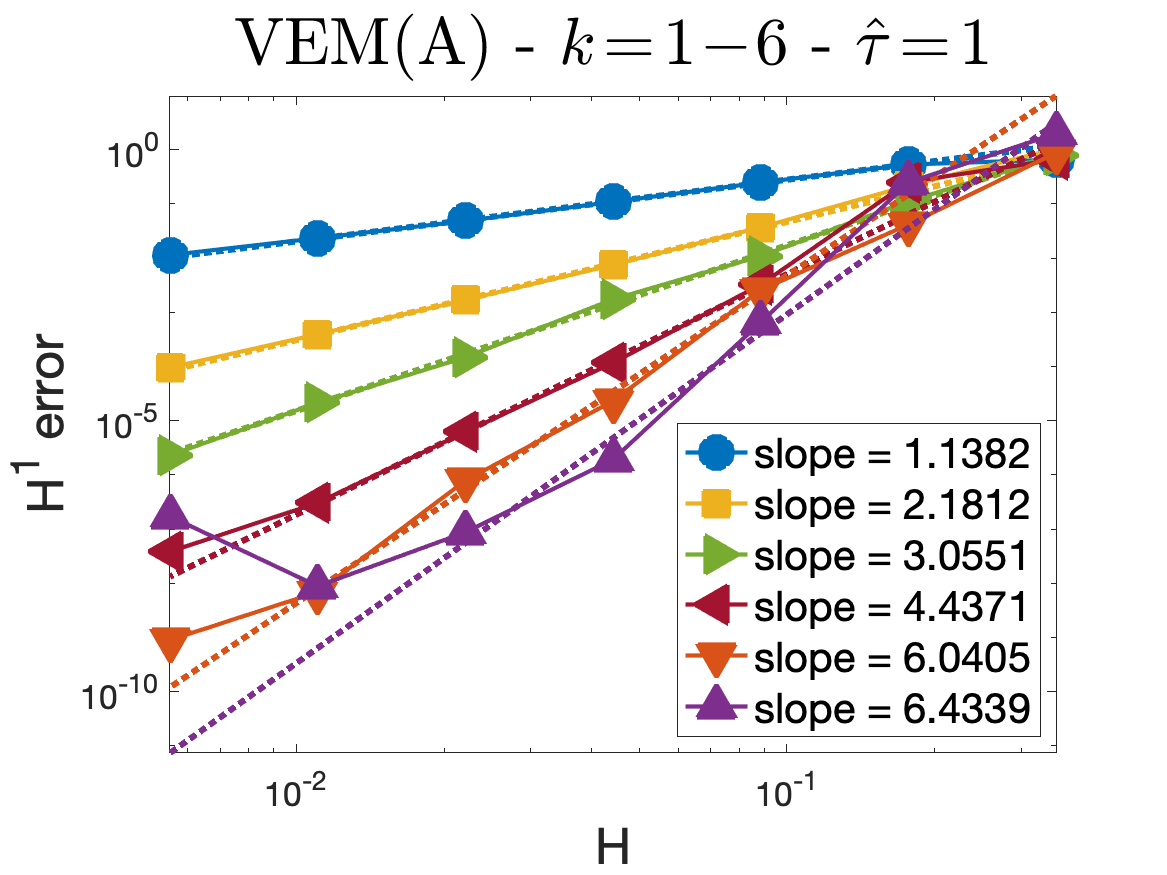} 
	\includegraphics[width=3.65cm]{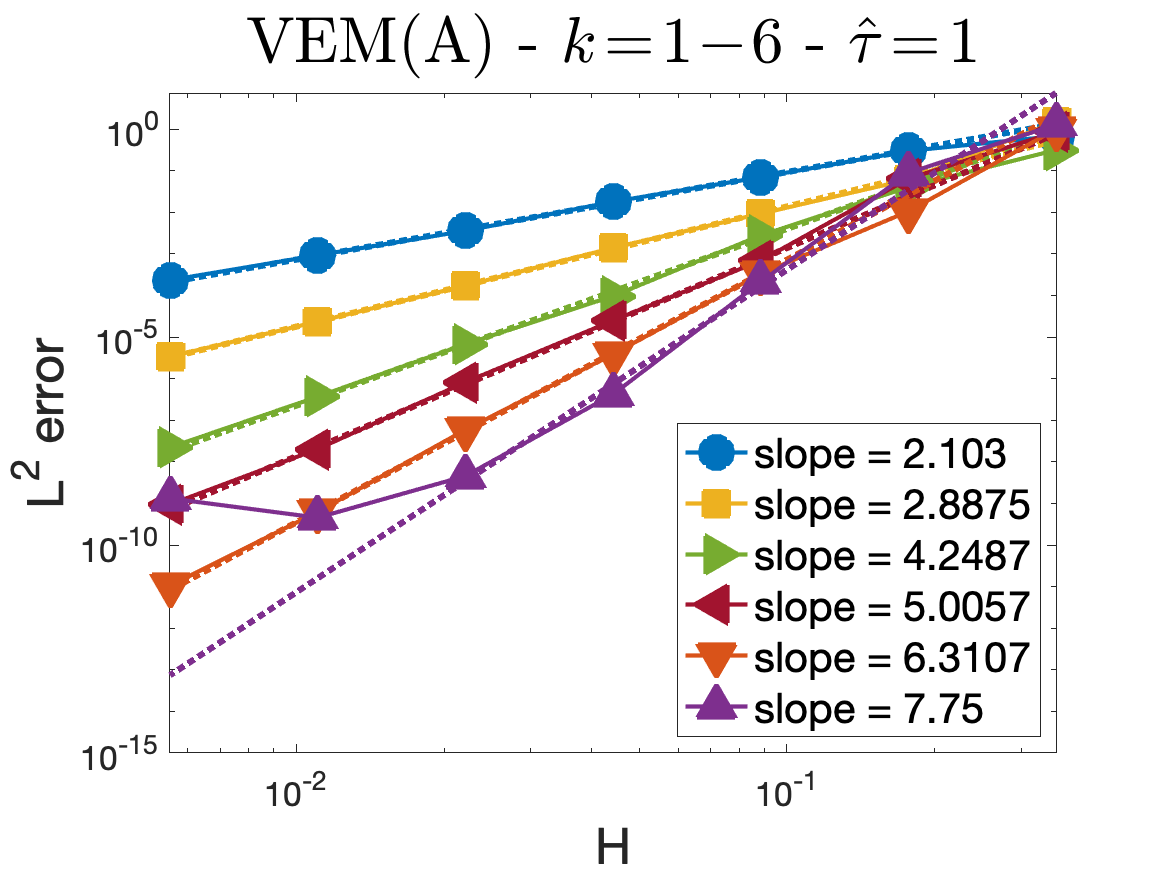} 
	\includegraphics[width=3.65cm]{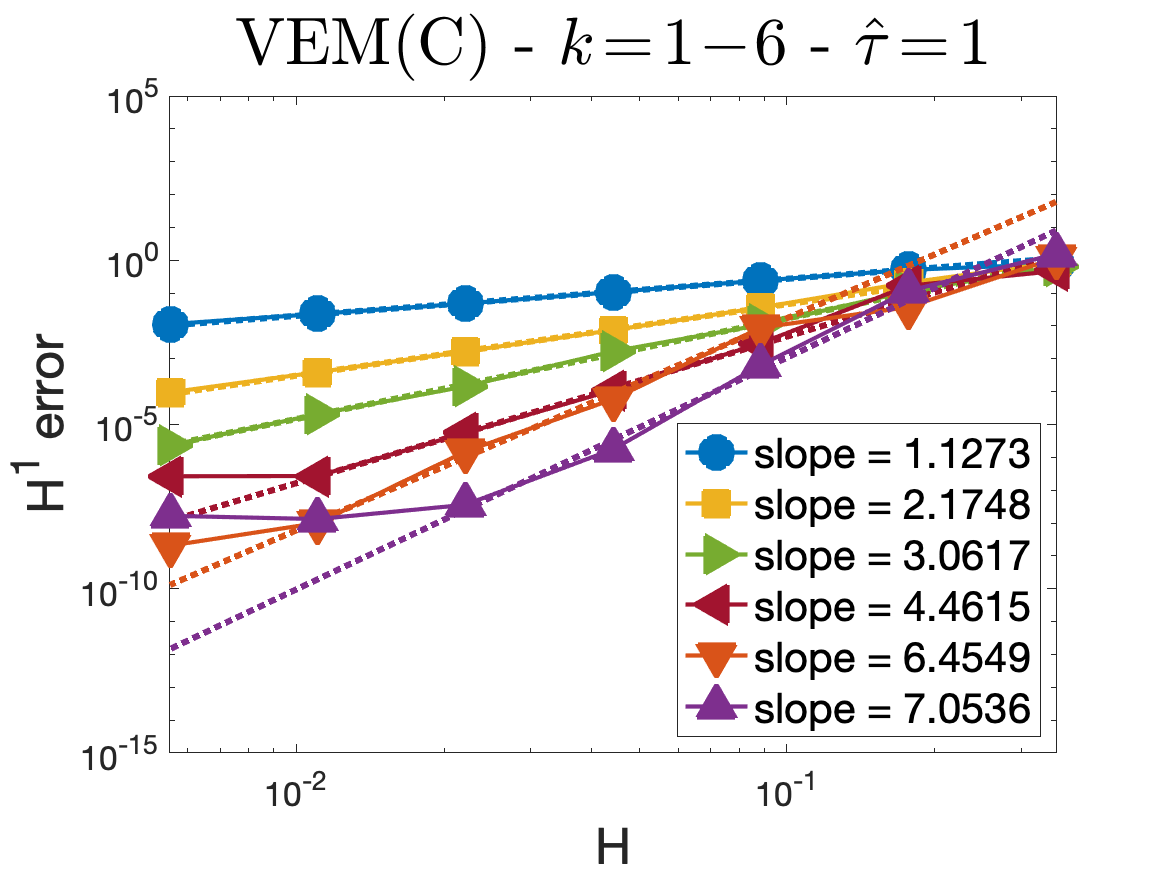} 
	\includegraphics[width=3.65cm]{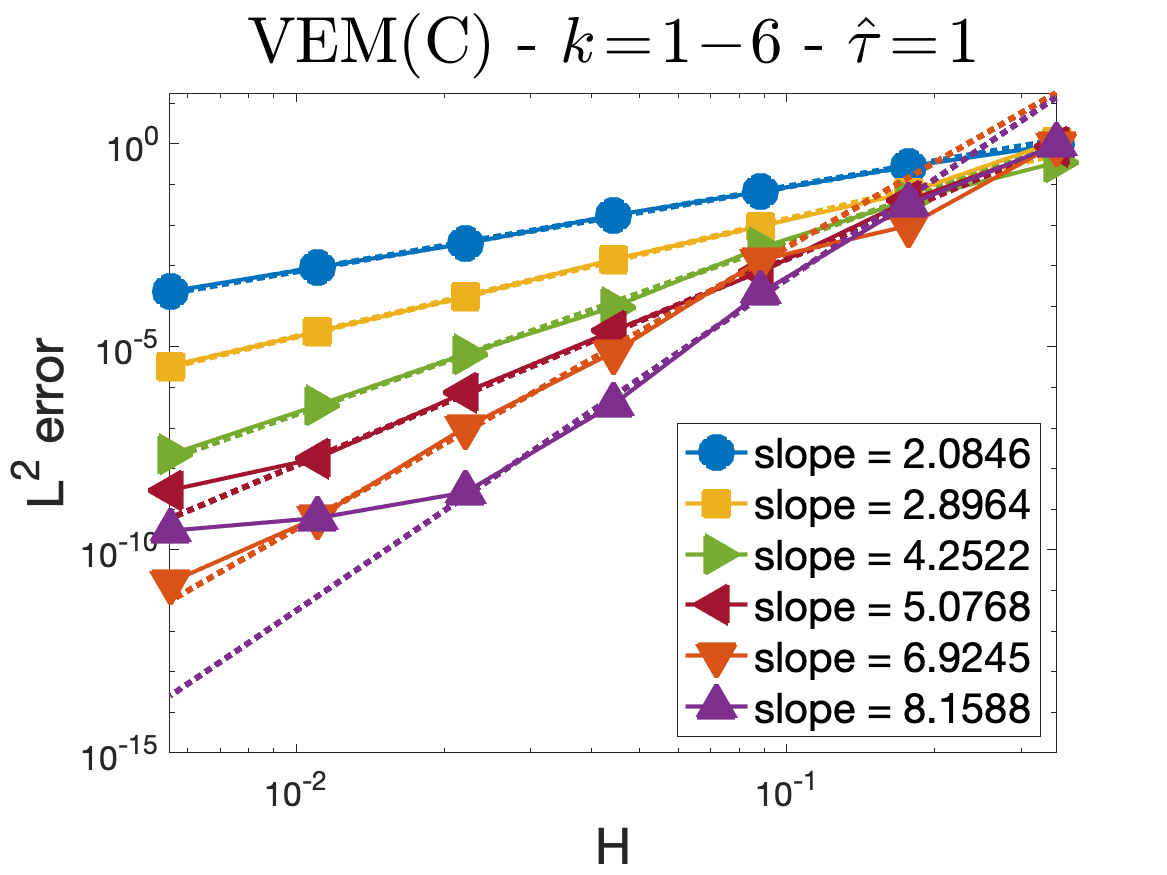} 
	\caption{Test case 1, test on standard triangular and quadrangular meshes: $H^1$  and $L^2$ error for $P_k$~(top) and $Q_k$~(middle) finite elements with boundary correction strategies \SBM~and \BHL, in comparison with VEM~(bottom) with boundary correction strategies \SBM~ and \BDT~ on the same quadrangular meshes (bottom). The stabilization parameter for Nitsche's method follows the recipe of \cite[Section 10.1]{SBMho}.}
	\label{fig:test1-g100-fem}
\end{figure}

\subsubsection{Effect of the choice of the stabilization parameter $\gamma$}
\begin{figure}[h]
	\centering
	\includegraphics[width=4.5cm]{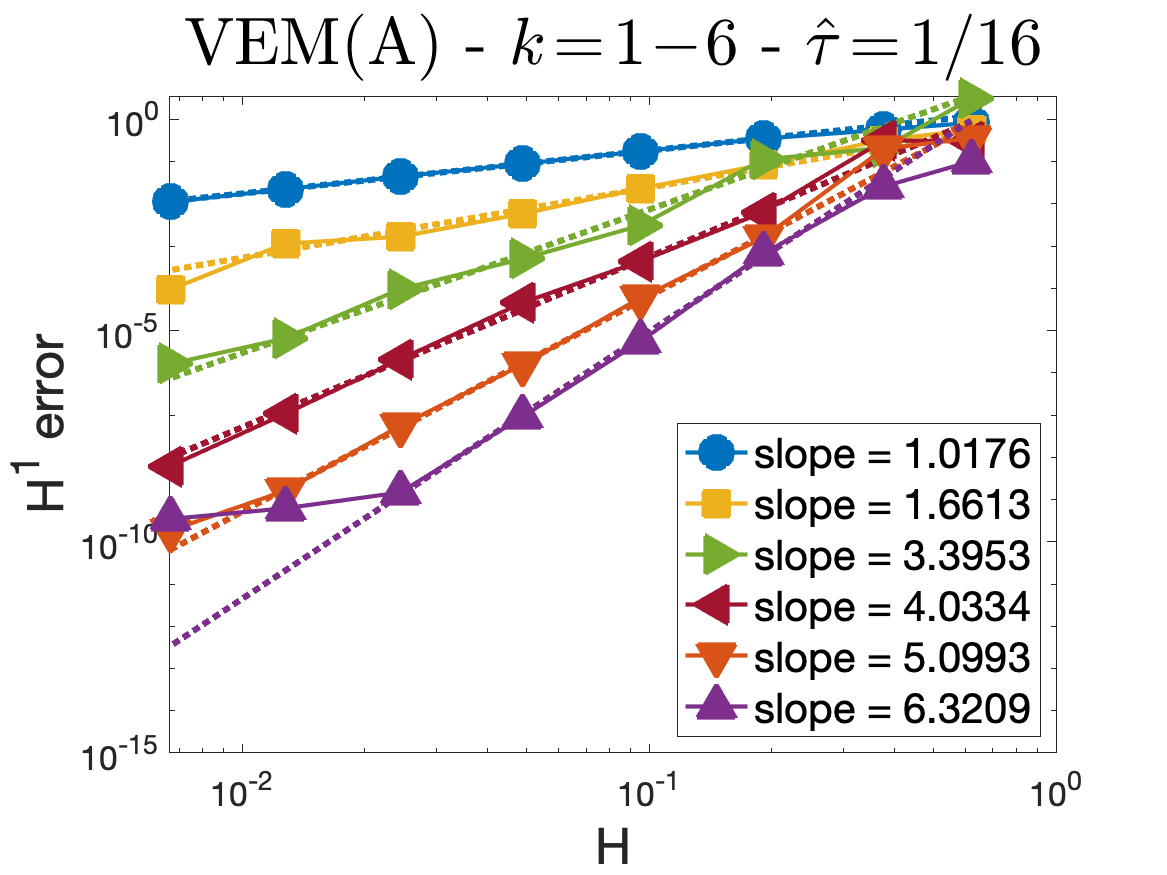}
	\includegraphics[width=4.5cm]{testVEM/figures/H1_Franke_SBM_1004.png}
	\includegraphics[width=4.5cm]{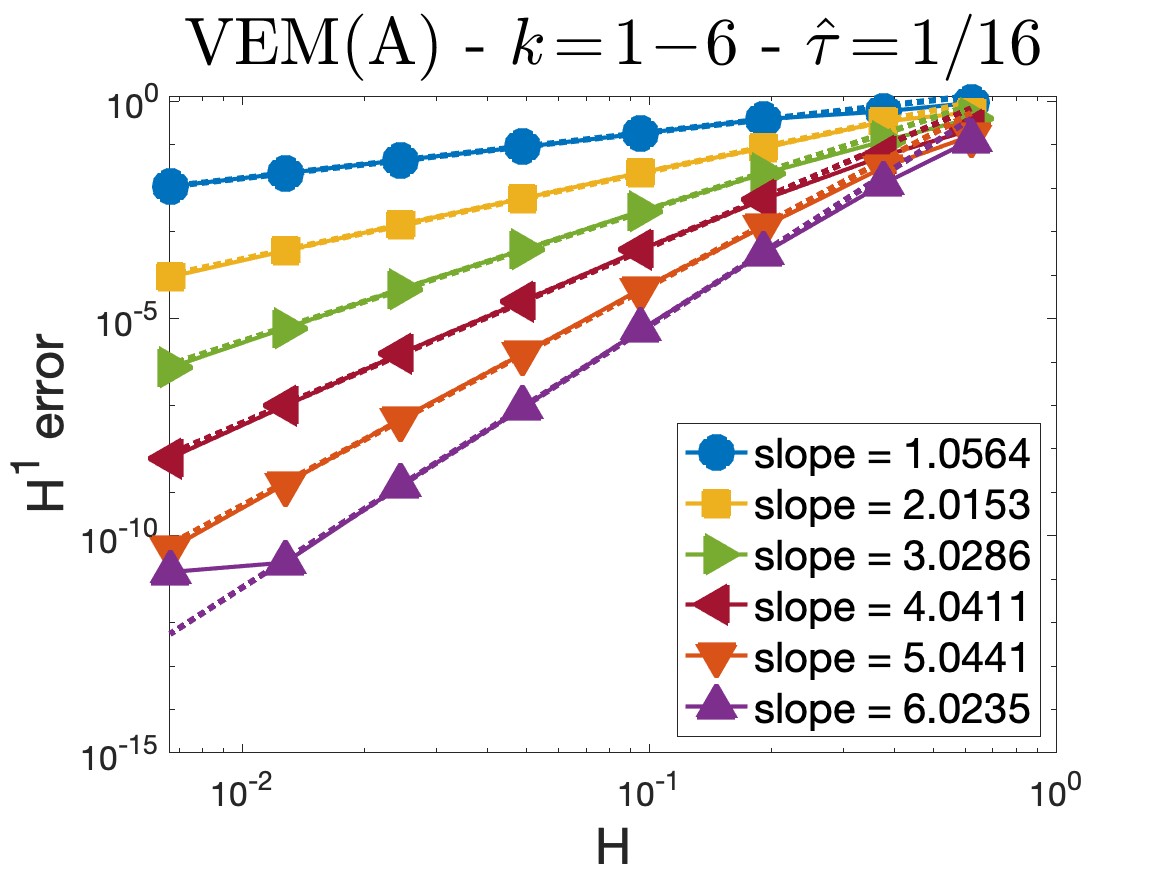}
	\\
	\includegraphics[width=4.5cm]{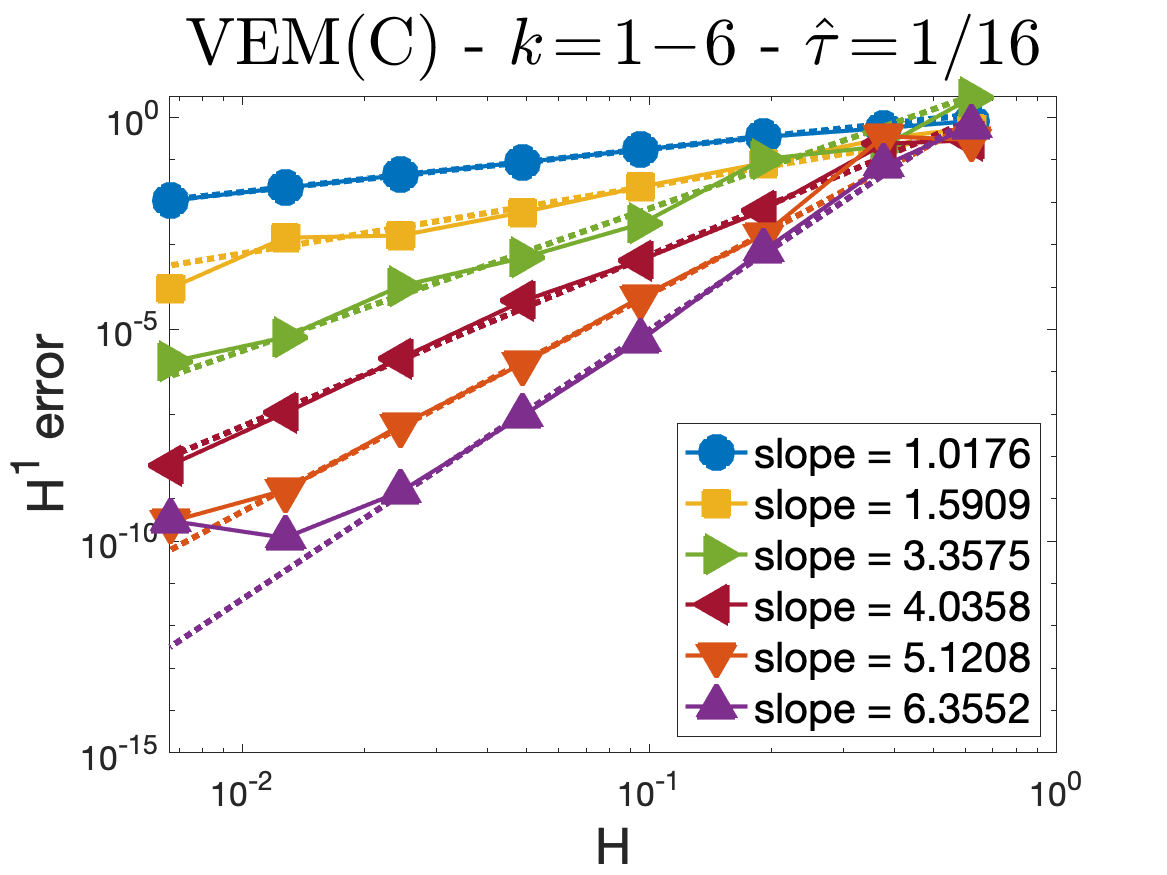}
	\includegraphics[width=4.5cm]{testVEM/figures/H1_Franke_BDT_1004.png}
	\includegraphics[width=4.5cm]{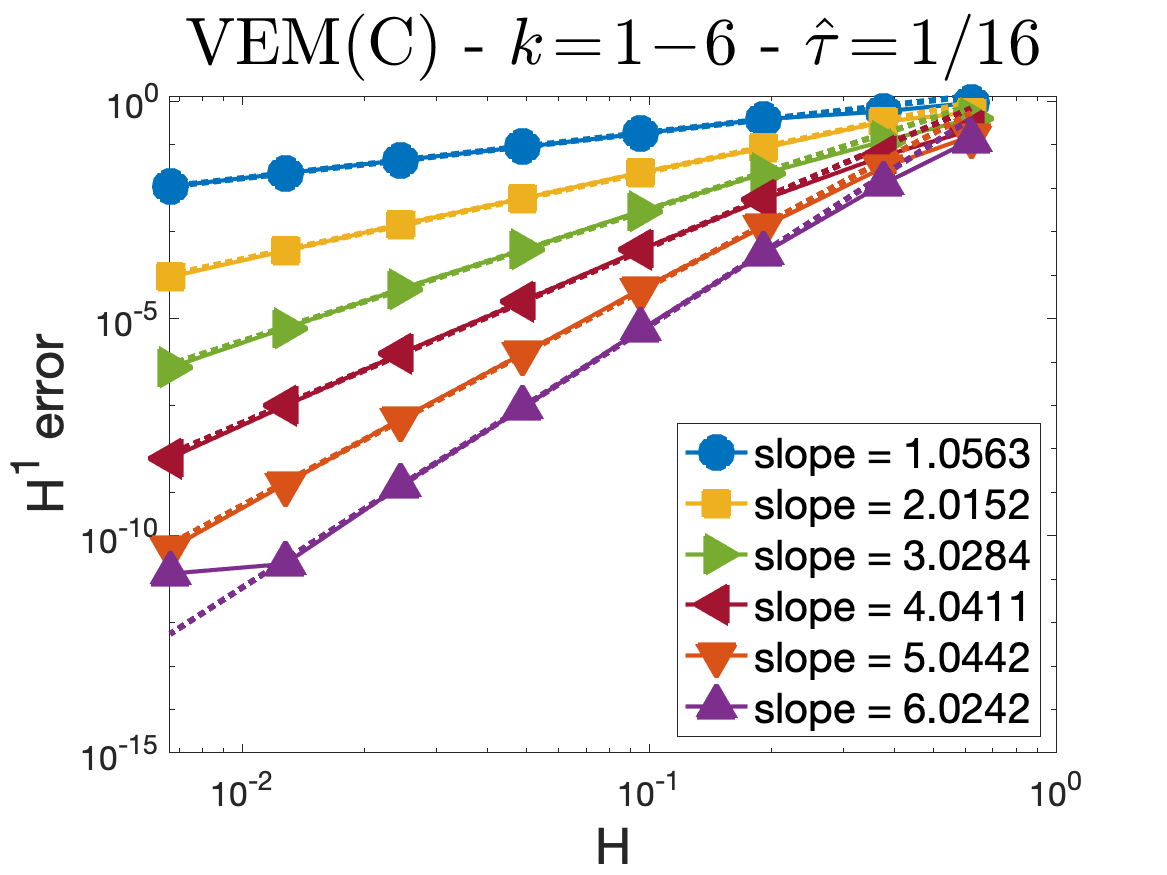}
	\caption{
		Test case 1 -  Effect of the choice of the Nitsche stabilization parameter $\gamma$. $H^1$ error for the VEM method with boundary correction strategies \SBM~(top) and \BDT~(bottom).  From left to right, $\gamma = 10$ (first  column), $\gamma = 100$ (central column), $\gamma = 150$ (last column). For all tests we have $\widehat \tau = 1/16$.}\label{fig:test1_gamma_H1}
\end{figure}

It is well known that the performance of Nitsche's method is affected by the choice of the stabilization parameter $\gamma$. The theoretical analysis of our method confirms that, also in our case, this has to be chosen large enough, that is, larger than $\gamma_0$ that depends on $k$ (tracking the dependence on $k$ in the analysis of the method suggests that $\gamma_0 \simeq \bar \gamma k^2$ for some positive $\bar \gamma$). To assess how sensitive to the choice of $\gamma$ the method actually is, we tested, for the case $\widehat \tau = 0.125$, three values of $\gamma$, namely, $\gamma = 10$, $\gamma = 100$ and $\gamma = 150$. Figure \ref{fig:test1_gamma_H1}  suggests that $\gamma = 10$ is too low for $k > 2$, while $\gamma = 100$ is good up to $k = 5$, but too low for $k=6$, for which $\gamma = 150$ yields instead, optimal results for both strategies \SBM~and \BDT. We also observe that increasing $\gamma$ does not seem to negatively affect the error, and then it seems reasonable, in the absence of a more precise analysis on the dependence of such a constant on the different parameters of the method
 to err on the side of caution and choose $\gamma$ larger rather than smaller.

{  \begin{remark}
		The virtual element method can be also combined with the penalty free shifted boundary method proposed in \cite{collins2023penalty}, and we believe that, provided the VEM stabilization constant $\cstar$ in condition \eqref{defSK} is large enough, the analysis presented here carries over to such a formulation.
\end{remark}}

\subsubsection{Robustness with respect to the parameter $\widehat\tau$} 
\begin{figure}
	\centering
	\includegraphics[width=4.9cm]{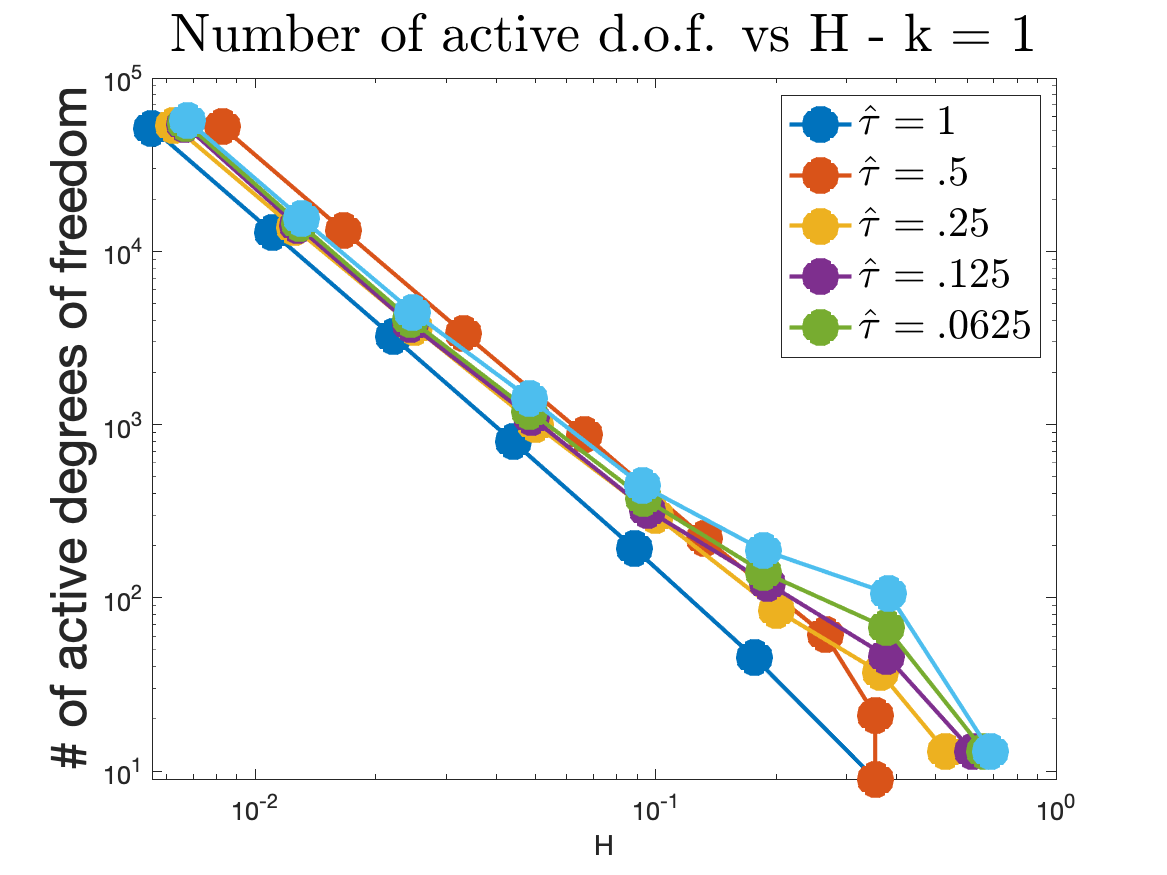}
	\includegraphics[width=4.9cm]{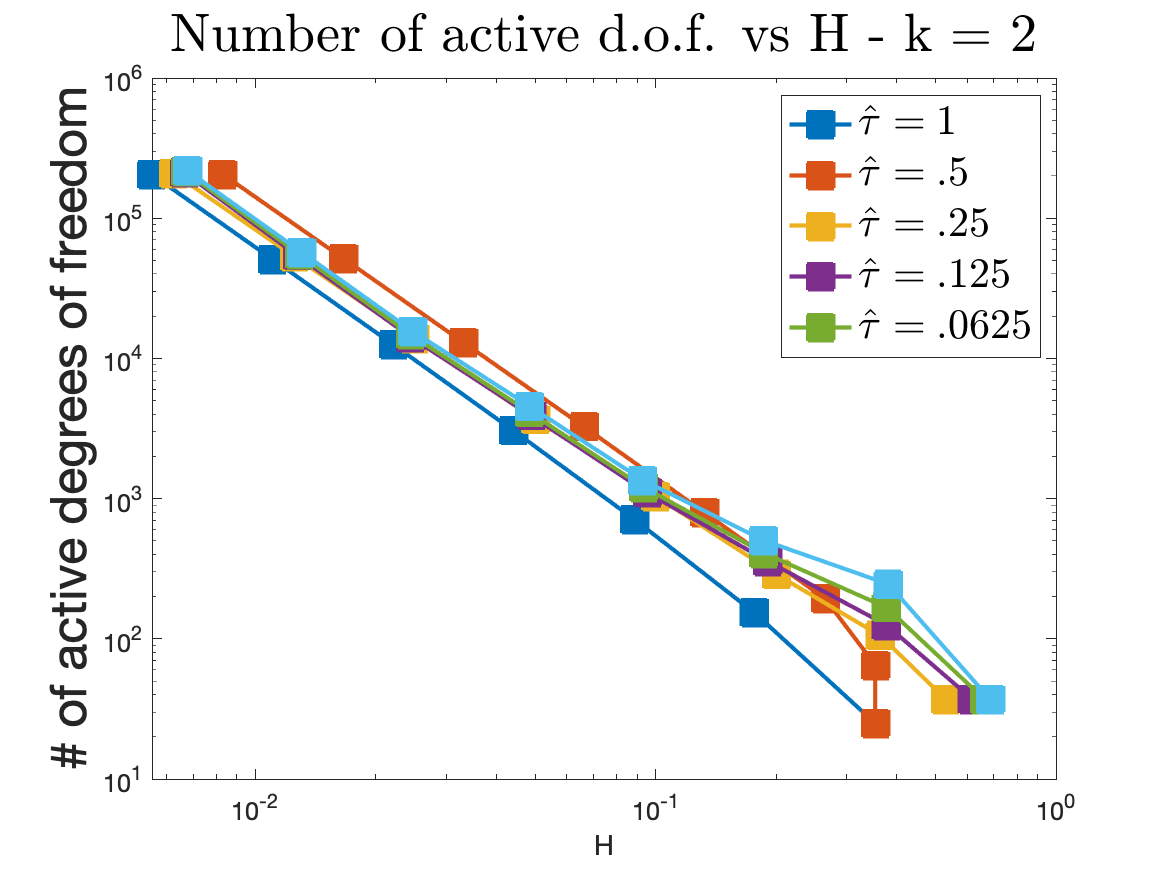}
	\includegraphics[width=4.9cm]{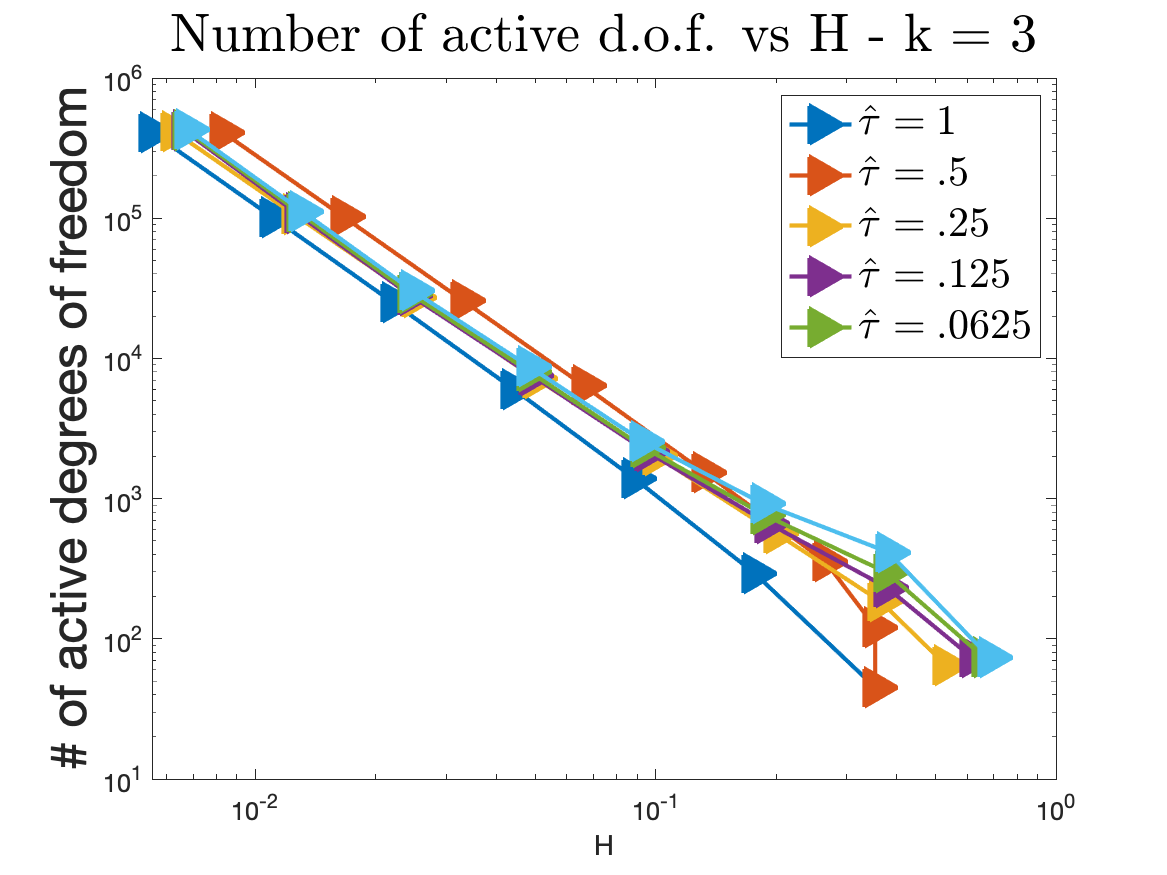}\\
	\includegraphics[width=4.9cm]{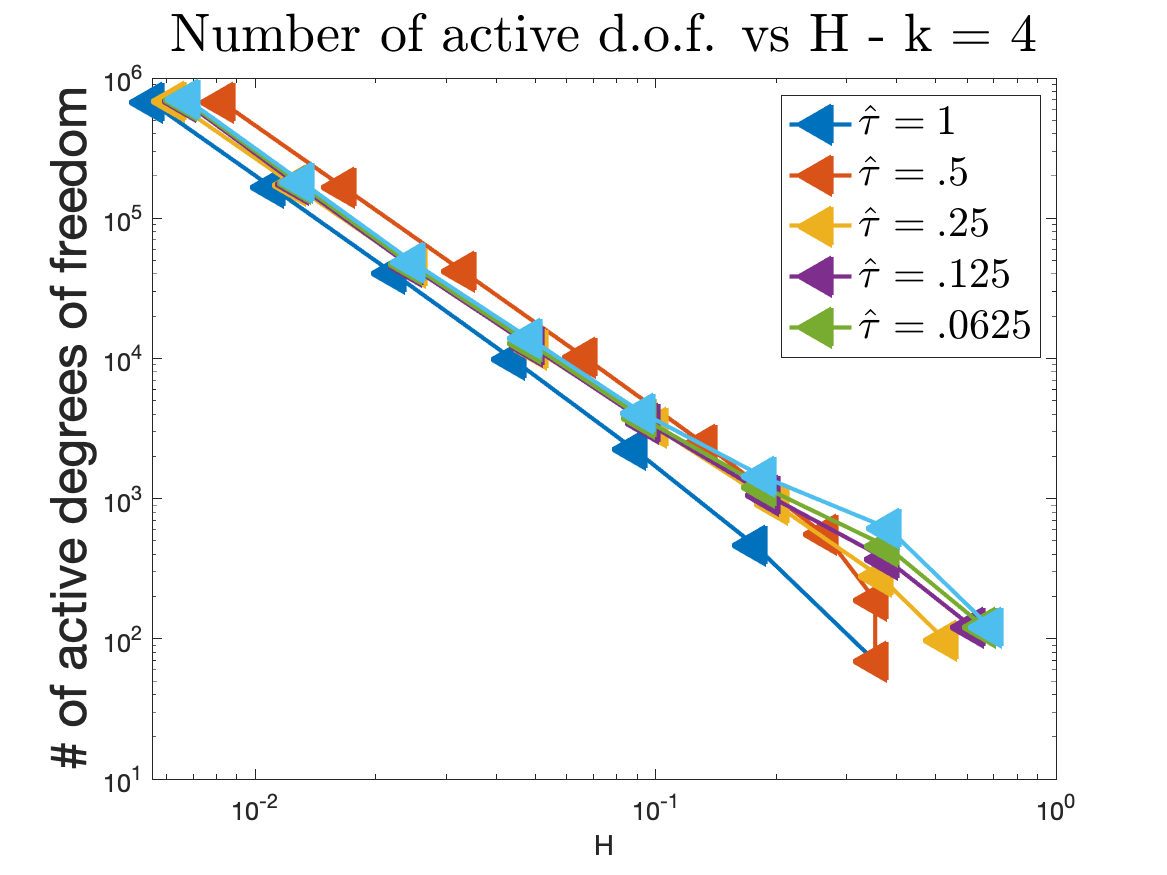}
	\includegraphics[width=4.9cm]{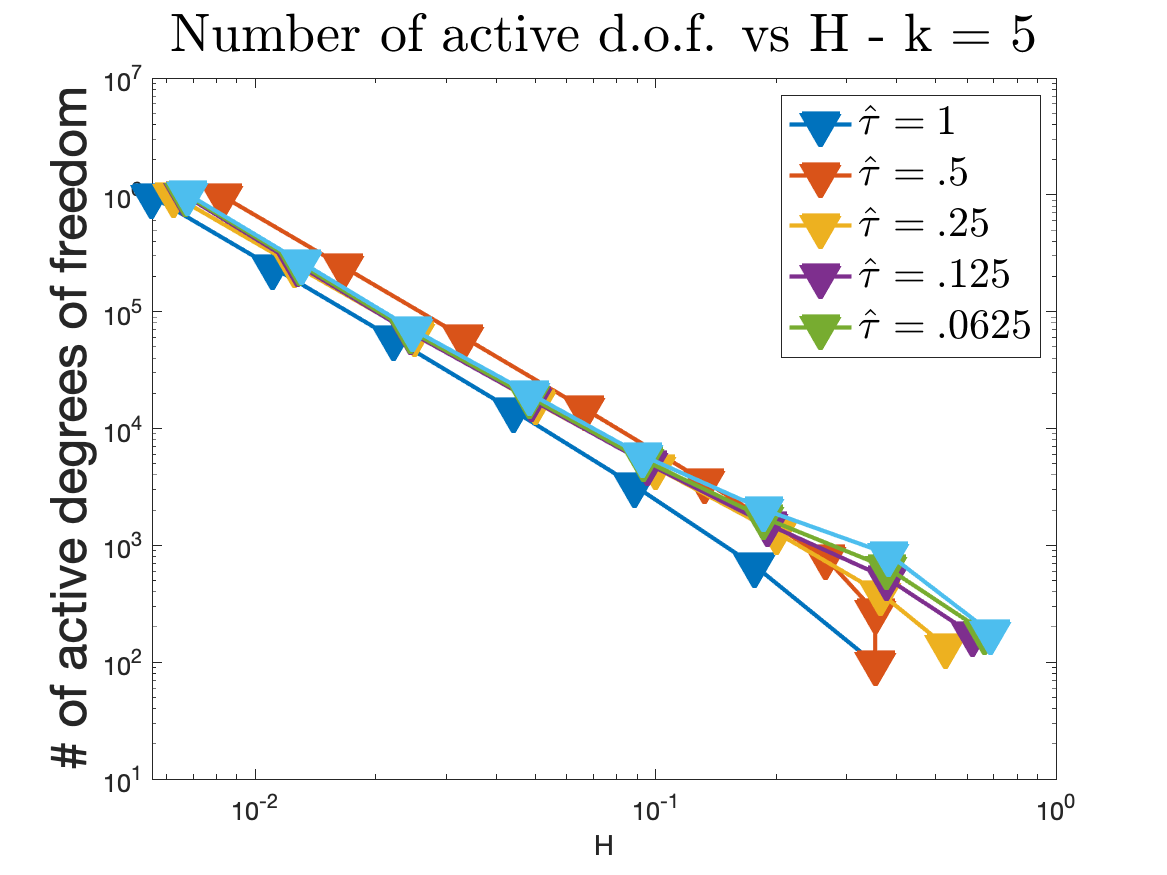}
	\includegraphics[width=4.9cm]{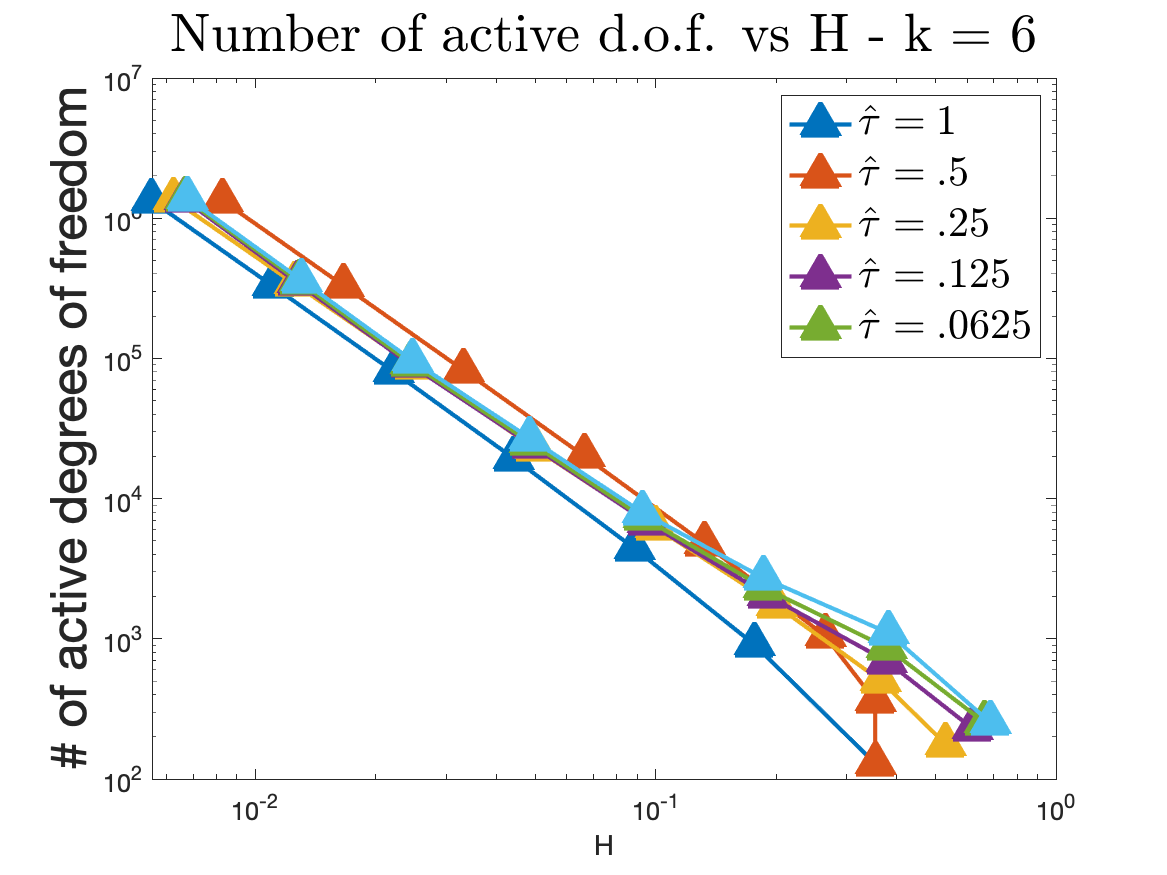}
	\caption{Number of active degrees of freedom  (independent of the chosen boundary correction strategy) for the meshes used for the tests on the disk domain, for different values of the order $k$ of the method and of the refinement parameter $\tau$.}\label{fig:robmesh}
\end{figure}

We aim at demonstrating that, by applying the static condensation procedure proposed in Section \ref{sec:Lazy}, it is possible to lower the value of the parameter $\widehat\tau$ without asymptotically increasing the number of degrees of freedom. In Figure \ref{fig:robmesh} we plot, for $k=1,\cdots,6$ the number of active degrees of freedoms retained after the elimination of the lazy degrees of freedom by the approach presented in Section \ref{sec:Lazy}, for the meshes used for all the tests of Sections \ref{sec:test1} and \ref{sec:test2}. Observe how, for $H$ sufficiently small and $\widehat\tau \leq .25$ the curves are, for all values of $k$ essentially superposed, that is, the number of degrees of freedom for a given value of $H$ does not increase as $\widehat\tau$ goes to $0$. For all $k$ it is therefore possible to choose such a parameter in such a way that the resulting method is stable, without increasing the size of the linear system to be solved.

{ 
	Next, we study the overhead resulting from the static condensation procedure. In Table  \ref{tab:3}  we report the ratio between the computational cost needed for carrying out the elimination of the lazy degrees of freedom and the one needed for the assembly of the stiffness matrix. We see that for  larger values of $H$ and/or  $\widehat \tau$, the overhead is negligible, while for smaller values of simultaneously $H$ and $\widehat \tau$ the overhead becomes larger, though still acceptable. We recall that both assembly of the stiffness matrix and elimination of the lazy degrees of freedom are highly parallelizable, as they can be carried out independently of each element and macro edge respectively.

We then consider the condition number  of the stiffness matrices resulting from our method. We test the matrices before and after the elimination of the lazy degrees of freedom, and we report, for the sake fof comparison, the matrix relative to Nitsche's method without boundary correction.  As we can see in Figure \ref{fig:cond}, the condition number with and without elimination of the lazy degrees of freedom are prectically superposed. The instability of the method for larger values of $k$ and $\widehat \tau$ is confirmed, while for smaller values of $\widehat \tau$ the condition of the boundary corrected method is smaller than the one for the plain Nitsche's method.

\begin{figure}
	\includegraphics[width=4.9cm]{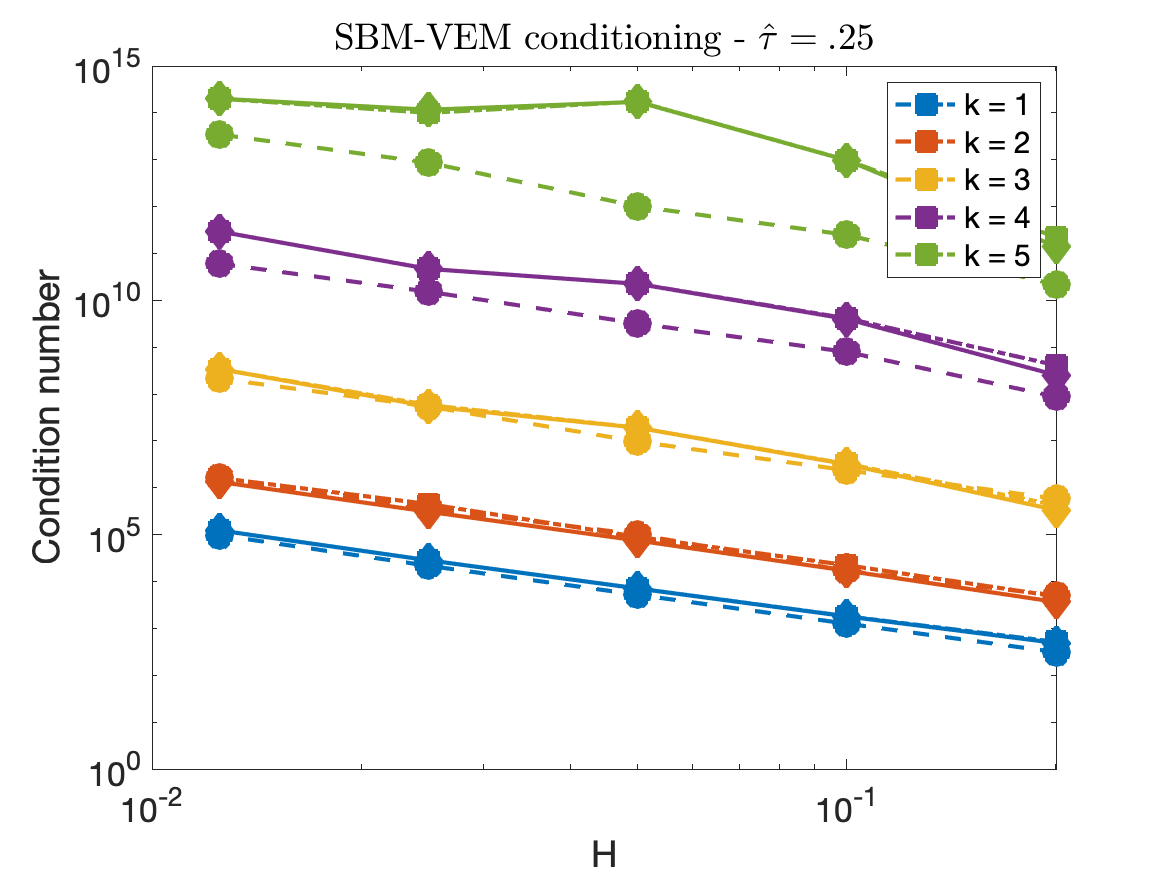}
	\includegraphics[width=4.9cm]{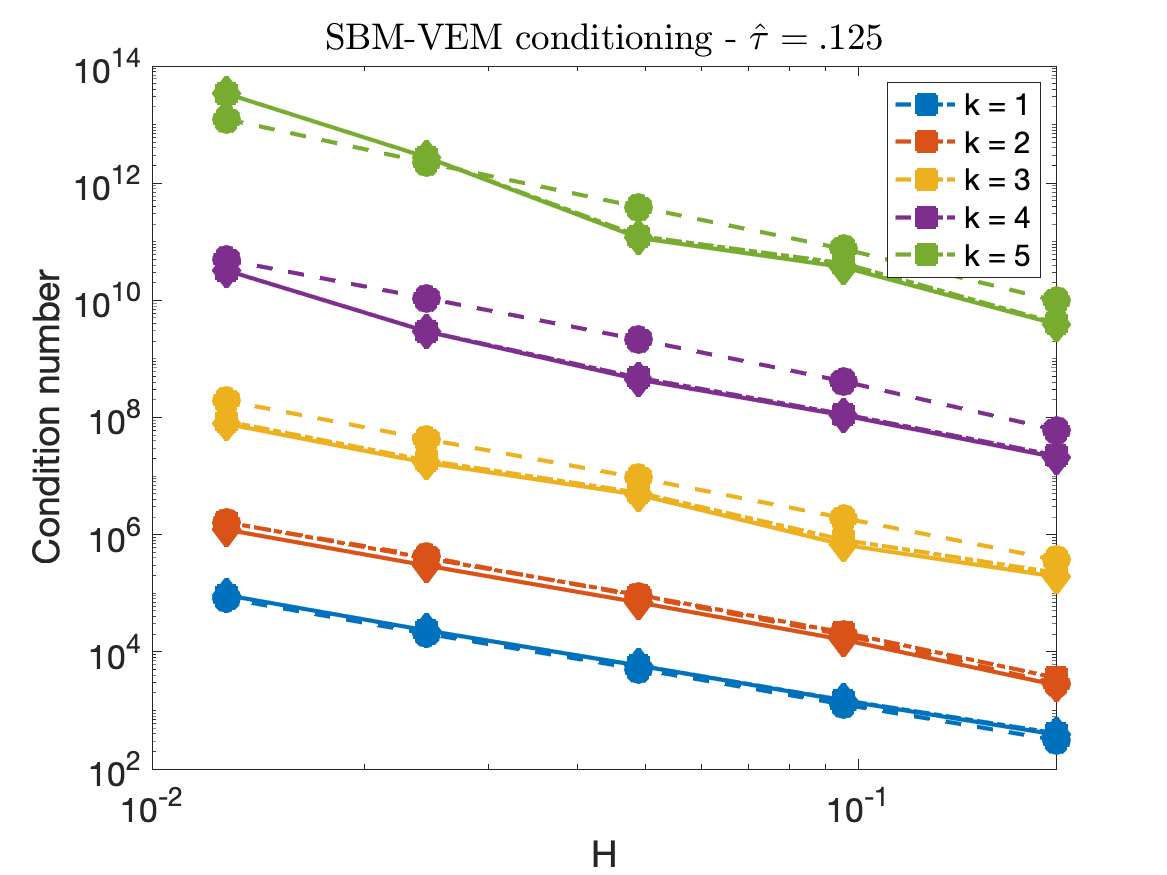}
	\includegraphics[width=4.9cm]{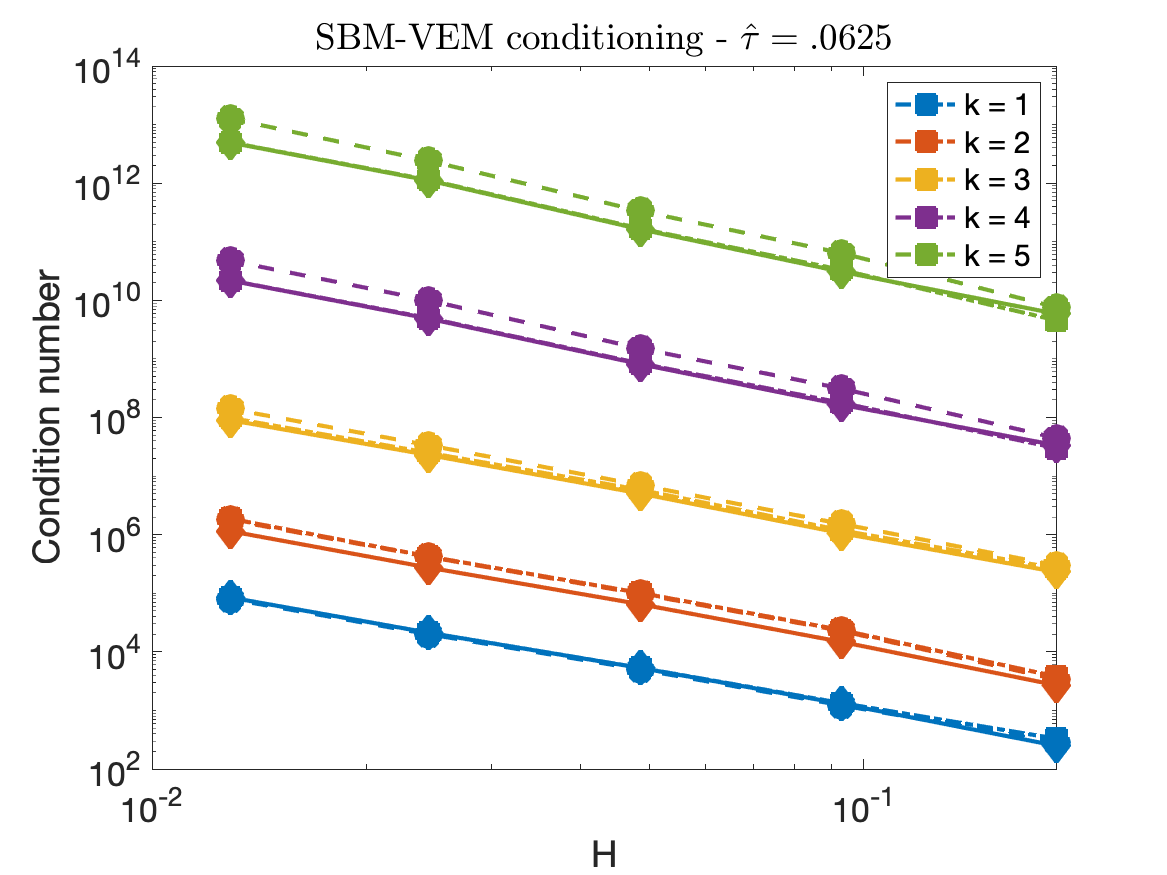}
	\caption{Condition number of the matrices resulting from the VEM\SBM~method, for different values of the parameter $\widehat \tau$ and of the order $k$ of the method, before (dash-dot line) and after (solid line) elimination of the lazy degrees of freedom. The condition number for the matrix relative to Nitsche's method (dashed line) are plotted for comparison.}\label{fig:cond}
\end{figure}
}

\subsection{Test 2 -- Bean domain}\label{sec:test2}

In this test we consider a {  non convex} domain with a corner with interior angle equals to $3\pi/2$, obtained as the union of a quarter of a disk with two half disks: 
\[
\Omega = \Omega_1 \cup \Omega_2 \cup \Omega_3\quad\text{ with }  \begin{cases}
\Omega_1 =	\{(x,y): (x^2 + y^2) < 1,\ x<0, y > 0\}, & \\
\Omega_2 =	\{(x,y): (x+.5)^2 + y^2 < 0.25, x<0\}, & \\
\Omega_3 =	\{(x,y): x + (y-.5)^2 < 0.25, x<0\}. & 
\end{cases}
\]
The load and the non homogeneous Dirichlet data are chosen so that  the solution (see Figure \ref{fig:solution}), expressed in polar coordinates, is 
\[
u(\rho,\theta) = \rho^{2/3} \sin(2 \theta/3), 
\]
a standard benchmark for corner singularities. We know that $u \in H^s(\Omega)$ for all $s < 3/2$, but $u \not \in H^{3/2}(\Omega)$. As we are in the presence of an interior angle, we follow an $hp$ strategy by resorting to geometric meshes, progressively refined in the proximity of the singular point $(0,0)$, while simultaneously increasing the order of the method. In the case of a polygonal domain, this strategy is expected to yield an exponential convergence of the $H^1(\Oh)$ error satisfying the bound
\[
\| u - u_h \|_{1,\Oh} \lesssim \exp(-\sqrt[3]{N}),
\]
where $N$ is the number of degrees of freedom \cite{beirao_hp_exponential}.  Observing that the parameter $\tau_\gamma$ and $\gamma_0$ given by Theorem \ref{cor:nitsche:curved} depend on $k$ and that a rough analysis suggests that $\tau_\gamma \simeq 1/k^2$ and $\gamma_0 \simeq k^2$, as the order of the method increases, we decrease $h$ so that the first condition is satisfied, and possibly modify also the elements that are not refined, enlarging them through the union of squared elements of the finer squared mesh with meshsize $h$. We refer to this adjustment of the elements close to the boundary as $\delta$ refinement. Simultaneously, we adjust the stabilization parameter and set $\gamma = \bar \gamma k^2$. 
 
\begin{figure}
	\centering
	\includegraphics[width=4.9cm]{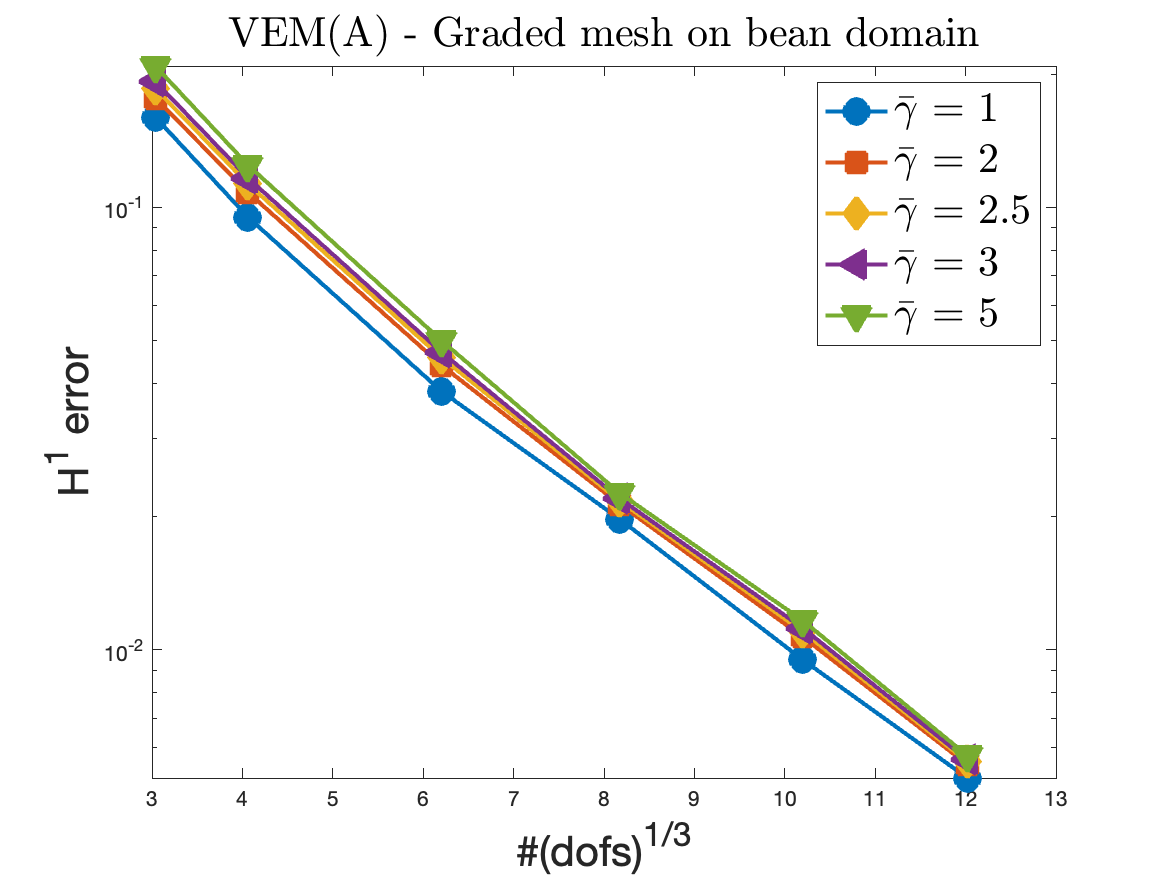}
	\includegraphics[width=4.9cm]{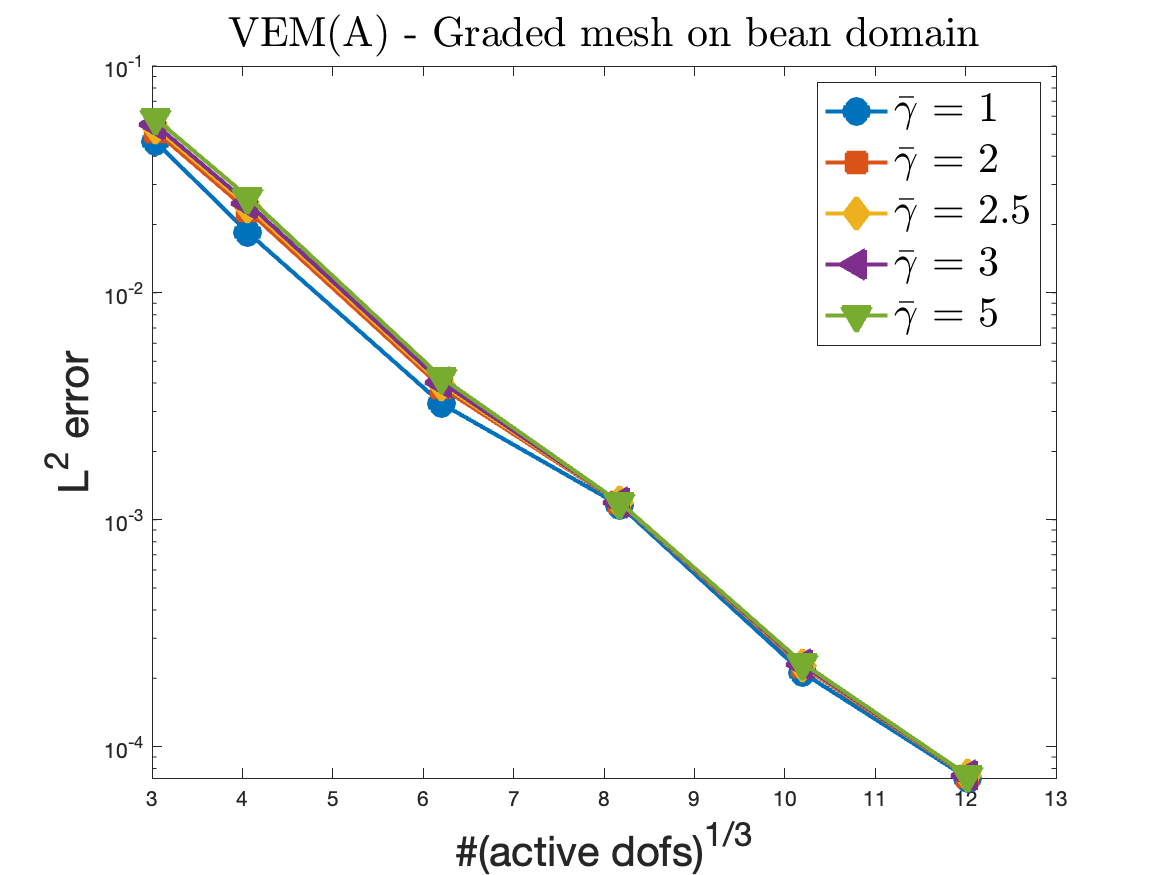}
	\\
	\includegraphics[width=4.9cm]{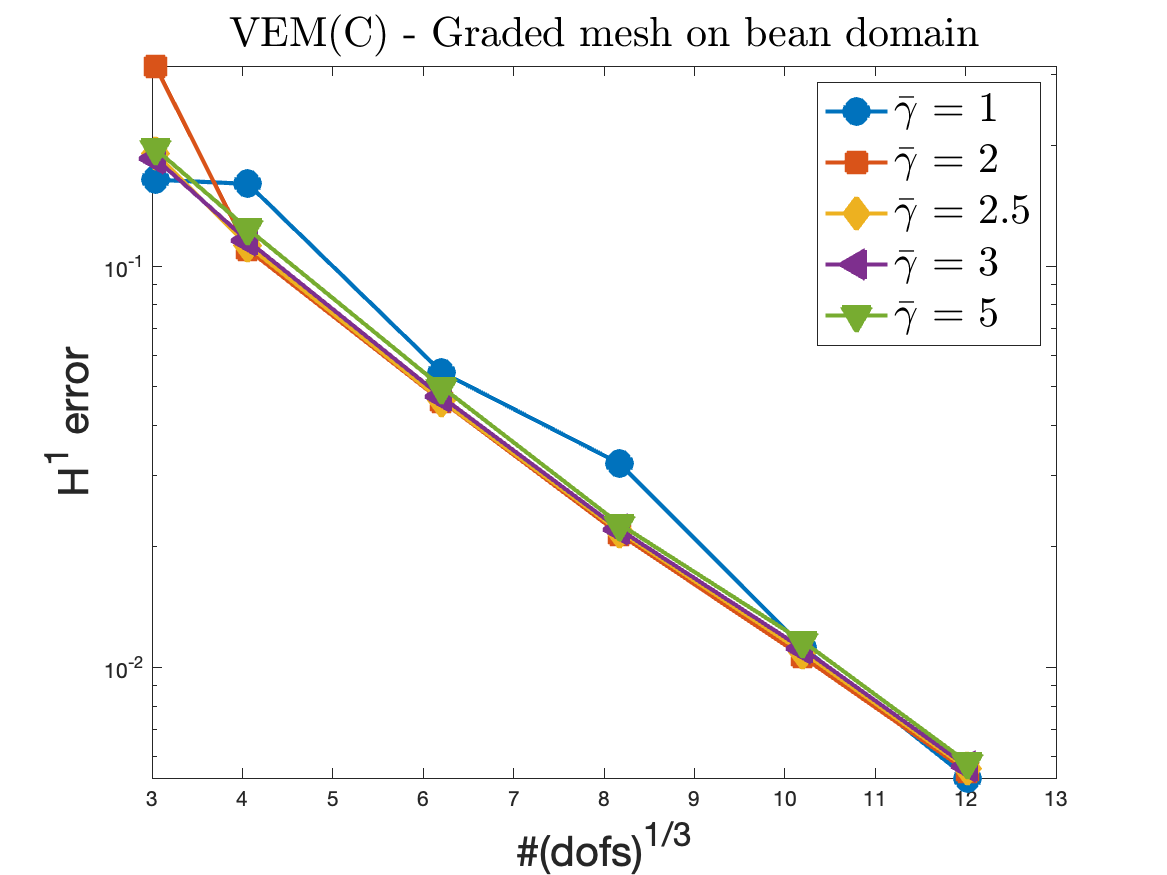}
\includegraphics[width=4.9cm]{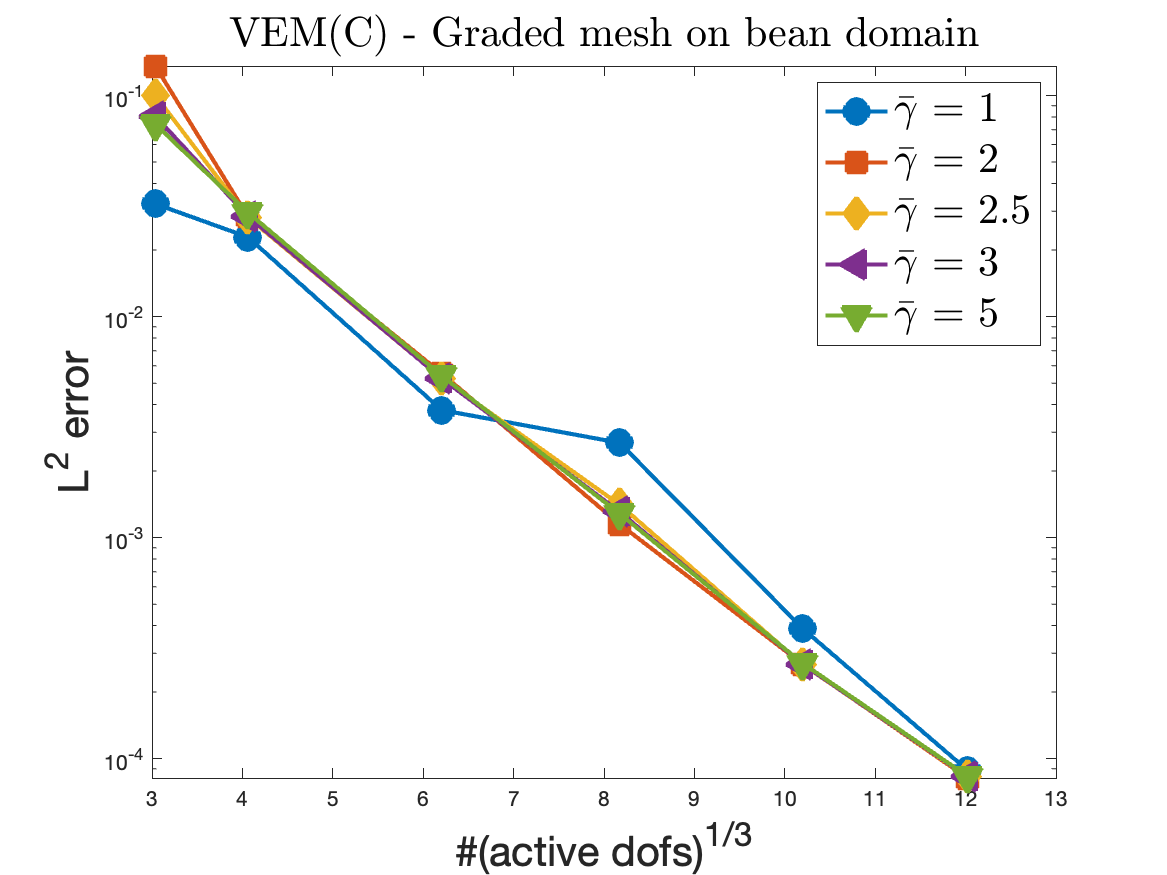}
\caption{Test case 2. Convegence of the VEM method with boundary correction strategies \SBM~and \BDT, for a singular solution, with a graded mesh and both $hp$ and $\delta$ refinement.}\label{fig:test2}. 
\end{figure}

Once again, we consider the virtual element method with both strategy \SBM~ and \BDT~ boundary correction, and we test different values of the parameter $\bar \gamma$. The results, presented in Figure \ref{fig:test2}, display a behaviour similar to the one obtained for polygonal domains, where no boundary correction is needed. This suggests that, in the framework of the virtual element method, that allows to adjust the distance of the approximate and true boundary by reducing the mesh size of the fine grid $\microTess$,  boundary correction approaches such as the shifted boundary method or the BDT method with closest point mapping are potentially applicable in the $hp$ framework.

\renewcommand{\vec}[1]{\mathbf{#1}}

\subsection{Application to elasticity}
\begin{figure}
	\centering
	\includegraphics[width=4.9cm]{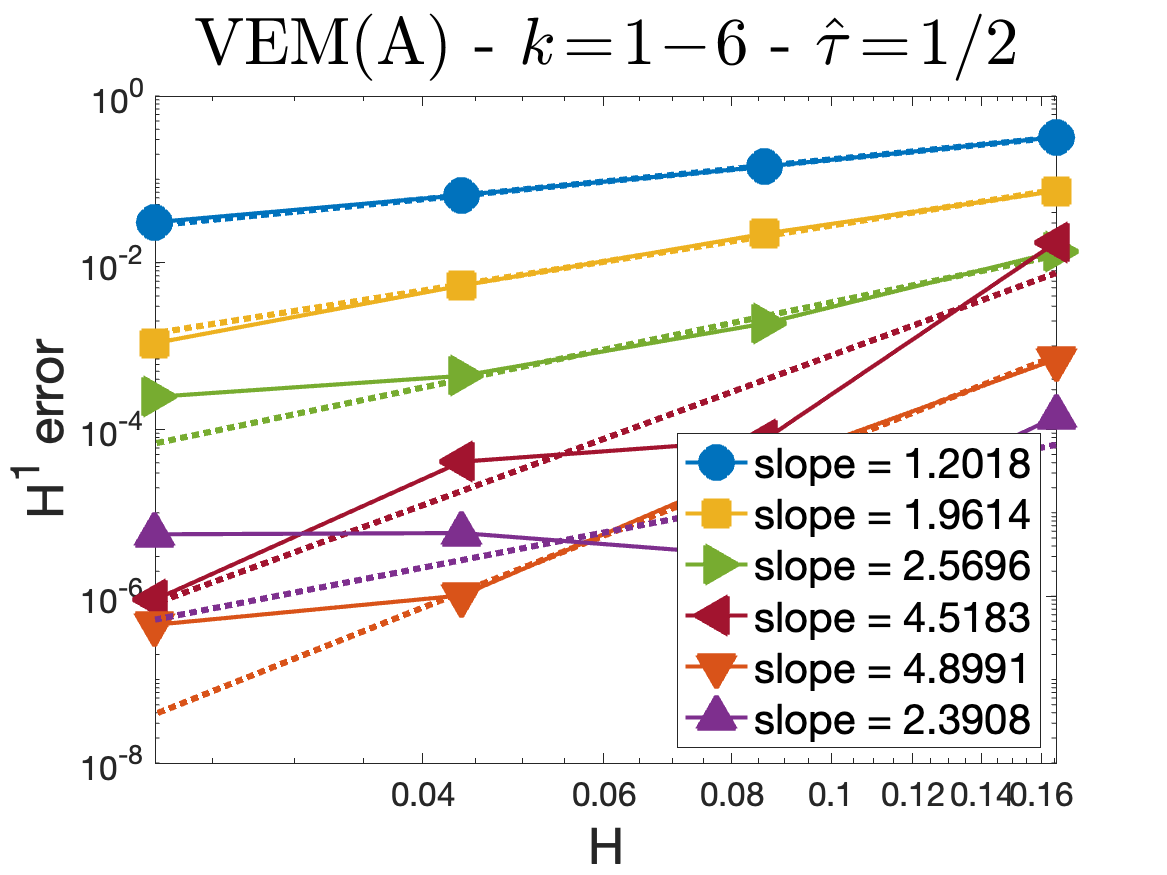}
	\includegraphics[width=4.9cm]{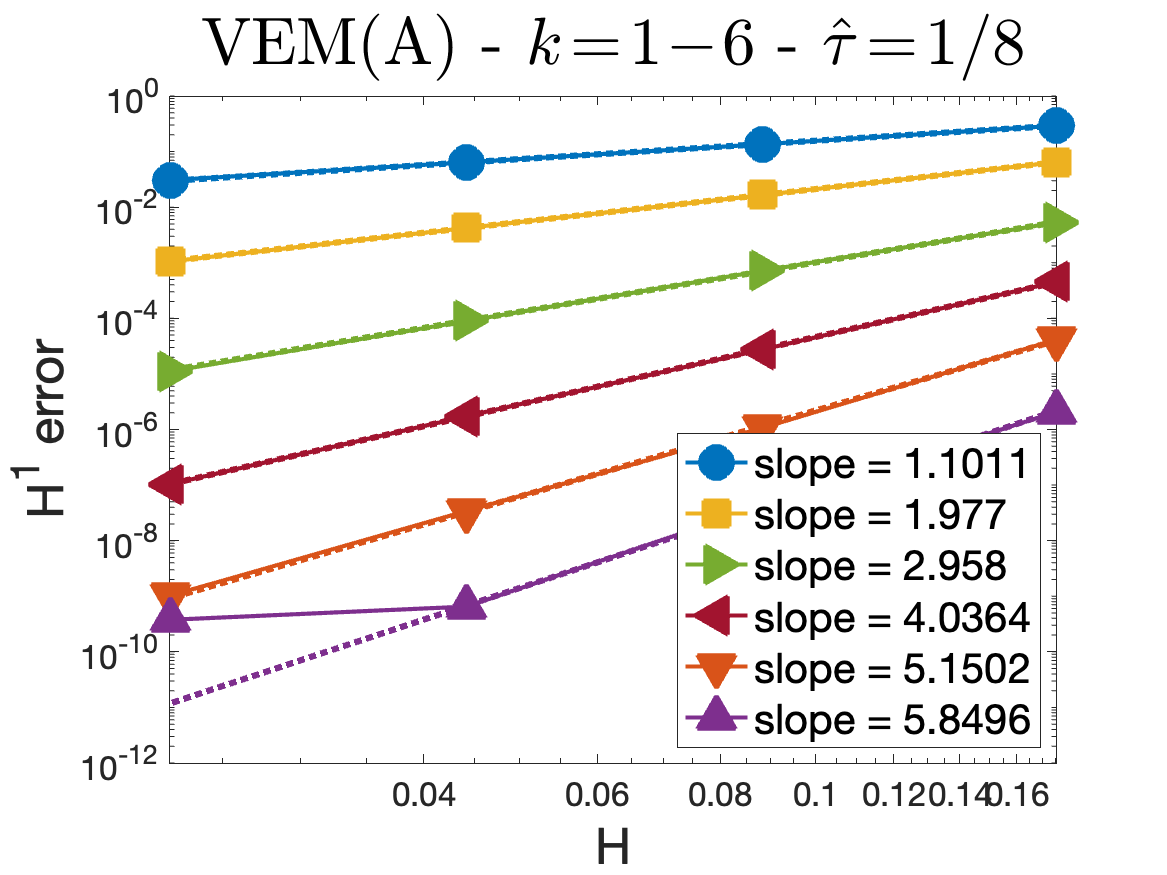}
	\includegraphics[width=4.9cm]{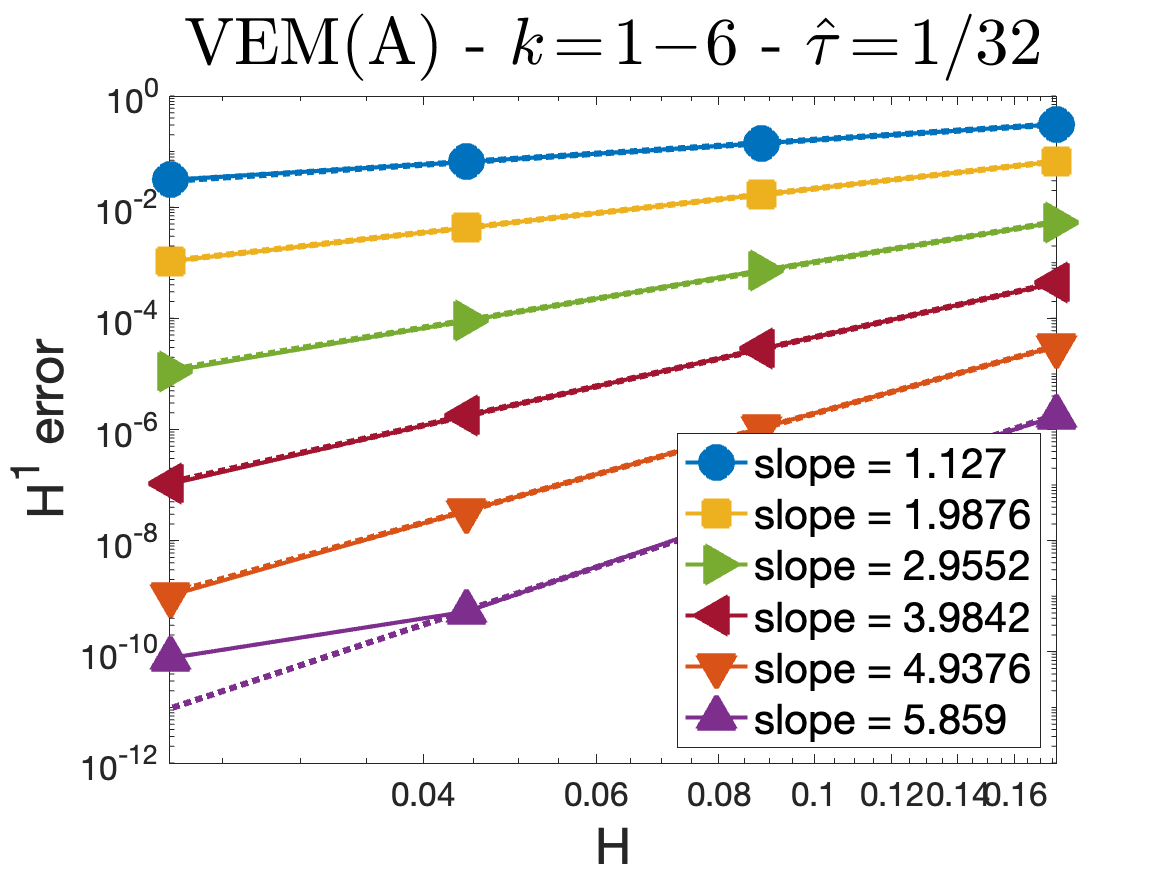}
	\caption{Convergence of the VEM\SBM~method on equation \eqref{eq:elasticity}.}\label{fig:elasticity}
\end{figure}

As it happens in the finite element framework, also in the virtual element framework boundary correction methods such as the ones considered in this paper can be extended to a wide variety of different problems. We test here the proposed approach on the following elasticity model problem:
\begin{equation}\label{eq:elasticity}
	-\mathrm{div}(\sigma) = \vec f\quad\text{in }\Omega, \qquad
	\vec u = \vec g \quad\text{on }\partial\Omega
\end{equation}
with $\sigma\colon\Omega\to\mathbf R^{2\times 2}$ is the stress tensor defined as
\[
\sigma = 2\mu\varepsilon(\vec u) + \lambda \mathrm{div}(\vec u)I,
\]
where $\mu$ and $\lambda$ are the constant Lam\'e parameters. To discretize such a problem, we take each component of the vector valued displacement space in the usual non-enhanced virtual space for the Poisson problem~\cite{Artioli2017}. The discrete problem then reads: find $\vec{u}_h \in (V_h)^2$ such that for all $\vec{v}  \in (V_h)^2$ it holds that

\begin{multline*}
	2\mu\int_{\Tess} \Pi^0_{k-1}\varepsilon(\vec u_h)\colon\Pi^0_{k-1}\varepsilon(\vec v) + \lambda\int_{\Tess} \Pi^0_{k-1}\mathrm{div}(\vec u_h)\cdot\Pi^0_{k-1}\text{div}(\vec v) + S_K(\vec u_h, \vec v)\\
	-\int_\bOh\left[\left(
	2\mu\varepsilon(\Pi^\nabla\vec u_h) + \lambda \mathrm{div}(\Pi^\nabla \vec u_h)\right) 
	\nuhx 
	\right] \cdot\vec v\\
	-\int_\bOh \left[\left(2\mu\varepsilon(\Pi^\nabla\vec v) + \lambda \text{div}(\Pi^\nabla \vec v)\right)\vec \nuhx 
	+ \gamma H^{-1}   (\Pi^\nabla\vec v+\widehat{\mathscr C}[\Pi^\nabla \vec v])
	\right] \cdot(\Pi^\nabla\vec u_h+\mathscr C[\Pi^\nabla \vec u_h])\\
	=
	\int_{\Tess}\vec f_h\cdot\vec v
	-\int_{\bOh} \left(\left(2\mu\varepsilon(\Pi^\nabla\vec v) + \lambda \text{div}(\Pi^\nabla \vec v)\right)\vec n
	+ \gamma H^{-1} (\Pi^\nabla\vec v+\widehat{\mathscr C}[\Pi^\nabla \vec v])
	\right)\cdot\vec g^*.
\end{multline*}
We tested our method in the VEM\SBM~version on a unit circle,  with Lam\'e coefficients $\mu = \lambda = 1$, and data taken in such a way that the continuous solution is $\vec u = (u_1,u_2)$ with $u_1(x,y) = \sin(\pi x y) \sin(\pi(x-y))$ and $u_2(x,y) = \cos(\pi(x+1))\sin(\pi x^2)$. The results reported in Figure \ref{fig:elasticity} show a behaviour similar to the one observed for the Poisson equation.

%
%
%
%
%

	\begin{table}
	\centering
	{\footnotesize\begin{tabular}{c|cccc||cccc} 
			& \multicolumn{4}{c}{$\widehat \tau = 1/8$} & \multicolumn{4}{c}{$\widehat \tau = 1/16$} \\ 
			$H$		&	0.6187 &     0.1907 &      0.0489 &      0.0127&	0.6629   &  0.1864    &  0.0485 5 &	    0.0129 \\ \hline
			$k=1$ &	 0.0113 & 0.0083 & 0.0171 & 0.0516&	 0.0060 & 0.0115  & 0.0271  & 0.0723\\ 
			$k=2$ &  0.0173  & 0.0149  & 0.0564 & 0.1600 & 0.0099 & 0.0254 & 0.0811  & 0.2015\\ 
			$k=3$ & 	 0.0204  & 0.0236  & 0.1028 & 0.2088& 	0.0129& 0.0326 & 0.1180  & 0.2820\\ 
			$k=4$&  0.0204 & 0.0316  & 0.1460  & 0.2584& 	0.0152  & 0.0388  & 0.1683  & 0.3340\\ 
			$k=5$&	 0.0229 & 0.0343  & 0.1895 & 0.1724&	0.0197 & 0.0477  & 0.2005  & 0.2078\\
			$k=6$&	 0.0257 & 0.0422 & 0.2275 & 0.1377 &		0.0219 & 0.0490  & 0.2208  & 0.1858 
		\end{tabular}
	}
	\caption{Computational overhead resulting from the elimination of the ``lazy'' degrees of freedom: ratio between the time needed for the elimination, and the time needed for assembling the stiffness matrix for $\widehat \tau = 1/8$ and $\widehat\tau = 1/16$.}\label{tab:3}
\end{table}

\section{Conclusions and perspectives}
We evaluated the performance of  a boundary corrected virtual element method, for the numerical solution of the Poisson equation on curved smooth domain approximated by a polygonal domain of the type that can be easily built out of images (i.e. domain obtained as the union of pixels). The use of polygonal elements obtained as agglomeration of pixels allows boundary correction methods such as the shifted boundary method to satisfy the assumptions guaranteeing stability and optimality of the error estimates for arbitrary values of the order of the method. Eliminating, by a cheap static condensation procedure, a large number of degrees of freedom that do not actively contribute to the consistency of the method, allows to retrieve a robust behaviour of the error as a function of the number of degrees of freedom, independently of the order of approximation, that can be arbitrarily high. The numerical results demonstrate the potential of the method. {  While the paper only deals with the two dimensional case, its stability and convergence properties carry over to three dimensions, as shown in the preliminary tests carried out in \cite[Section 5.3]{VEM_weakly}, for $k=1,2$. This suggests that the method is well suited to be applied also in the three dimensional case, though in such a case, in the elimination of the ``lazy'' degrees of freedom extra care will be needed to deal with the ones lying on the wirebasket.} The method is well suited to be eventually coupled with an adaptive reconstruction of the smooth continuous boundary from imaging data, which, together with the extension to three dimensions, will be the focus of a forthcoming paper.

\appendix
\section{Proof of Lemma \ref{lem:trace}}\label{appendix:A}
We let $\interval{a}{b}$, $a<b$, denote the interval of extrema $a$ and $b$.
We start by observing that, under our assumptions,  the elements $K$ in our mesh are extension domains (\cite{Jones81}). We can then embed them in a square of diameter $\simeq \h$, which, for simplicity, we will assume to be the square $\widehat S = \interval0\h^2$, and, for all $v \in H^1(K)$ there exists $\widehat v \in H^1(\widehat S)$ with $\| \widehat v \|_{1,\widehat S} \lesssim \| v \|_{1,K}$, the constant in the inequality independent of $\h$ and $\mh$. Moreover, a Poincar\'e inequality holds in $K$, of the form
\[
\inf_{\alpha \in \mathbb{R}}\| v - \alpha \|_{0,K} \lesssim \h  | v |_{1,K}
\]
 with constants independent of $\h$ and $\mh$.
For $v \in C^1(\bar K)$ we can write
\[
\| v \|_{0,\bOh}^2 = \sum_{e\in \Ekh } \| v \|_{0,e}^2 + \sum_{e\in \Ekv } \| v \|_{0,e}^2,
\]
where $\Ekh  \subset \EK$ and $\Ekv  \subset \EK$ are the sets of, respectively, horizontal and vertical edges of $K$.
Let us consider the contribution of vertical edges. We have, with $e = \{x_e \} \otimes \interval{y_e^0}{y_e^1}$,

\begin{multline*}
\sum_{e\in \Ekv } \| v \|_{0,e}^2  = \sum_{e\in \Ekv } \int_{y^0_e}^{y_e^1} | v(x_e,y) |^2\,dy \lesssim 
\sum_{e\in \Ekv } \int_{y^0_e}^{y_e^1} | |K|^{-1} \int_K  v(\sigma,\tau)\,d\sigma\,d\tau |^2\,dy \\
+ 
\sum_{e\in \Ekv } \int_{y^0_e}^{y_e^1} | v(x_e,y) - |K|^{-1} \int_K  v(\sigma,\tau)\,d\sigma\,d\tau |^2\,dy = I + II.
\end{multline*}
We can write
\[
I\lesssim \sum_{e\in \Ekv } \int_{y^0_e}^{y_e^1}  |K|^{-1} \int_K  | v(\sigma,\tau) |^2\,d\sigma\,d\tau\,dy \leq \frac{| \partial K|}{|K|} \| v \|_{0,K}^2 \lesssim \h ^{-1} \| v \|^2_{0,K}.
\]
Now we have 
\begin{multline*}
II = \sum_{e\in \Ekv } \int_{y^0_e}^{y_e^1} |  |K|^{-1} \int_K  (v(x_e,y) - v(\sigma,\tau))\,d\sigma\,d\tau |^2\,dy\\ \leq 
\sum_{e\in \Ekv } \int_{y^0_e}^{y_e^1} |  |K|^{-1} \int_K  (v(x_e,y) - \widehat v(\sigma,y))\,d\sigma\,d\tau |^2\,dy \\+ \sum_{e\in \Ekv } \int_{y^0_e}^{y_e^1} |  |K|^{-1} \int_K  (\widehat v(\sigma,y) - v(\sigma,\tau))\,d\sigma\,d\tau |^2\,dy = III + IV.
\end{multline*}
We can bound $III$ as
\begin{multline*}
	III\lesssim \sum_{e\in \Ekv } \int_{y^0_e}^{y_e^1}   |K|^{-1} \int_K | (v(x_e,y) - \widehat v(\sigma,y))|^2\,d\sigma\,d\tau \,dy \\
	=
	 \sum_{e\in \Ekv } \int_{y^0_e}^{y_e^1}   |K|^{-1} \int_K | \int_{\sigma}^{x_e}  \partial_x \widehat v(\xi,y) \,d\xi |^2\,d\sigma\,d\tau \,dy \\
	 \lesssim \sum_{e\in \Ekv } \int_{y^0_e}^{y_e^1}   |K|^{-1} \int_K | x_e - \sigma | \int_{\sigma}^{x_e} | \partial_x \widehat v(\xi,y)|^2 \,d\xi \,d\sigma\,d\tau \,dy \\
	 \leq | K |^{-1} \int_K \h \sum_{e\in \Ekv } \int_{y^0_e}^{y_e^1}  \int_{0}^{\h} | \partial_x \widehat v(\xi,y)|^2 \,d\xi  \,dy \,d\sigma\,d\tau \\ \lesssim \h \sum_{e\in \Ekv } \int_{y^0_e}^{y_e^1}  \int_{0}^{\h} | \partial_x \widehat v(\xi,y)|^2 \,d\xi  \,dy \leq
	    H \| \partial_x \widehat v \|_{0,\widehat S}^2,
\end{multline*}
while $IV$ is bound as
\begin{multline*}
IV\lesssim \sum_{e\in \Ekv } \int_{y^0_e}^{y_e^1}   |K|^{-1} \int_K | (v(\sigma,y) - \widehat v(\sigma,\tau))|^2\,d\sigma\,d\tau \,dy \\
=
\sum_{e\in \Ekv } \int_{y^0_e}^{y_e^1}   |K|^{-1} \int_K | \int_{\tau}^{y}  \partial_y \widehat v(\sigma,\zeta) \,d\zeta |^2\,d\sigma\,d\tau \,dy \\
\lesssim \sum_{e\in \Ekv } \int_{y^0_e}^{y_e^1}   |K|^{-1} \int_{\widehat S} | y - \tau | \int_{\tau}^{y} | \partial_y \widehat v(\sigma,\zeta)|^2 \,d\zeta \,d\sigma\,d\tau \,dy \\
\lesssim \sum_{e\in \Ekv } \int_{y^0_e}^{y_e^1}  \h  |K|^{-1} \int_0^\h \int_0^\h \int_0^\h | \partial_y \widehat v(\sigma,\zeta)|^2 \,d\zeta \,d\sigma\,d\tau \,dy \lesssim 
H \| \partial_y \widehat v \|_{0,\widehat S}^2,
\end{multline*}
finally yielding
\[
II \lesssim III + IV \lesssim \h | \widehat v |^2_{1,\widehat S} \lesssim \h | v |^2_{1,K}.
\]
The contribution of horizontal edges is bound by the same argument, allowing to conclude that \eqref{trace} holds for $v$ smooth. 
The result for a generic $v \in H^1(K)$ is obtained by density. The bound \eqref{tracep} for polynomials is a direct  consequence of the combination of the previous bound with the inverse inequality
\[
\| p \|_{1,S(\xK,\h_K)} \lesssim \h ^{-1} \| p \|_{0,S(\xK, \alpha_1\h)}.
\]

\

In order to prove \eqref{trace-r}, we rely on the triangulation  $\TriK$ provided by Assumption \ref{additional}. For each edge $e \in \EK$ we let $T_e \in \TriK$ denote the triangle having $e$ as an edge. Thanks to the shape regularity of the triangulation we can write
\(
| v |_{r-1/2,e}  \lesssim | v |_{r,T_e},
\)
the implicit constant in the inequality only depending on the constant $\alpha_1$.
Then
\[
\sum_{e \in \EK} | v |_{r-1/2,e}^2 \lesssim \sum_{e \in \EK} | v |^2_{r,T_e} \lesssim \sum_{T \in \TriK} | v |^2_{r,T} \lesssim | v |_{r,K}^2.
\]

\section{Proof of Lemma \ref{lem:interpolation}}\label{appendix:C} 
Given $u \in H^s(K)$, $2 \leq s \leq k+1$, we aim at  constructing a quasi interpolant $\uI\in \VEMK$ such that we can prove an optimal estimate on $u - \uI$, robust in $h$ and $H$. We start by recalling that for all $v \in H^{1/2}(\bK)$ it holds that
\[
\inf_{{\phi \in H^1(K)}\atop{\phi=v\text{ on }\bK}} | \phi |_{1,K} = | \Harm{(v)} |_{1,K},
\]
where $\Harm(u)$ denotes the harmonic lifting of $u$. We can then consider, for $H^{1/2}(\bK)$ the non standard norm 
\[
\vvvert v \vvvert_{1/2,\bK} = | \Harm(v) |_{1,K} + \| v \|_{0,\bK}.
\]

\newcommand{\twI}{\widetilde w_I}

Letting  $\vvvert \cdot \vvvert_{-1/2,\bK}$ be defined as
\[
\vvvert \phi \vvvert_{-1/2,\bK} = \sup_{{v \in H^{1/2}(\bK)}\atop{\int_{\bK} v = 0}} \frac{\int \phi v}{\vvvert v \vvvert_{1/2,\bK}},
\]
we easily see that, for $v \in H^1(K)$ with $\Delta v \in L^2(K)$, $\Delta v$ average free, it holds, uniformly in $h$ and $H$, that
\begin{equation}\label{bounddnv}
\vvvert \nabla v \cdot \nK \vvvert_{-1/2,\bK}  \lesssim \| \Delta v \|_{(H^1(K)/\mathbb{R})'}  + | v |_{1,K}.
\end{equation}

\
 
For $u \in H^s(K)$, $2 \leq s \leq k+1$ we now let $\tuI \in \widetilde V^{K,\p}$ be defined as
\begin{gather*}
\tuI(x_i) = u(x_i) \text{ for all node $x_i$ of $\Tess$}, \ \text{ and }\ \int_K \Delta(u - \tuI) q
= 0 \text{ for all  $q \in\Poly{k}(K)$},
\end{gather*}
and we define $\uI \in \VEMK$ as
\begin{gather*}
	\uI = \tuI  \text{ on $\bK$}\qquad \text{ and } \qquad \int_K (\uI - \tuI) q = 0 \text{ for all $q \in \Poly{k-2}(K)$}.
\end{gather*}
We claim that $\uI$ thus defined satisfies \eqref{VEMapproxI}.  To prove our claim we now let $\tVK$ denotes the standard  finite element space of continuous piecewise polynomial functions of degree at most $k$ defined on the auxiliary triangulation $\TriK$ given by Assumption \ref{additional}. We let $\twI$ denote the interpolant of $u$ in $\tVK$. We have, by standard finite element approximation estimates
\begin{equation}\label{FEMapprox}
	| u - \twI |_{1,K} \lesssim H^{s-1} | u |_{s,K}.
\end{equation}
We then observe that we have
\[
\vvvert{u - \tuI} \vvvert_{1/2,\partial K} =
\vvvert u - \twI \vvvert_{1/2,\partial K} \lesssim | u - \twI |_{1,K} + \| u - \twI \|_{0,\partial K} \lesssim H^{s-1} | u |_{s,K}.
\]
Let us bound $| u - \tuI |_{1,K}$. Integrating by parts, using \eqref{bounddnv}
 as well as an Aubin Nitsche type duality trick to bound the $(H^1(K)/\mathbb{R})'$ norm of $\Delta (u - \tuI)$, plus some standard polynomial interpolation bound on $\partial K$, we can write
\begin{multline*}
	| u - \tuI |_{1,K}^2 = - \int_K (\Delta(u - \tuI)) (u - \tuI) + \int_{\bK} \nabla(u-\tuI)\cdot\nK (u - \tuI)\\[0.8mm]
	\lesssim \| \Delta(u - \tuI) \|_{(H^1(K)/\mathbb{R})'} | u - \tuI |_{1,K} + \vvvert \nabla(u-\tuI)\cdot\nK \vvvert_{-1/2,\bK} \vvvert u - \tuI \vvvert_{1/2,\bK}
\\[2mm]
	\lesssim \h ^{s-1} | u |_{s,K} | u - \tuI |_{1,K} + H^{2s - 2} | u |_{s,K}^2,
	\end{multline*}
yielding, for $\varepsilon > 0$ arbitrary and for some positive constants $C, C'$ only depending on the shape regularity parameters and on $s$,
\[ 
| \nabla (u - \tuI) |_{1,K}^2 \leq \frac C {2 \varepsilon} \h ^{2s-2} | u |_{s,K}^2 + \frac {C\varepsilon} {2} | u - \tuI |_{1,K}^2.
\]
Choosing $\varepsilon = 1/C$ yields the bound
\[
| u - \tuI |_{1,K} \lesssim \h^{s-1} | u |_{s,K}.
\]

{We now need to bound $| \tuI - \uI |_{1,K}$. Letting $\proj{\ell} : L^2(K) \to \Poly{\ell}(K)$ denote the orthogonal projection onto the space of polynomials of degree at most $\ell$ on $K$, integrating by parts and using the fact that both $\Delta \tuI$ and $\Delta \uI$ are polynomials of degree at most $k$ and that $\tuI - \uI$ is orthogonal to $\Poly{k-2}(K)$, we have
\begin{multline*}
	| \tuI - \uI |^2_{1,K} = - \int_{K} \underbrace{\Delta(\tuI - \uI)}_{\in \Poly{k}\cap\Poly{k-2}^\perp} (\tuI - \uI)  \\=
\int_{K} (\Delta (\tuI - \uI) - \underbrace{\proj{k-2}(\Delta (\tuI - \uI))}_{=0})(\tuI - \PinablaK \tuI) 
	\end{multline*}
where we could replace $\uI$ with $\PinablaK \tuI$ thanks to definition of $\uI$. Indeed, for $\uI \in \VEMK$ we have by definition that $\int_{K} \uI q = \int_K \PinablaK \tuI q$ for all polynomials $q \in \Poly{k} \cap \Poly{k-2}^\perp$,
where $\tuI \in \tVK$ is any element of $\tVK$ satisfying $\int_K \tuI p = \int_K \uI p$ for all $p \in \Poly{k-2}$ (see \eqref{defVEMK}).
Then, adding and subtracting $u$ at the second factor on the right hand side, we have
\begin{equation*}
		| \tuI - \uI |_{1,K} \lesssim 
	 \| (1 - \proj{k-2})\Delta(\tuI - \uI) \|_{(H^1(K)/\mathbb{R})'} (| \PinablaK \tuI - u|_{1,K} + | u - \tuI |_{1,K}).
	\end{equation*}
Using an Aubin-Nitsche's duality argument, a polynomial approximation bound and an inverse inequality for polynomials, we bound, for $\ell = k, k-2$
\begin{multline*}
\| (1 - \proj{\ell})\Delta(\tuI - \uI) \|_{(H^1(K)/\mathbb{R})'} \lesssim H \| (1 - \proj{\ell})\Delta(\tuI - \uI) \|_{0,K}\\[1.3mm] \lesssim H \| \Delta(\tuI - \uI) \|_{0,K} \lesssim | \tuI - \uI |_{1,K}.
\end{multline*}
Moreover, adding and subtracting $u$ and using a polynomial approximation bound, we have
\begin{multline*}
| \PinablaK \tuI - u |_{1,K} \leq | \PinablaK (\tuI - u) + (u - \PinablaK u) |_{1,K}\\
\leq | \tuI - u |_{1,K} + | u - \PinablaK u |_{1,K} \lesssim \h^{s-1} | u |_{s,K}.
\end{multline*}}
Combining the different bounds and using a triangular inequality yields the desired bound.

\section{Proof of Theorem \ref{cor:nitsche:curved}}\label{appendix:B}

From Lemma \ref{lem:trace},	we immediately obtain that, for all $\phi \in \VEMK$
\[
\| \dnh \PinablaK (\phi) \|_{0,\bK} \leq \|  \nabla \PinablaK (\phi) \|_{0,\bK} \lesssim \h^{-1/2} \| \nabla \PinablaK (\phi)\|_{0,K},
\]
as well as
\[
\|	\partial^j_\sigma  \PinablaK (\phi) \|_{0,\bK}  \lesssim \h^{1/2-j} \| \nabla \PinablaK (\phi)\|_{0,K}. 
\]


Then, it is not difficult to prove that the following bounds hold
\begin{gather}
	\| \Corrstar{\Pnphi} \|_{0,\partial\Oh} + \| \Corrhat{\Pnphi} \|_{0,\partial\Oh}  \leq \Cuno \frac h H \sqrt{H} | \Pnphi |_{1,\Th},\label{boundcorrections}\\
	\| \dnh \PinablaK (\phi) \|_{0,\bOh} \leq \Cdue H^{-1/2} | \Pnphi |_{1,\Th},\\[1mm]
	\| \phi - \Pnphi \|_{0,\bOh}  \leq \Ctre \h^{1/2} | \phi - \Pnphi |_{1,\Tess}
	\end{gather}
(we recall that the boundary correction operators $\Corrstar{\cdot}$ and $\Corrhat{\cdot}$ are  defined in \eqref{definecorrections}).


\

Let now $\Ah$ be defined as
\begin{multline}\label{defAH}
\Ah(\phi,\psi) = 	\ah (\phi,\psi)  
- \int_{\bOh} \dnh  \Pn(\phi)\psi
\\ - \int_\bOh  \left(
\Pn(\phi)  + \Corrstar{\Pnphi} \right) \left( \dnh \Pn(\psi)  - \gamma \h^{-1} \Corrhat{\Pnpsi} \right).
\end{multline}

\newcommand{\ext}[1]{\widehat{\mathscr{E}}[#1]}
\newcommand{\extstar}[1]{\mathscr{E}[#1]}
\newcommand{\diff}[1]{\widehat{\mathscr{D}}[#1]}

Continuity of the bilinear form $\Ah$ with respect to the norm $\vvvert \cdot \vvvert_{\Oh}$ follows from the above bounds. Let us prove that, provided $\mh/\h < \tau$ with  $\tau$ small enough, the bilinear form $\Ah$ is also coercive. Letting \[\ext{w} = w + \Corrhat{w} = \sum_{j=1}^\khat \frac{\delta^j} {j!} \dsh^j w, \quad \text{and}\quad \diff{w} = \Corrstar{w}-\Corrhat{w}= \sum_{j=\khat+1}^\kstar \frac{\delta^j} {j!} \dsh^j w,\] we can write, for $\epsilon > 0$ arbitrary
\begin{multline}
\label{eq:coecivity}
\Ah(\phi,\phi)  \geq | \Pn(\phi) |_{1,\Tess}^2 + \beta\cstar | \phi - \Pn(\phi) |_{1,\Tess}^2 + \gamma \h ^{-1} \| \ext{\Pn(\phi)} \|_{0,\bOh}^2  \\[1mm] - \int_{\bOh} \dnh \Pnphi (\phi - \Pnphi) -2\int_{\bOh} \dnh \Pnphi \ext{\Pnphi} - \int_{\bOh} \diff{\Pnphi} \dnh\Pnphi \\ + \int_{\bOh} \dnh \Pnphi \Corrhat{\Pnphi} + \gamma H^{-1} \int_{\bOh} \diff{\Pnphi}\ext{\Pnphi}
\\ \geq
\big(\frac 1 2 - {\Cdue^2}\epsilon - \frac \gamma 2 \Cuno^2 \tau^2 - (\Cuno^2 + \Cdue^2) \tau\big)| \Pn(\phi) |_{1,\Tess}^2 + \big(\beta\cstar  - \frac{\Ctre^2}{4\epsilon}\big)| \phi - \Pnphi |_{1,\Tess}^2 \\+ \big(
\frac \gamma 2 - 2{\Cdue^2} 
\big) H^{-1} \| \ext{\Pnphi} \|_{0,\bOh}^2.
\end{multline}
We now choose $\epsilon = 1/(4\Cdue^2)$ and let $\beta_0 =\Ctre^2/(4 \cstar\epsilon)$ and $\gamma_0 = 4\Cdue^2$, in such a way that for $\beta > \beta_0$ and $\gamma > \gamma_0$  we have $\beta\cstar - \Ctre^2/(4\epsilon) > 0$ and
$\gamma/2 - 2\Cdue^2 > 0$. For   $\gamma > \gamma_0$, let now $\tau_0(\gamma)$ denote the only positive solution of 
the equation $\frac 1 2 - \frac \gamma 2 \Cuno^2 \tau^2 - (\Cuno^2 + \Cdue^2) \tau = 0$. We easily see that for $\tau < \tau_0(\gamma)$, the coefficients of the first terms on the right hand side of \eqref{eq:coecivity} is strictly positive and the bilinear form $\Ah$ is, therefore, coercive with respect to the norm $\vvvert \cdot \vvvert_{\Oh}$. An unique discrete solution $\uh$ does thus exist.

\

Let $u_I$ denote the  VEM interpolant given by Lemma \ref{lem:interpolation} and $u_\pi \in \pwPoly$ the $L^2(\Omega)$ projection of $u$ onto the space of discontinuous piecewise polynomials, and set $d_h = u_I-u_h$. As in \cite{VEM_curvo}, we obtain 
\begin{equation}\label{B5}
\vvvert u_I - u_h \vvvert^2_{\Oh} \lesssim  |E1| + | E2 | + | E3 | + | E4 |  + | E5 | + | E6| + |E7|,
\end{equation}			
with 
\begin{gather*}
E1 = a_h(u_I-u_\pi,d_h), \qquad E2 = \sum_{K\in\Th} a^K(u_\pi - u,d_h),\\ E3 =  \int_{\bOh} \dnh(u - \Pn (u_I))\ext{\Pn(d_h)} ,\qquad 
E4 =  \int_{\bOh} \dnh(u - \Pn (u_I))\Corrhat{\Pn(d_h)} 
\\
E5 = \int_{\bOh} \dnh(u - \Pn (u_I))(d_h - \Pn(d_h)), \qquad 
E6 = \int_{\Oh}(f - \Pi^0 f) d_h,  \\
E7 = \gamma H^{-1} \int_\bOh (
 \tg - \extstar{\Pn(u)})(
 \gamma \h^{-1} \ext{\Pn(d_h)}
 - \dnh\Pn(d_h))
\end{gather*}

\

We observe that we have
\[
\| \dnh \Pn(d_h) - \gamma \h ^{-1} \Pn(d_h)) \|_{0,\bOh} \lesssim  \h ^{-1/2} \vvvert d_h \vvvert_{\Oh},
\]	
as well as
\[
\| d_h - \Pn(d_h)\|_{0,\bOh} \lesssim  \h^{1/2} | d_h - \Pn(d_h) |_{1,\Tess} \lesssim \h ^{1/2} \vvvert d_h \vvvert_{\Oh}.
\]
In view of these bound, of \eqref{boundcorrections}, and of the definition of the norm $\vvvert \cdot \vvvert_{\Oh}$,  all terms at the right hand side of \eqref{B5} are bound as in \cite{basicVEM} and \cite{VEM_curvo}, yielding
\[
\vvvert u_I - u_h \vvvert^2_{\Oh}  \lesssim \Big( \h ^{k} | u |_{k+1,\Omega} + \h^ {-1/2 } \mh^{\kstar+1} \| u \|_{\kstar+1,\infty,\Omega} \Big) \vvvert \uI - \uh \vvvert_{\Oh}.
\]
We obtain the desired bound by dividing both sides by $ \vvvert \uI - \uh \vvvert_{\Oh}$.

\bibliographystyle{elsarticle-num} 
\bibliography{biblio}

\begin{thebibliography}{10}
\expandafter\ifx\csname url\endcsname\relax
  \def\url#1{\texttt{#1}}\fi
\expandafter\ifx\csname urlprefix\endcsname\relax\def\urlprefix{URL }\fi
\expandafter\ifx\csname href\endcsname\relax
  \def\href#1#2{#2} \def\path#1{#1}\fi

\bibitem{SBMho}
N.~M. Atallah, C.~Canuto, G.~Scovazzi, The high-order shifted boundary method
  and its analysis, Computer Methods in Applied Mechanics and Engineering 394
  (2022) 114885.

\bibitem{VEM_weakly}
S.~Bertoluzza, M.~Pennacchio, D.~Prada, Weakly imposed {D}irichlet boundary
  conditions for 2{D} and 3{D} virtual elements, Computer Methods in Applied
  Mechanics and Engineering 400 (2022) 115454.

\bibitem{strang1973change}
G.~Strang, A.~E. Berger, The change in solution due to change in domain, in:
  Partial differential equations, 1973, pp. 199--205.

\bibitem{ramiere2008convergence}
I.~Rami{\`e}re, Convergence analysis of the q1-finite element method for
  elliptic problems with non-boundary-fitted meshes, International Journal for
  Numerical Methods in Engineering 75~(9) (2008) 1007--1052.

\bibitem{microFEM1}
N.~K. Knowles, J.~Kusins, M.~P. Columbus, G.~S. Athwal, L.~M. Ferreira,
  Experimental {DVC} validation of heterogeneous micro finite element models
  applied to subchondral trabecular bone of the humeral head~(9) (2022).
\newblock \href {https://doi.org/10.1002/jor.25229}
  {\path{doi:10.1002/jor.25229}}.

\bibitem{EdgeDetectionSurvey22}
R.~Sun, T.~Lei, Q.~Chen, Z.~Wang, X.~Du, W.~Zhao, A.~K. Nandi, Survey of image
  edge detection, Frontiers in Signal Processing 2 (2022).

\bibitem{parvizian2007finite}
J.~Parvizian, A.~D{\"u}ster, E.~Rank, Finite cell method, Computational
  Mechanics 41~(1) (2007) 121--133.

\bibitem{burman2015cutfem}
E.~Burman, S.~Claus, P.~Hansbo, M.~G. Larson, A.~Massing, Cutfem: discretizing
  geometry and partial differential equations, International Journal for
  Numerical Methods in Engineering 104~(7) (2015) 472--501.

\bibitem{Schillinger}
D.~Schillinger, M.~Ruess, The finite cell method: A review in the context of
  higher-order structural analysis of cad and image-based geometric models,
  Archives of Computational Methods in Engineering 22~(3) (2015) 391--455.

\bibitem{BDT}
J.~H. Bramble, T.~Dupont, V.~Thom{\'e}e, Projection methods for {D}irichlet's
  problem in approximating polygonal domains with boundary-value corrections,
  Math. Comp. 26~(120) (1972) 869--879.

\bibitem{SBMmechanics}
N.~M. Atallah, C.~Canuto, G.~Scovazzi, The shifted boundary method for solid
  mechanics, International Journal for Numerical Methods in Engineering
  122~(20) (2021) 5935--5970.

\bibitem{SBMfreesurface}
O.~Colom{\'e}s, A.~Main, L.~Nouveau, G.~Scovazzi, A weighted shifted boundary
  method for free surface flow problems, Journal of Computational Physics 424
  (2021) 109837.

\bibitem{SBMhyperbolic}
T.~Song, A.~Main, G.~Scovazzi, M.~Ricchiuto, The shifted boundary method for
  hyperbolic systems: Embedded domain computations of linear waves and shallow
  water flows, Journal of Computational Physics 369 (2018) 45--79.

\bibitem{SBMreduced}
E.~N. Karatzas, G.~Stabile, L.~Nouveau, G.~Scovazzi, G.~Rozza, A reduced-order
  shifted boundary method for parametrized incompressible navier--stokes
  equations, Computer Methods in Applied Mechanics and Engineering 370 (2020)
  113273.

\bibitem{SBM1}
A.~Main, G.~Scovazzi, The shifted boundary method for embedded domain
  computations. part {I}: Poisson and stokes problems, Journal of Computational
  Physics 372 (2018) 972--995.

\bibitem{SBM2}
N.~M. Atallah, C.~Canuto, G.~Scovazzi, The second-generation shifted boundary
  method and its numerical analysis, Computer Methods in Applied Mechanics and
  Engineering 372 (2020) 113341.

\bibitem{burman2020dirichlet}
E.~Burman, P.~Hansbo, M.~Larson, Dirichlet boundary value correction using
  lagrange multipliers, BIT Numerical Mathematics 60~(1) (2020) 235--260.

\bibitem{babuska}
I.~Babu\v{s}ka, The finite element method with {L}agrangian multipliers,
  Numer.Math. 20 (1973) 179--192.

\bibitem{stenberg1995some}
R.~Stenberg, On some techniques for approximating boundary conditions in the
  finite element method, Journal of Computational and applied Mathematics
  63~(1-3) (1995) 139--148.

\bibitem{burman2018cut}
E.~Burman, P.~Hansbo, M.~Larson, A cut finite element method with boundary
  value correction, Mathematics of Computation 87~(310) (2018) 633--657.

\bibitem{cheung2019optimally}
J.~Cheung, M.~Perego, P.~Bochev, M.~Gunzburger, Optimally accurate higher-order
  finite element methods for polytopial approximations of domains with smooth
  boundaries, Mathematics of Computation 88~(319) (2019) 2187--2219.

\bibitem{cockburn2012solving}
B.~Cockburn, M.~Solano, Solving dirichlet boundary-value problems on curved
  domains by extensions from subdomains, SIAM Journal on Scientific Computing
  34~(1) (2012) A497--A519.

\bibitem{liu2022weak}
Y.~Liu, W.~Chen, Y.~Wang, A weak galerkin mixed finite element method for
  second order elliptic equations on 2d curved domains, arXiv preprint
  arXiv:2204.01067 (2022).

\bibitem{VEM_curvo}
S.~Bertoluzza, M.~Pennacchio, D.~Prada, High order {VEM} on curved domains.,
  Atti Accad. Naz. Lincei Rend. Lincei Mat. Appl. 30 (2019) 391--412.

\bibitem{basicVEM}
L.~{Beir\~ao}~da Veiga, F.~Brezzi, A.~Cangiani, G.~Manzini, L.~D. Marini,
  A.~Russo, Basic principles of virtual element methods, Math. Models Methods
  Appl. Sci. 23~(1) (2013) 199--214.

\bibitem{hitchVEM}
L.~{Beir\~ao}~da Veiga, F.~Brezzi, L.~D. Marini, A.~Russo, The {hitchhiker}'s
  guide to the virtual element method, Math. Models Methods Appl. Sci. 24~(8)
  (2014) 1541--1573.

\bibitem{3DVEM}
B.~Ahmad, A.~Alsaedi, F.~Brezzi, L.~D. Marini, A.~Russo, Equivalent projectors
  for virtual element methods, {Comput. Math. Appl.} 66~(3) (2013) 376--391.

\bibitem{3DVEM2}
L.~{Beir\~ao}~da Veiga, F.~Dassi, A.~Russo, High-order virtual element method
  on polyhedral meshes, {Comput. Math. Appl.} 74 (2017) 1110--1122.

\bibitem{VEM_second_order}
L.~{Beir\~ao}~da Veiga, F.~Brezzi, L.~D. Marini, A.~Russo, Virtual element
  method for general second-order elliptic problems on polygonal meshes, Math.
  Models Methods Appl. Sci. 26~(4) (2016) 729--750.

\bibitem{beirao_parab}
G.~Vacca, L.~{Beir\~ao}~da Veiga, Virtual element methods for parabolic
  problems on polygonal meshes, Numer. Methods Partial Differential Equations
  31~(6) (2015) 2110--2134.

\bibitem{navier-stokes2d}
L.~{Beir\~ao}~da Veiga, C.~Lovadina, G.~Vacca, Virtual elements for the
  {Navier--Stokes} problem on polygonal meshes, SIAM Journal on Numerical
  Analysis 56~(3) (2018) 1210--1242.

\bibitem{stokes3d}
L.~{Beir\~ao}~da Veiga, C.~Lovadina, G.~Vacca, Divergence free virtual elements
  for the stokes problem on polygonal meshes, ESAIM: M2AN 51~(2) (2017)
  509--535.

\bibitem{Antonietti_VEM_Stokes}
P.~F. Antonietti, L.~{Beir\~ao}~da Veiga, D.~Mora, M.~Verani, A stream virtual
  element formulation of the {Stokes} problem on polygonal meshes, SIAM J.
  Numer. Anal. 52~(1) (2014) 386--404.

\bibitem{beirao_elastic}
L.~{Beir\~ao}~da Veiga, C.~Lovadina, D.~Mora, A virtual element method for
  elastic and inelastic problems on polytope meshes, Comput. Methods Appl.
  Mech. Engrg. 295 (2015) 327 -- 346.

\bibitem{beirao_linear_elasticity}
L.~{Beir\~ao}~da Veiga, F.~Brezzi, L.~D. Marini, Virtual elements for linear
  elasticity problems, SIAM J. Numer. Anal. 51~(2) (2013) 794--812.

\bibitem{VEM_3D_elasticity}
A.~L. Gain, C.~Talischi, G.~H. Paulino, On the virtual element method for
  three-dimensional linear elasticity problems on arbitrary polyhedral meshes,
  Comput. Methods Appl. Mech. Engrg. 282 (2014) 132--160.

\bibitem{Antonietti_VEM_Cahn}
P.~F. Antonietti, L.~B. da~Veiga, S.~Scacchi, M.~Verani, A {$C^1$} virtual
  element method for the {Cahn-Hilliard} equation with polygonal meshes, SIAM
  J. Numer. Anal. 54~(1) (2016) 34--56.

\bibitem{perugia_Helmholtz}
I.~Perugia, P.~Pietra, A.~Russo, A plane wave virtual element method for the
  {Helmholtz} problem, ESAIM Math. Model. Numer. Anal. 50~(3) (2016) 783--808.

\bibitem{ABPV:minsurf}
P.~Antonietti, S.~Bertoluzza, D.~Prada, M.~Verani, The virtual element method
  for a minimal surface problem, Calcolo 57~(4) (2020) 39.

\bibitem{VEM_Laplace_Beltrami}
M.~{Frittelli}, I.~{Sgura}, Virtual element method for the {Laplace-Beltrami}
  equation on surfaces, ESAIM Math. Model. Numer. Anal. 52~(3) (2018) 965 --
  993.

\bibitem{Moraetal15}
D.~Mora, G.~Rivera, R.~Rodr{\'i}guez, A virtual element method for the
  {S}teklov eigenvalue problem, Math. Models Methods Appl. Sci. 25~(8) (2015)
  1421--1445.

\bibitem{wriggers2017efficient}
P.~Wriggers, B.~Reddy, W.~Rust, B.~Hudobivnik, Efficient virtual element
  formulations for compressible and incompressible finite deformations,
  Computational Mechanics 60~(2) (2017) 253--268.

\bibitem{brenner2021ac}
S.~Brenner, L.~Sung, Z.~Tan, A ${C}^1$ virtual element method for an elliptic
  distributed optimal control problem with pointwise state constraints,
  Mathematical Models and Methods in Applied Sciences (2021).

\bibitem{VEM_mixed}
L.~{Beir\~ao}~da Veiga, F.~Brezzi, L.~D. Marini, A.~Russo, Mixed virtual
  element methods for general second order elliptic problems on polygonal
  meshes, ESAIM Math. Model. Numer. Anal. 50~(3) (2016) 727--747.

\bibitem{de2016nonconforming}
B.~Ayuso~de Dios, K.~Lipnikov, G.~Manzini, The nonconforming virtual element
  method, ESAIM: Mathematical Modelling and Numerical Analysis 50~(3) (2016)
  879--904.

\bibitem{BMPP_VEM_nonconforming_stab}
S.~Bertoluzza, G.~Manzini, M.~Pennacchio, D.~Prada, Stabilization of the
  nonconforming virtual element method, Comput. Math. with Appl. 116 (2022)
  25--47.

\bibitem{MASCOTTO_Trefftz}
L.~Mascotto, I.~Perugia, A.~Pichler, A nonconforming {Trefftz} virtual element
  method for the {Helmholtz} problem: Numerical aspects, Computer Methods in
  Applied Mechanics and Engineering 347 (2019) 445--476.

\bibitem{curved_Trefftz}
A.~Anand, J.~S. Ovall, S.~E. Reynolds, S.~Wei\ss{}er, Trefftz finite elements
  on curvilinear polygons, SIAM Journal on Scientific Computing 42~(2) (2020)
  A1289--A1316.

\bibitem{VEM_discrete_fracture}
M.~F. Benedetto, S.~Berrone, S.~Scial\'o, A globally conforming method for
  solving flow in discrete fracture networks using the virtual element method,
  Finite Elem. Anal. Des. 109 (2016) 23 -- 36.

\bibitem{BEIRAODAVEIGA2015327}
L.~{Beir\~ao}~da Veiga, C.~Lovadina, D.~Mora, A virtual element method for
  elastic and inelastic problems on polytope meshes, Computer Methods in
  Applied Mechanics and Engineering 295 (2015) 327--346.

\bibitem{CHI2017148}
H.~Chi, L.~{Beir\~ao}~da Veiga, G.~H. Paulino, Some basic formulations of the
  virtual element method (vem) for finite deformations, Computer Methods in
  Applied Mechanics and Engineering 318 (2017) 148--192.

\bibitem{wriggers2016virtual}
P.~Wriggers, W.~Rust, B.~Reddy, A virtual element method for contact,
  Computational Mechanics 58~(6) (2016) 1039--1050.

\bibitem{BEIRAODAVEIGA2021}
L.~{Beir\~ao}~da Veiga, F.~Dassi, G.~Manzini, L.~Mascotto, Virtual elements for
  {Maxwell's} equations, Computers {\&} Mathematics with Applications (2021).

\bibitem{Nitsche}
J.~Nitsche, {\"U}ber ein variationsprinzip zur {L}\"osung von
  {D}irichlet-{P}roblemen bei {V}erwendung von {T}eilr\"aumen, die keinen
  {R}andbedingungen unterworfen sind, Abhandlungen aus dem Mathematischen
  Seminar der Universit\"at Hamburg 36 (1970) 9--15.

\bibitem{beirao_stab}
L.~{Beir\~ao}~da Veiga, C.~Lovadina, A.~Russo, Stability analysis for the
  virtual element method, Math. Models Methods Appl. Sci. 27~(13) (2017)
  2557--2594.

\bibitem{atallah2021analysis}
N.~Atallah, C.~Canuto, G.~Scovazzi, Analysis of the shifted boundary method for
  the poisson problem in domains with corners, Mathematics of Computation
  90~(331) (2021) 2041--2069.

\bibitem{VEM_curved_beirao}
L.~{Beir{\~a}o da Veiga}, A.~{Russo}, G.~{Vacca}, The virtual element method
  with curved edges, ESAIM Math. Model. Numer. Anal. 53~(2) (2019) 375--404.

\bibitem{VEM_curved_brezzi}
L.~Beir{\~a}o~da Veiga, F.~Brezzi, L.~Marini, A.~Russo, Polynomial preserving
  virtual elements with curved edges, Mathematical Models and Methods in
  Applied Sciences 30~(08) (2020) 1555--1590.

\bibitem{VEMcurvoext1}
F.~Dassi, A.~Fumagalli, D.~Losapio, S.~Scial{\`o}, A.~Scotti, G.~Vacca, The
  mixed virtual element method on curved edges in two dimensions, Computer
  Methods in Applied Mechanics and Engineering 386 (2021) 114098.

\bibitem{VEMcurvoapp1}
F.~Aldakheel, B.~Hudobivnik, E.~Edoardo~Artioli, L.~{Beir\~ao}~da Veiga,
  P.~Wriggers, Curvilinear virtual elements for contact mechanics, Computer
  Methods in Applied Mechanics and Engineering 372 (2020) 113394.

\bibitem{VEMcurvoapp2}
F.~Dassi, A.~Fumagalli, I.~Mazzieri, A.~Scotti, G.~Vacca, A virtual element
  method for the wave equation on curved edges in two dimensions, Journal of
  Scientific Computing 90~(1) (2022) 1--25.

\bibitem{VEMcurvoapp3}
E.~Artioli, L.~{Beir\~ao}~da Veiga, F.~Dassi, Curvilinear virtual elements for
  {2D} solid mechanics applications, Computer Methods in Applied Mechanics and
  Engineering 359 (2020) 112667.

\bibitem{guyan1965reduction}
R.~J. Guyan, Reduction of stiffness and mass matrices, AIAA journal 3~(2)
  (1965) 380--380.

\bibitem{da2016serendipity}
L.~B. Da~Veiga, F.~Brezzi, L.~D. Marini, A.~Russo, Serendipity nodal vem
  spaces, Computers \& Fluids 141 (2016) 2--12.

\bibitem{wriggers2020serendipity}
P.~Wriggers, B.~Hudobivnik, F.~Aldakheel, Serendipity virtual elements for
  general element shapes, Computer Methods in Applied Mechanics and Engineering
  (2020).

\bibitem{beirao2018serendipity}
L.~Beir{\~a}o Da~Veiga, F.~Brezzi, F.~Dassi, L.~D. Marini, A.~Russo,
  Serendipity virtual elements for general elliptic equations in three
  dimensions, Chinese Annals of Mathematics, Series B 39 (2018) 315--334.

\bibitem{chen2023stabilization}
A.~Chen, N.~Sukumar, Stabilization-free serendipity virtual element method for
  plane elasticity, Computer Methods in Applied Mechanics and Engineering 404
  (2023) 115784.

\bibitem{franke1979critical}
R.~Franke, A critical comparison of some methods for interpolation of scattered
  data, Tech. rep., Naval Postgraduate School Monterey CA (1979).

\bibitem{harari1992c}
I.~Harari, T.~J. Hughes, What are c and h?: Inequalities for the analysis and
  design of finite element methods, Computer methods in applied mechanics and
  engineering 97~(2) (1992) 157--192.

\bibitem{collins2023penalty}
J.~H. Collins, A.~Lozinski, G.~Scovazzi, A penalty-free shifted boundary method
  of arbitrary order, Computer Methods in Applied Mechanics and Engineering 417
  (2023) 116301.

\bibitem{beirao_hp_exponential}
L.~{Mascotto}, L.~{Beir{\~a}o da Veiga}, A.~{Chernov}, A.~{Russo}, Exponential
  convergence of the hp virtual element method with corner singularities,
  Numer. Math. (2018) 138--581.

\bibitem{Artioli2017}
E.~Artioli, L.~Beir\~ao~da Veiga, C.~Lovadina, E.~Sacco, Arbitrary order 2d
  virtual elements for polygonal meshes: part i, elastic problem, Computational
  Mechanics 60 (2017) 355--377.

\bibitem{Jones81}
P.~Jones, Quasiconformal mappings and extendability of functions in {S}obolev
  spaces, Acta Math. 147 (1981) 71--88.

\end{thebibliography}

%
%

\end{document}